\documentclass[11pt,english]{article}
\usepackage{lmodern}
\usepackage[T1]{fontenc}
\usepackage[latin9]{inputenc}

\usepackage{geometry}
\usepackage{pdflscape} 
\usepackage{rotating}  
\usepackage{mathrsfs,amsmath,amsthm,mathtools,esint,bm}
\usepackage{amssymb,latexsym,wasysym}
\usepackage{physics}
\usepackage{authblk}
\usepackage{xfrac} 
\usepackage{listings} 
\usepackage{color} 

\usepackage{graphicx}
\usepackage{graphics}
\graphicspath{{./figures/}}

\usepackage{subcaption}
\usepackage{xfp}
\makeatletter
\define@key{Gin}{sqrtofarea}{%
    \def\Gin@req@sizes{%
      \edef\Gin@scalex{\fpeval{#1/sqrt(\Gin@nat@height*\Gin@nat@width)}}%
      \let\Gin@scaley\Gin@exclamation
      \Gin@req@height\Gin@scalex\Gin@nat@height
      \Gin@req@width\Gin@scalex\Gin@nat@width
      }%
  \@tempswatrue}
\makeatother

\usepackage{afterpage}

\usepackage{titlesec}

\titlespacing*{\section} {0pt}{3ex plus 1ex minus .2ex}{1.5ex plus .2ex}
\titlespacing*{\subsection} {0pt}{2.5ex plus 1ex minus .2ex}{1.25ex plus .2ex}
\titlespacing*{\subsubsection}{0pt}{2.25ex plus 1ex minus .2ex}{1ex plus .2ex}
\titlespacing*{\paragraph} {0pt}{2.5ex plus 1ex minus .2ex}{1em}

\usepackage[labelfont=bf,tableposition=top,figureposition=bottom]{caption}

\DeclareCaptionLabelSeparator*{spaced}{\\[2ex]}
\usepackage{booktabs}
\usepackage{multirow}

\usepackage{rotating}

\usepackage{enumitem}
\usepackage{pdfpages}
\usepackage{comment}

\usepackage{hyperref}
\usepackage{url}
\usepackage{color}

\usepackage{float}
\usepackage{placeins}
\usepackage{lscape}
\usepackage[english]{babel}

\usepackage[lined,boxed,linesnumbered,ruled,vlined]{algorithm2e} 
\SetAlFnt{\smaller}
\usepackage{floatpag}

 \usepackage{algpseudocode}

\usepackage{lineno}

\usepackage[toc,page]{appendix}
\usepackage{titlesec}

\usepackage{pdfpages}
\usepackage{eso-pic}
\usepackage{everyshi}

\usepackage[noabbrev,capitalise]{cleveref}

\usepackage{mathpazo}

\makeatletter
\renewcommand\normalsize{%
\@setfontsize\normalsize{9}{12.5}\selectfont}   
\makeatother
\normalsize

\usepackage{setspace}

\newcommand{\mathleft}{\@fleqntrue\@mathmargin0pt}
\usepackage{makecell}
\usepackage{longtable}
\makeatletter
\newcommand{\mathcenter}{\@fleqnfalse}
\newcommand{\mathrigth}

\makeatother

\usepackage{etoolbox}

\makeatletter
\patchcmd{\maketitle}{\@fnsymbol}{\@alph}{}{}  
\makeatother

\usepackage{booktabs}

\newcommand{\N}{\mathbb{N}}

\def\minimize{\mathop{\rm minimize}}

\newcommand{\st}{\textnormal{s.t.}}

\def\argmin{\mathop{\rm arg\,min}}

\newcommand*{\noun}[1]{\textsc{#1}}

\providecommand{\keywords}[1]
{
  \small	
  \textbf{\textit{Keywords:}} #1
}

\algnewcommand{\algorithmicand}{\textbf{ and }}
\algnewcommand{\algorithmicor}{\textbf{ or }}
\algnewcommand{\OR}{\algorithmicor}
\algnewcommand{\AND}{\algorithmicand}
\algnewcommand{\var}{\texttt}

\newsavebox{\smlmats}% Box to store smallmatrix content
\savebox{\smlmats}{$\left(\begin{smallmatrix}\sigma^{2}&0\\0&\sigma^{2}\end{smallmatrix}\right)$}

  \providecommand{\claimname}{Claim}
  \providecommand{\lemmaname}{Lemma}
\providecommand{\theoremname}{Theorem}
  \providecommand{\definitionname}{Definition}

\theoremstyle{plain}

\theoremstyle{remark}

\theoremstyle{definition}

\theoremstyle{plain}

\theoremstyle{plain}
\newtheorem*{sublemma*}{Sub-lemma}
\theoremstyle{plain}

\DeclareMathAlphabet{\mathpzc}{OT1}{pzc}{m}{it}  %(for some strange looking math font)

\newcommand{\footcite}[1] {\cite{#1}\let\thefootnote\relax\footnote{\cite{#1} \bibentry{#1}.}}

\newcommand{\twofootcite}[2] {\cite{#1,#2}\let\thefootnote\relax\footnote{\cite{#1} \bibentry{#1}.}\let\thefootnote\relax\footnote{\cite{#2} \bibentry{#2}.}}
\makeatletter
\let\save@mathaccent\mathaccent
\newcommand*\if@single[3]{%
  \setbox0\hbox{${\mathaccent"0362{#1}}^H$}%
  \setbox2\hbox{${\mathaccent"0362{\kern0pt#1}}^H$}%
  \ifdim\ht0=\ht2 #3\else #2\fi
  }
\newcommand*\rel@kern[1]{\kern#1\dimexpr\macc@kerna}
\newcommand*\widebar[1]{\@ifnextchar^{{\wide@bar{#1}{0}}}{\wide@bar{#1}{1}}}
\newcommand*\wide@bar[2]{\if@single{#1}{\wide@bar@{#1}{#2}{1}}{\wide@bar@{#1}{#2}{2}}}
\newcommand*\wide@bar@[3]{%
  \begingroup
  \def\mathaccent##1##2{%
    \let\mathaccent\save@mathaccent
    \if#32 \let\macc@nucleus\first@char \fi
    \setbox\z@\hbox{$\macc@style{\macc@nucleus}_{}$}%
    \setbox\tw@\hbox{$\macc@style{\macc@nucleus}{}_{}$}%
    \dimen@\wd\tw@
    \advance\dimen@-\wd\z@
    \divide\dimen@ 3
    \@tempdima\wd\tw@
    \advance\@tempdima-\scriptspace
    \divide\@tempdima 10
    \advance\dimen@-\@tempdima
    \ifdim\dimen@>\z@ \dimen@0pt\fi
    \rel@kern{0.6}\kern-\dimen@
    \if#31
      \overline{\rel@kern{-0.6}\kern\dimen@\macc@nucleus\rel@kern{0.4}\kern\dimen@}%
      \advance\dimen@0.4\dimexpr\macc@kerna
      \let\final@kern#2%
      \ifdim\dimen@<\z@ \let\final@kern1\fi
      \if\final@kern1 \kern-\dimen@\fi
    \else
      \overline{\rel@kern{-0.6}\kern\dimen@#1}%
    \fi
  }%
  \macc@depth\@ne
  \let\math@bgroup\@empty \let\math@egroup\macc@set@skewchar
  \mathsurround\z@ \frozen@everymath{\mathgroup\macc@group\relax}%
  \macc@set@skewchar\relax
  \let\mathaccentV\macc@nested@a
  \if#31
    \macc@nested@a\relax111{#1}%
  \else
    \def\gobble@till@marker##1\endmarker{}%
    \futurelet\first@char\gobble@till@marker#1\endmarker
    \ifcat\noexpand\first@char A\else
      \def\first@char{}%
    \fi
    \macc@nested@a\relax111{\first@char}%
  \fi
  \endgroup
}
\makeatother
\definecolor{mygreen}{RGB}{28,172,0} % color values Red, Green, Blue
\definecolor{mylilas}{RGB}{170,55,241}

\usepackage{multicol}
\usepackage[table]{xcolor}
\definecolor{mygray}{rgb}{.651,  .651,  .651}

\title{Last Mile Delivery with Drones and Sharing Economy}

\author[1]{Mehdi Behroozi \footnote{Corresponding author: m.behroozi@neu.edu; Phone: (617) 373-2032; Address: 449 Snell Engineering Center, 360 Huntington Ave, Boston, MA 02115. M.B. gratefully acknowledges the support of Northeastern University's Tier-2 Grant as well as an NSF Planning Grant for this research.}}
\author[2]{Dinghao Ma}

\affil[1,2]{Department of Mechanical and Industrial Engineering, Northeastern University, Boston, MA, USA}
\affil[1,2]{E-mail addresses: m.behroozi@neu.edu; ma.di@northeastern.edu}

\date{}

\begin{document}

\maketitle

\begin{abstract}
We consider a combined system of regular delivery trucks and crowdsourced drones, available via a sharing economy platform, to provide a technology-assisted crowd-based last-mile delivery experience. We develop analytical models and methods for a system in which package delivery is performed by a big truck carrying many packages to a neighborhood or a town in a metropolitan area and then the packages are assigned to crowdsourced drone operators to deliver them to their final destinations. We develop several optimization models for various cases of the problem and use a combination of heuristic algorithms to solve this NP-hard problem. Finally, we present computational results for the models and algorithms, and conduct an exhaustive sensitivity analysis to check the influence of different parameters and assumptions on the final results. We also provide extensive managerial insights on the benefits of our proposed model if implemented in practice.
\end{abstract}

\keywords{Last Mile Delivery; Drone Delivery; Crowdsourcing; Sharing Economy; Decentralized Logistics}

\section{Introduction}
\label{sec:Introduction}

The number of deliveries and the revenue obtained from delivery operations have been growing continuously and rapidly during the last two decades, thanks to the exponential growth of e-commerce. However, the efficiency of delivery operations remains a big challenge. 
The last mile of delivery process has consistently been one of the most expensive (nearly or even more than 50\% of the total cost), least efficient, and most polluting part of the entire parcel delivery supply chain \cite{joerss2016parcel,gevaers2011characteristics}; the fact that Amazon Flex has been paying \$18-\$25 for Uber-like package delivery services \cite{AmazonFlex}, while they have not increased their hourly wages to \$15 up until just recently \cite{AmazonWages},  speaks to the expensiveness of the last-mile delivery operations. 

The expensiveness of last-mile delivery is due to several factors including the facts that it is a labor-intense operation, it is a scattered operation serving different individual customers at dispersed places, which often results in underutilized carrier capacity  and the least economy of scale in the whole delivery operation, and that such deliveries are usually very time-consuming because of road congestion, accessibility of the destinations, and most importantly unattended deliveries. The rapidly increasing importance of same day and same hour delivery in our lives will make this operation even more inefficient. Such deliveries are mostly used for low-value high-frequency products such as grocery items for which the shipping cost could quickly become disproportionate in the eyes of consumers. 

The advancement of technology can revolutionize the conventional delivery practices and boost the efficiency. Among these advancements are the recent efforts to adopt \textit{autonomous vehicles}, \textit{unmanned aerial vehicles} (UAVs), \textit{automated guided vehicles} (AGVs), and other \textit{droids} in package delivery operations. The integration of autonomous and semi-autonomous technologies into the last-mile delivery operations in a centralized or decentralized manner have the potential to remove or mitigate the long-lasting factors such as pooling and routing inefficiencies that have been contributing to the expensiveness of last-mile delivery. Utilizing drones, for example, can lead to significant savings in the delivery operation and have been explored for a few years \cite{AmazonTestingDrone}. According to \cite{JungSunghun2017AoAP}, Amazon's shipping costs in the period between 2016 and 2019 was around \$6,400M, \$8,025M, \$10,025M, \$12,375M, and Prime Air operations could save about \$1,044M, \$2,677M, \$5,147M, \$8,549M from these costs.

\begin{figure}[t]
   \begin{center}
      \includegraphics[width=0.7\textwidth]{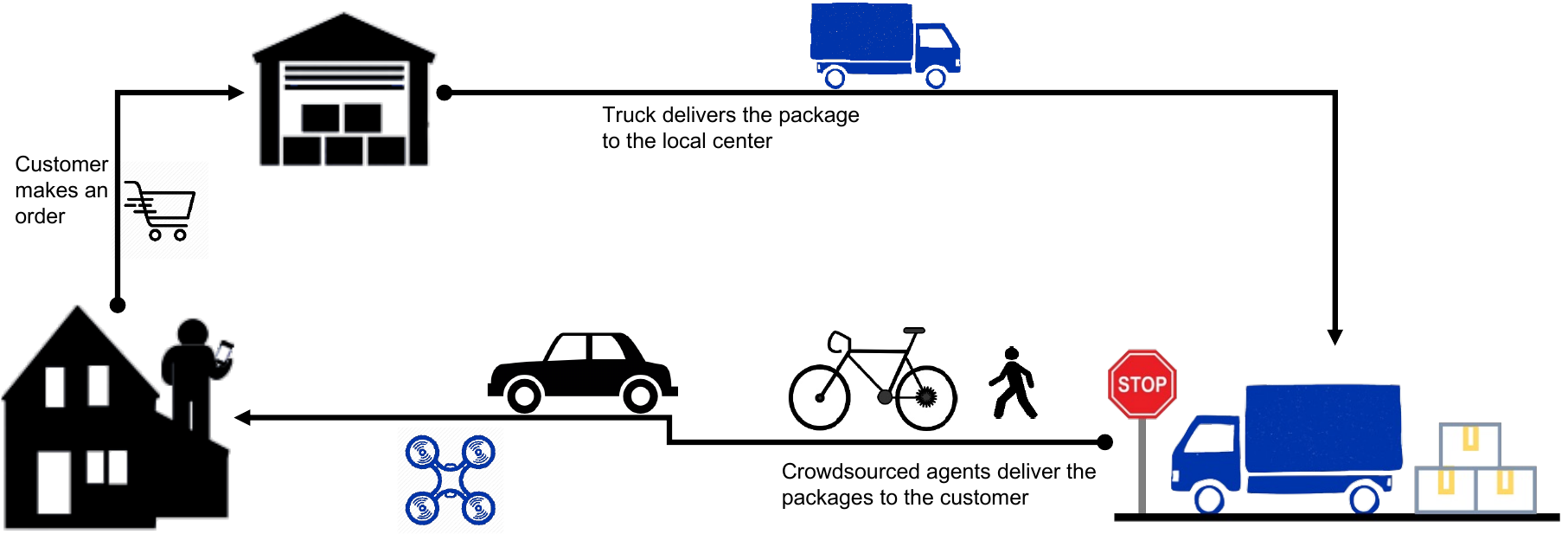}
   \caption{Sharing economy platforms, crowdsourced agents, and last-mile delivery.}
      \label{fig:selmd}
   \end{center}
\end{figure}

In an earlier work  \cite{carlsson2016household}, we have shown that for a centralized delivery system to be competitive with the decentralized household shopping model, a very large portion of the population have to adopt the centralized system and shows inefficient pooling as the primary cause of inefficient last-mile delivery. In this paper, we analyze the impact of decentralization, in particular crowdsourcing of the last part of the last-mile delivery operations when integrated with new technologies, on the efficiency of pooling and clustering customers.
To reduce the cost of the last mile delivery operation and reducing the delivery time (customer's waiting time), we analyze the traditional last mile delivery operation when conducted in coordination with crowdsourced vehicles, bikes, or drones, as illustrated in Figure \ref{fig:selmd}. These crowdsourced delivery operators can work in the context of sharing economy and form one side of a two-sided market platform. In fact, this provides a combination of centralized and decentralized systems. 

In this paper, we combine the crowdsourced autonomous delivery vessels with regular delivery trucks, vans, cars, and bikes 
to provide a technology-assisted last-mile delivery experience and a better and smoother transition to a fully autonomous parcel delivery ecosystem. We develop analytical models and methods for one of these intermediary systems, in which package delivery is performed by a big truck carrying a large number of packages to a neighborhood or a town in a metropolitan area and then assign the packages to crowdsourced delivery agents who operate \emph{drones}\footnote{The technology around civil applications of drones is developing rapidly, that can be observed by a quick look into the websites of major companies in this field. Many companies have developed different types of drones for video recording or delivering products. DJI is a famous manufacturer of drones for customers to record videos or take photos. Phantom 4 Pro V2.0 is one of its recent products which costs \$1599 and can travel 30 minutes with a full battery. Its speed is around 45 mph and weighs 1375 g. Amazon, Google, Domino, and DHL are all developing their own drones for delivering packages. Amazon's Prime Air can carry products with a weight of 2.26 kg and fly for 30 minutes with 50 mph speed, while it can be controlled 10 miles away from the base. One could compare this to the Google Wing that has a capacity of 1.5 kg and can fly at 75 mph within 8.7 miles from the base. DHL's commercial drones are the biggest ones with 2200 mm diagonal size and can fly at 43.5 mph within 7.5 miles from the base carrying packages up to 3 kg.} to deliver them to their final destinations. Other applications such as post disaster response operations in which trucks cannot access all customers and have to rely on drones in the neighborhood can also be considered for this problem. We introduced this problem in our preliminary work  \cite{behroozi2020crowdsourced}. In this paper, we expand that and study and analyze this problem in a more comprehensive way. To the best of our knowledge, besides our earlier extended abstract, this is the first work studying this problem. A combination of heuristic algorithms is used to solve this NP-hard problem and an exhaustive sensitivity analysis is done to check the influence of different parameters and assumptions such as speed ratio of drones and trucks, the number of drones in the service region, and the distribution of the customers. The simulation results show significant savings in the total delivery cost under reasonable assumptions.

\subsection{Related Work}
\label{subsec:Literature}
Relevant work in the literature are mainly in two categories: drone-aided delivery and crowdsourced delivery, with the latter being more relevant to our work. The literature of drone-aided delivery is so young but has been expanding rapidly. It can be divided into four main categories; travelling salesman problem (TSP) with drones, vehicle routing problem (VRP) with drones, Done(-only) delivery problem, and carrier-vehicle with drones \cite{macrina2020drone}. Murray and Chu \cite{murray2015flying} introduced routing problems for two combined systems of truck and drone, namely flying sidekick travelling salesman problem (FSTSP) and parallel drone scheduling TSP. The first problem and its variants have been called with several names including HorseFly Problem, TSPD, Truck-and-Drone TSP, etc. Since this first paper in 2015 the literature on such problems have been growing rapidly (see \cite{ponza2016optimization,Ferrandez2016optimization,carlsson2018coordinated,agatz2018optimization,bouman2018dynamic,poikonen2020mothership,salama2020joint,bruni2022logic} for example).
Although these problems resemble some similarities with the problem we study in this paper, but they are essentially solving a different problem. The main difference is that all of these problems are defined under a centralized delivery system, while our problem combines a centralized truck delivery system and decentralized drone delivery system, and to the best of our knowledge is the first to do so. We will show that the drone(-only) delivery problem is the centralized delivery version of a special case of our problem when the range and speed of drones are large enough, while the pure TSP model is the centralized alternative for the case where the flying range and speed of drones are rather small. For a review on the literature of centralized drone-aided delivery problems please see \cite{macrina2020drone,chung2020optimization,moshref2021applications,madani2022hybrid}.

Crowdsourced delivery, in essence, can be categorized under pickup and delivery problem (PDP), where the goal is to transport parcels from origins to destinations while minimizing costs. The PDP has been studied for almost half a century; see \cite{berbeglia2007static,cordeau2008recent} for overview on this body of literature. Unlike the traditional PDP setting, our problem has on-demand features, which is also somewhat related to the literature of dynamic PDP \cite{berbeglia2010dynamic}, and the delivery operation is partially conducted on a sharing economy platform. Sharing economy indicates a system in which people share access to goods and services as opposed to ownership and it has been extensively studied \cite{heydari2023reengineering,rojanakit2022sharing,hossain2020sharing,standing2019implications}. However, many aspects of the application of sharing economy systems, especially those with (semi-)autonomous technologies, in logistics and delivery services has received less attention \cite{boysen2021last,savelsbergh2022challenges,duman2023shared}. The paper \cite{sadilek2013crowdphysics} introduced  the idea of ``crowdphysics'' -- crowdsourcing tasks that require people to collaborate and synchronize both in time and physical space and discussed the concept using a crowd-powered delivery service. The paper \cite{rouges2014crowdsourcing} introduced a typology of the early startup businesses.
A first framework for crowd logistics was proposed in \cite{mladenow2015crowdsourcing}. The paper \cite{oliveira2019crowd} developed a more comprehensive framework for crowd-based operations in city-level logistics. The authors discussed different kinds of crowd-based logistics including crowdsourced delivery, cargo-hitching (i.e., integration of freight and passenger transport), receiving packages through pickup locations, providing storage by the crowd, and returning (reverse flow) logistics. There have been a number of theoretical and experimental research related to each of these topics, out of which we review those relevant to crowdsourced delivery.

In 2013 Wal-Mart started to consider the idea the idea of encouraging individuals/shoppers in a store who are willing to deliver packages for online customers on their way back home \cite{WalMartOccasionalCouriers}. Later, Amazon started considering the crowd as couriers \cite{AmazonNextDroneYou}. Many researchers explored the crowdsourced delivery service with crowdsourced drivers coming to the package center to pickup the packages and deliver them to their destination. Among the relevant theoretical research, the paper \cite{archetti2016vehicle} studied this idea and formulated it as a vehicle routing problem with occasional drivers (VRPOD). The paper \cite{dayarian2020crowdshipping} studied this problem more in-depth under both static and dynamic settings, in which the latter incorporates probabilistic information about future online orders and in-store customer arrivals. The paper \cite{gdowska2018stochastic} proposed a stochastic optimization model for the same problem. 
A similar problem with stochastic supply of crowd vehicle and time windows was analyzed in \cite{torres2022vehicle}. The authors presented a two-stage stochastic model, an exact branch and price algorithm, a column generation heuristic.
The paper \cite{arslan2019crowdsourced} also considered a group of ad-hoc drivers for performing parcel delivery tasks and proposed a route planning system that automatically makes the matching between ad-hoc drivers and these tasks. The platform also operates a fleet of dedicated vehicles for delivery tasks that cannot be completed by the ad-hoc drivers. The authors present a rolling horizon framework and an exact solution approach based on a matching formulation to solve the problem. They also compared their results with the traditional delivery system and concluded that the use of crowdsourced drivers can significantly reduce the costs.

On the experimental side, the paper \cite{suh2012leveraging} tried to leverage social networks to find available crowd couriers. They compared three delivery systems with respect to their greenhouse gas emission: the traditional delivery system, a pickup location delivery system (PLS) in which parcels are delivered to a pickup center and customers pick them up on their way to home, work, etc, and finally a socially networked PLS (SPLS) where the last-mile delivery from the pickup locations are done by agents socially connected to the customers. They showed that the PLS system is worse than the traditional delivery system but SPLS significantly reduces the emission.
The paper \cite{devari2017crowdsourcing} also studied the idea of creating a social network-supported last mile delivery and developed a simulation model for it, in which they exploit the social networks of the online customer to find potential couriers in their networks.
A trial with crowd participants for book delivery in a Finish library was conducted in \cite{paloheimo2016transport} that showed crowdsourced delivery reduces vehicle miles traveled. 
The paper \cite{miller2017crowdsourced} used a survey to analyze potential driver behavior in choosing to change status from pure commuters to traveler-shippers. Meanwhile, the paper \cite{punel2017modeling} also  created a survey to study the determinants of crowd-shipping acceptance among drivers. The paper \cite{chen2017crowdsourced} developed an agent-based simulation model for the crowdsourced last-mile delivery service with the existence of central pickup location/warehouse and identified the important factors influencing its performance. They set the packages to have three types of accepted delivery time windows and three kinds of crowdsourced agents, including driving, cycling, or walking, delivering one package each time. 
In contrast to an entirely crowdsourced fleet, the performance of a hybrid fleet of privately owned delivery vehicles and crowdsourced assets was simulated in \cite{castillo2022hybrid}.

Relatively fewer papers in the literature have considered a cooperative delivery system with cooperation between traditional delivery trucks  and crowdsourced carriers or between the crowdsources themselves. The paper \cite{qi2018shared} evaluated the use of shared mobility for last mile delivery services in coordination with delivery trucks. They tried to minimize the combined transportation and outsourcing cost of the trucks and shared mobility. They also considered minimizing greenhouse gas emissions as one of their objectives. 
They used an analytical model and found that crowdsourced shared mobility is not as economically scalable as the conventional truck-only system with respect to the operating costs, due to the fact that the payment to crowdsourced drivers is route specific and accounts for the competition with the ride-share service market. However, they state that a transition towards this model can create other economic and operational benefits such as reducing the truck fleet size, avoiding high-demand areas and peak hours, and adjusting vehicle loading capacities. 

Cooperation is usually conducted through intermediary depots. A crowd-delivery problem with a number of pop stations as intermediary between trucks and crowdsource workers was modeled as a network min-cost flow problem in \cite{wang2016towards} in which workers are assigned to their closest pop station.
The paper \cite{kafle2017design} considered cyclists and pedestrians as crowdsources who are close to customers and can perform the last-leg parcel delivery and the first-leg parcel pickup with relaying the parcels with a truck carrier at some relay points selected from a pre-defined set of such points. They formulated the problem as a mixed integer non-linear program to simultaneously select crowdsources and relay points and find the truck routes and schedule. The crowdsource routes and their feasibility are determined by the crowdsources themselves before making a bid (that includes the job they are willing to do and the price) and the platform only selects the the best offered bids. 
Pedestrians were also considered in \cite{akeb2018building} as crowd delivery partners assigned to their closest relay point determined by a circle packing scheme.
A similar problem was discussed in \cite{huang2019decision} and was formulated as a vehicle routing problem (VRP) type model, in which customers can be served by either an operator's truck or a crowdsourcing partner that received the package at a relay point.
The paper \cite{macrina2020crowd} also considered a VRPOD with time windows and intermediate depots. The authors also considered the scenario in which the crowdsourced drivers are forbidden at the central depot and can only visit the transshipment centers. Similar transshipment points with occasional drivers were considered in a pickup and delivery problem in \cite{voigt2022crowdsourced}. The same problem was considered in \cite{vincent2022crowd} assuming the existence of pickup locations where the packages could be dropped off and the customers can walk and pick up their packages there. A multi-objective version was analyzed in \cite{lan2022multi}. Another collaborative system between delivery truck and crowdsourced drivers was considered in a recent paper \cite{elsokkary2023crowdsourced} for which the authors present a heuristic algorithm that solves the involved routing and assignment problems separately. 

\newpage
The idea of relaying packages have also been used in a setting that crowd agents could cooperate with each other. The paper \cite{chen2018multi} analyzed a multi-driver multi-parcel matching problem (MDMPMP) based on the multi-hop ride sharing problem, in which parcels may be transported by a single or by multiple drivers, being relayed between drivers on the path from parcels origin to its destination. Crowd drivers have preset routes from their origins to their destinations from which can deviate slightly to pick up or drop off a package.

There is a lack of a study on the design of a cooperative delivery system with a truck and autonomous or semi-autonomous crowdsourced carriers. In this paper, we fill this gap.

\subsection{Regulations}
\label{subsec:Regulations}
Before we get into the details of the problem and our solution approach, it is necessary to discuss the applicability of the studied problem. One important area that needs to be considered in drone-related operations is their regulations. In 2016, the U.S. government announced that commercial use of drones needs to be registered by the Federal Aviation Administration(FAA) and Department of Transportation(DOT) and need to follow the new regulations \cite{Dnr}. The person who operates the drones needs either to be at least 16 years old and get a remote pilot certificate with a small unmanned aircraft systems (UAS) rating or to be directly supervised by someone who already has the certificate \cite{Pr2016}.

The regulations state that all drones including their carried materials or products must weigh less than 55 lbs which is 25 kg, the limit speed of drones is 100 mph, and the maximum altitude is 400 feet above ground level. The drones need to be within the visual line-of-sight of the remote pilot during the flight and need to be close enough to receive the command from the pilot. Meanwhile, all the drones must only be operated in daylight which is from 30 minutes before official sunrise to 30 minutes after official sunset. They need to yield the way to other aircraft and the minimum weather visibility is 3 miles from a control station. One person can only operate one drone at a time and cannot manipulate a moving vehicle unless it is a sparsely populated area. Moreover, drones cannot carry hazardous materials \cite{SUMAR}.

Other regulations are mainly related to the safety issues that drones need to be checked before their launching and the object being carried by the drone is securely attached. On the other hand, most of the restrictions can be waived if the applicant demonstrates that the operation can safely be conducted under the terms of a certificate of waiver \cite{SUMAR}.

\section{Problem Statement}
\label{sec:ProblemStatement}
\begin{figure}[t]
\centering
      \includegraphics[width=0.65\textwidth]{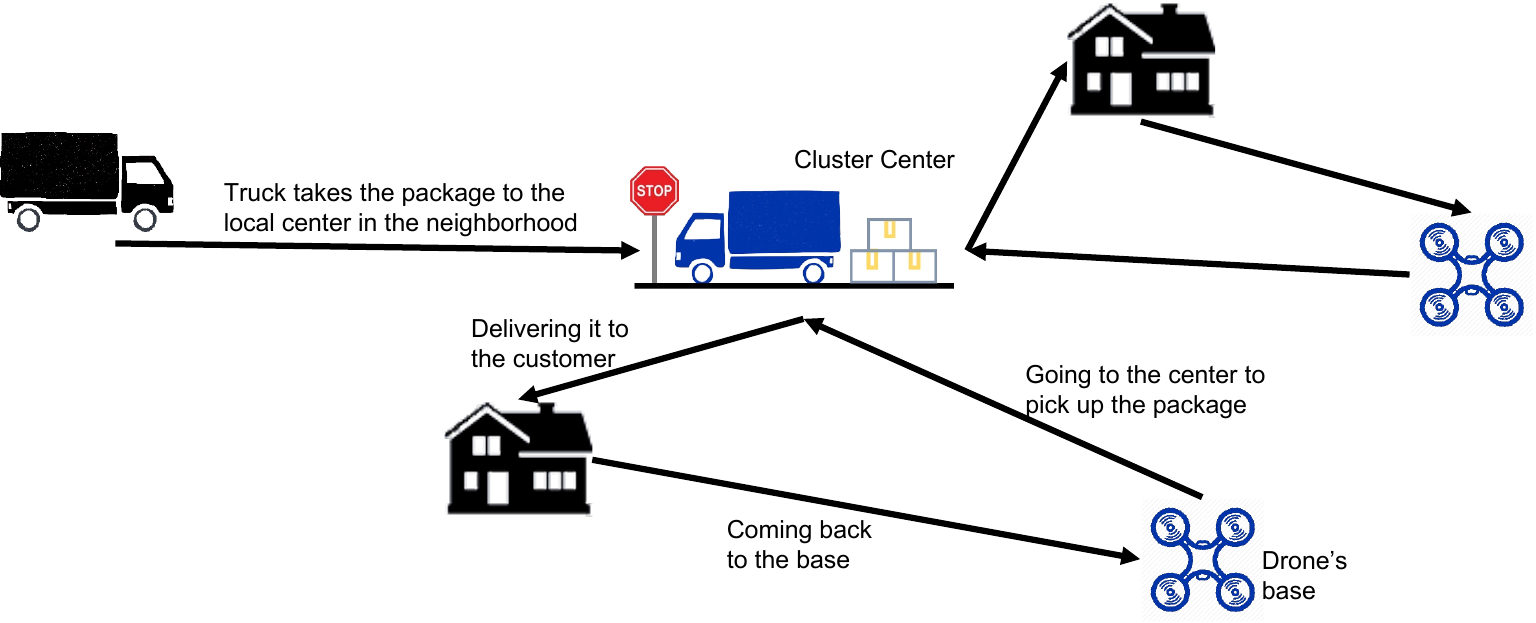}
      \protect\caption{ \label{fig:csdd}A schematic representation of the proposed crowdsourced drone delivery.}
\end{figure}

Consider a region $R$ as a residential area in which one truck has to go through all the neighborhoods to deliver some packages. There are also private drone operators in the area that could deliver packages from the truck to their final destination (households). When the truck stops at a neighborhood corner, the crowd-based drones, after receiving an order from the courier, will fly from their base to that corner to pick up the packages, deliver them to the customers and go back to their base, i.e., operator's house, for recharging its battery. The objective is to design a coordinated system between the truck and these crowdsourced drones in a way that minimizes the total time spent on fulfilling the demand of all customers in that area. 
 Figure \ref{fig:csdd} shows an illustration of this cooperative delivery between a truck and crowd-sourced drones. 

In this problem we assume that:
\begin{enumerate}
          \item Each drone can only carry one item at a time.
          \item The charging time for drones at its home base is 0 (i.e., changing quickly to a new battery).
          \item Each drone base location launches only one drone.
          \item The speed of all drones is set to be fixed and equal to each other. 
          \item There is no weight limit for a drone to carry the package. So all packages of a customer will be delivered in one trip.
          \item Only one truck is used and its capacity is sufficient for the entire demand in the region.
          \item If there is no drone nearby, the truck will serve all the customers in that neighborhood.
          \item Each drone can be only used at one truck's stop (cluster center).
          \item The truck will wait at the pickup center to track the drones until the driver receives the confirmation that the packages are delivered to the customers by the drone and send necessary instructions to drone operators in case of delivery failures.
          \item The compensation mechanism for drone delivery tasks is assumed to be set by the platform separately. 
\end{enumerate}

Depending on whether there are one or multiple pickup centers, or whether drones need to go back to their home bases after each delivery or not, we define four problems as followings:
\begin{description}
\item [Problem I--(One Center with Recharging):] Problem with one pickup center with the assumption that after each delivery, drones should stop at their home bases (for recharging/replacing battery) before revisiting the pickup center for another delivery.
\item [Problem II--(One Center with Revisiting):] Problem with one pickup center with the assumption that after each delivery, drones could revisit the pickup center for another delivery before going back to their home bases.
\item [Problem III--(Multiple Centers with Recharging):] Problem with multiple pickup centers with the assumption that after each delivery, drones should stop at their home base (for recharging/replacing battery) before revisiting the pickup centers for another delivery.
\item [Problem IV--(Multiple Centers with Revisiting):] Problem with multiple pickup centers with the assumption that after each delivery, drones could revisit the pickup centers for another delivery before going back to their home base.
\end{description}     
\begin{table}[!h]
\protect\caption{The four considered problems.}
\label{tab:4Problems}
\begin{center}
\begin{tabular}{|c|c|c|}  
      \hline
      & \textbf{Recharging} & \textbf{Revisiting} \\
      \hline
      \textbf{One Center} & Problem I & Problem II \\
      \hline
     \textbf{Multiple Centers} & Problem III & Problem IV \\
       \hline
   \end{tabular}
\end{center}
\vspace{-12pt}
\end{table}

\Cref{tab:4Problems} summarizes the conditions of these four problems. In the followings we will model these four problems. In the rest of this section, we delve into each of these problems and model them mathematically. The following notations are commonly used in all of the presented models:  we denote the set of customer nodes with $C$ and set of drone nodes with $D$ and let $v_D$ be the speed of drones, $v_T$ be the speed of trucks, $L$ be the longest distance (range) a drone can travel, and $M$ to be a sifficiently large number. Note that setting $L=0$ reduces problems III and IV to TSP.

\subsection{Problem I--(One Center with Recharging)} 	
\label{subsec:ProblemI}
In this problem we have one pickup center and assume that drones should go back to their home base for recharging after each delivery. There is no truck route and truck (or a store, warehouse, or distribution center) operates as a depot for a fleet of crowdsourced drones to pick up the package and deliver them to the customers. The problem with one center is important to study because it sets a a basis for the multiple centers problems and also it helps to understand the dynamics inside a cluster, introduced later, in a better way. The insight behind our algorithm is partly driven by this sub-problem. Before we present our model for this problem we define the sets, parameters, and variables in Table \ref{tab:Sets-Parameters-Variables-ProblemI}. 

\begin{table}[!h]
\protect\caption{Sets, parameters, and variables used in the one center model.}
\label{tab:Sets-Parameters-Variables-ProblemI}
\begin{center}
\begin{tabular}{c|p{0.72\textwidth}}  
      \hline
      \hline
      \textbf{Parameters/Variables} & \textbf{Description} \\
      \hline
      $d_{ij}$ & Route length going from node $i\in D$ to the center node (0) and then to node $j\in C$ and back to node $i$\\
      \hline
       $x_{ij}$ & Binary decision variable. It is $1$ when a drone travels from node $i\in D$ to the center node 0 and then to node $j\in C$ and back to node $i$. It is $0$ otherwise. \\
       $q_{i}$ & Total travel time of the drone with base at node $i$\\
       $Q$ & Maximum time spent by all drones (makespan) \\
       \hline
   \end{tabular}
\end{center}
\end{table}

 A mixed integer linear programming (MILP) model for this problem can be written as follows: 
  \begin{eqnarray}
      \minimize \qquad  Q & &  \quad \qquad \st \nonumber  \\
        \sum_{i\in D} x_{ij}  & = &  1\,, \qquad \forall j\in C 		\label{eq:Problem1SatisfyAll} \\
       d_{ij}x_{ij}  & \leq &  L\,,  \qquad \forall i\in D, \; \forall j\in C 		\label{eq:Problem1DroneRange} \\
      \frac{1}{v_D} \sum_{j \in C} d_{ij}x_{ij} & \leq & q_i \,, \qquad \forall i \in D 	\label{eq:Problem1DroneTime} \\
       0 \;\; \leq \;\; q_i  & \leq &  Q\,, \qquad \forall i\in D 						\label{eq:Problem1TotalTime} 	\\
       x_{ij} & \in & \{0, 1\}\,, \quad \forall i\in D, j\in C   \nonumber 
  \end{eqnarray}
The objective is to minimize the maximum time of drone routes (makespan). Constraint \eqref{eq:Problem1SatisfyAll} makes sure all customers are visited once. Constraint \eqref{eq:Problem1DroneRange} ensures that the distance traveled by each drone does not exceed the maximum travel distance (range) of drones. Constraints \eqref{eq:Problem1DroneTime} finds the time spent by each drone and Constraint \eqref{eq:Problem1TotalTime} calculates the maximum time among all drones. An alternative and less efficient but more flexible (adjustable) model for this problem is provided in \cref{sec:AltModelProblemI} of the Online Supplement of this paper.

\subsection{Problem II--(One Center with Revisiting)} 	
\label{subsec:ProblemII}
In this problem we consider one pickup center (no truck route) and after each delivery, drones could revisit the pickup center for another delivery before going back to their home bases. 
We first note that a feasible drone route in this case is a route that starts from a drone node $i\in D$, satisfies the demand of one or multiple customer nodes after picking their packages up at the center node, and comes back to the starting drone node while not violating the drone range $L$.
Let $\Pi_i$ denote the set of feasible drone routes starting from drone node $i\in D$ and let $\Pi_i(j) \subset \Pi_i$ be the set of routes that contain the customer node $j\in C$. We use $k$ to index these routes. For any route $\pi_{ik} \subset \Pi_i$ let $d_{\pi_{ik}}$ denote its total length. We define our decision variable as
\begin{flalign}
    x_{\pi_{ik}}  & =  \left\{ \begin{array}{ll}
                1 & \text{if route } \pi_{ik} \in \Pi_i \text{ is selected};\\
                0 & \mbox{otherwise}.\end{array} \right. && \nonumber
\end{flalign}

Then we can form an MILP formulation:
\begin{eqnarray}
 \minimize \qquad  Q \qquad & &  \qquad \qquad \st \nonumber  \\
 	& & \nonumber \\
    \sum_{i} \sum_{\pi_{ik} \in \Pi_i(j)}  x_{\pi_{ik}} & =  & 1\,, \quad \qquad \forall j\in C   \label{eq:Problem2SatisfyAll} \\
    d_{\pi_{ik}} x_{\pi_{ik}} & \leq & L  \;\; \quad \qquad \forall i\in D,\; \pi_{ik} \in \Pi_i  \label{eq:Problem2DroneRange} \\
   \frac{1}{v_D} \sum_{\pi_{ik} \in \Pi_i} d_{\pi_{ik}} x_{\pi_{ik}} & \leq & q_i\,, \quad \qquad \forall i \in D \label{eq:Problem2DroneTime} \\
   0 \;\; \leq \;\; q_i & \leq & Q\,, \quad \qquad \forall i\in D \label{eq:Problem2TotalTime} \\
    x_{\pi_{ik}} & \in & \{0, 1\}\,, \;\;\quad \forall i\in D, \; \pi_{ik} \in \Pi_i \nonumber %\\
\end{eqnarray}
Constraint \eqref{eq:Problem2SatisfyAll} ensures all customers are visited once. 
Constraints \eqref{eq:Problem2DroneTime} and \eqref{eq:Problem2TotalTime} find the maximum time among all drone routes. 
An alternative and less efficient but more flexible (adjustable) model for this problem is provided in \cref{sec:AltModelProblemII} of the Online Supplement of this paper.

\subsection{Problem III--(Multiple Centers with Recharging)} 	
\label{subsec:ProblemIII}
In this problem we can have multiple pickup centers and drones should go back to their home base for recharging after each delivery. All the customer and drone nodes are in a 
 region $R$. We first assume that truck nodes or centers can be placed anywhere in $R$ (we will change this later). This implies that all deliveries are done by the drones and the truck carries package through different neighborhoods and stops at one center in each neighborhood where it passes the packages to the drones. 
 Although the pickup centers are not necessarily the same as customers locations but if distance between pickup centers and some customers are less than or equal to a threshold, we assume that those customers are served by the truck.
The sets, parameters, and variables used in our model are defined in Table \ref{tab:Sets-Parameters-Variables-ProblemIII-CentersAnywhere}. Note that some of the variables are infinite dimensional. 
\begin{table}[!h]
\protect\caption{Sets, parameters, and variables used in the general (multi-center) model.}
\label{tab:Sets-Parameters-Variables-ProblemIII-CentersAnywhere}
\begin{center}
\begin{tabular}{c|p{0.72\textwidth}}  
      \hline
      \hline
       \textbf{Parameters/Variables} & \textbf{Description} \\
      \hline
      $d_{ij}$ & Distance between node $i\in D$ and $j\in C$ \\
      $\bm{y}_i$ & Vector of location of drone $i \in D$ \\
      $\bm{z}_j$ & Vector of location of customer $j \in C$ \\
      $\tau$ & Distance threshold for serving the customer by the truck (if the distance between truck's stop and a customer is less than or equal to this threshold) \\
      \hline
       $N_T$ & Total number of truck centers \\
       $\bm{x}_m$ & Vector of location of center $m$ with $\bm{x}_m \in R$ and $m = 1, \dots, N_T$ \\
       $\xi_{i\bm{x}_mj}$ &  
       				$D \times R \times C \xrightarrow{} \{0,1\}$. It is $1$ when drone travels from node $i\in D$ to center $\bm{x}_m$ then to node $j\in C$ and back to node $i$. It is $0$ otherwise. \\
	$\zeta_{\bm{x}_mj}$ & $R \times C \xrightarrow{} \{0,1\}$. It is $1$ when customer $j\in C$ is served by the truck that has stopped at location $\bm{x}_m$ and 0 otherwise \\
       $q_{i\bm{x}_m}$ & Total distance travelled by drone node $i$ assigned to center $\bm{x}_m$\\
       $Q_{\bm{x}_m}$ & Maximum distance among all $q_{i\bm{x}_m}$\\
       $\gamma_{(\bm{x}_m\bm{x}_n)}$ & 
       							$R \times R \xrightarrow{} \{0,1\}$. It is $1$ if the path from center $\bm{x}_m$ to $\bm{x}_n$ is used by the truck \\
       $u_{\bm{x}_m}$ & A dummy variable \\
       \hline
   \end{tabular}
\end{center}
\end{table}

We can model this problem in its most general sense as an infinite dimensional optimization problem that follows: \newpage
\begin{eqnarray}
 \minimize \; \frac{1}{v_T}  \iint_{\bm{x}_m\,, \bm{x}_n \in R} \|\bm{x}_m-\bm{x}_n\| \gamma_{(\bm{x}_m\bm{x}_n)} \, dA & + & \iint_{\bm{x}_m \in R} Q_{\bm{x}_m} \, dA     \qquad \qquad \st \nonumber \\ 
    \sum_{i\in D}\iint_{\bm{x}_m \in R} \xi_{i\bm{x}_mj}\, dA & = & 1\,, \qquad \qquad \forall j\in C \label{eq:Problem3CentersAnywhere-AllCustomersServed} \\
        \|\bm{x}_m - \bm{z}_j\|  & \leq &  \tau + M (1- \zeta_{\bm{x}_mj})  \,, \quad \forall j \in C\,, \; \forall \bm{x}_m \in R \label{eq:Problem3CentersAnywhere-TruckServing} \\
    (\|\bm{x}_m - \bm{y}_i\| + \|\bm{x}_m - \bm{z}_j\| + d_{ij}) \, \xi_{i\bm{x}_mj} & \leq &L\,, \qquad \qquad \forall i \in D\,, \; \forall j \in C\,, \; \forall \bm{x}_m \in R  \label{eq:Problem3CentersAnywhere-DroneRange} \\
   \frac{1}{v_D} \sum_{j\in C} (\|\bm{x}_m - \bm{y}_i\| + \|\bm{x}_m - \bm{z}_j\| + d_{ij}) \, \xi_{i\bm{x}_mj} \, (1-\zeta_{\bm{x}_mj}) & \leq & q_{i\bm{x}_m} \,, \qquad \forall i \in D\,, \; \forall \bm{x}_m \in R \label{eq:Problem3CentersAnywhere-CenterDroneTime} \\
    q_{i\bm{x}_m}  -  \frac{1}{v_D} \max_{j} d_{ji} \xi_{i\bm{x}_mj} & \leq & Q_{\bm{x}_m}\,, \qquad \forall i\in D\,, \; \forall \bm{x}_m \in R  \label{eq:Problem3CentersAnywhere-CenterMakespan} \\
    \iint_{\bm{x}_m \in R,\bm{x}_m\neq \bm{x}_n}{\gamma_{(\bm{x}_m\bm{x}_n)}} dA & = & 1\,, \qquad \qquad \forall \bm{x}_n \in R \label{eq:Problem3CentersAnywhere-TSPassignment1} \\
    \iint_{\bm{x}_n \in R,\bm{x}_n\neq \bm{x}_m}\gamma_{(\bm{x}_m\bm{x}_n)} dA & = & 1\,, \qquad \qquad \forall \bm{x}_m \in R \label{eq:Problem3CentersAnywhere-TSPassignment2}  \\
    \iint_{\bm{x}_m\,, \bm{x}_n \in R} \gamma_{\bm{x}_m\bm{x}_n} dA & = & N_T\,, \label{eq:Problem3CentersAnywhere-CentersNumber} \\
    u_{\bm{x}_m}-u_{\bm{x}_n} + N_T \; \gamma_{\bm{x}_m\bm{x}_n} & \leq & N_{T}-1 \,, \qquad \forall \; 2\leq u_{\bm{x}_m} \neq u_{\bm{x}_n} \leq N_T \label{eq:Problem3CentersAnywhere-TSPsubtour} \\
    d_{i\bm{x}_mj}, \; d_{\bm{x}_m\bm{x}_n}, \; q_{i\bm{x}_m}, \; Q_{\bm{x}_m} & \geq & 0\,, \qquad \qquad \; \forall  i\in D\,, \; \forall j\in C\,, \;  \forall \bm{x}_m \in R \nonumber \\
    \zeta_{\bm{x}_mj}, \; \xi_{i\bm{x}_mj}, \; \gamma_{\bm{x}_m\bm{x}_n} & \in & \{0,1\} \,, \qquad \forall j\in C\,, \; \forall \bm{x}_m\,, \bm{x}_n \in R \nonumber \\
    N_T, \; u_{\bm{x}_m}  & \in & \N \,,  \qquad \qquad \forall \bm{x}_m \in R \nonumber %\\
\end{eqnarray}
The objective function minimizes the sum of the time that truck takes to travel between the pickup centers on its route and the total time that it takes to serve clusters of customers with drones at these stopping points. 
Constraint \eqref{eq:Problem3CentersAnywhere-AllCustomersServed} ensures all customer are visited once. 
Constraint \eqref{eq:Problem3CentersAnywhere-TruckServing} controls whether a customer is served by the truck.
Constraint \eqref{eq:Problem3CentersAnywhere-DroneRange} enforces the drone battery limit. Constraints \eqref{eq:Problem3CentersAnywhere-CenterDroneTime} and \eqref{eq:Problem3CentersAnywhere-CenterMakespan} calculates the time cost in each center, excluding the time spent on the last edge of last route of each drone. 
Constraint \eqref{eq:Problem3CentersAnywhere-TSPassignment1}-\eqref{eq:Problem3CentersAnywhere-TSPsubtour} are travelling salesman problem (TSP) constraints for the truck as well as to calculate number of pickup centers.

This is obviously a complex problem to solve and to make it easier we need some relaxations. To do this, we assume that the truck only stops at a customer location and while stopping there that location will serve as a pickup center for drones as well. Luckily, this simple and very reasonable relaxation, makes the problem much easier to model and solve. This assumption is confirmed to be valid and beneficial by findings of \cite{lee2015dynamic} that showed the benefits of having dedicated drivers in parallel to crowd drivers in the context of dynamic ride-sharing, where the dedicated drivers, in addition to satisfying part of the demand, could also serve riders for whom no ad-hoc drivers are available. 
The sets, parameters, and variables used in this new model are defined in Table \ref{tab:Sets-Parameters-Variables-ProblemIII-CentersOnCustomers}.
\begin{table}[!h]
\protect\caption{Sets, parameters, and variables used in the general (multi-center) model.}
\label{tab:Sets-Parameters-Variables-ProblemIII-CentersOnCustomers}
\begin{center}
\begin{tabular}{c|p{0.72\textwidth}}  
      \hline
      \hline
      \textbf{Parameters/Variables} & \textbf{Description} \\
      \hline
      $d_{pp'}$ & Distance between node $p\in C$ and $p'\in C$ \\
      $d^p_{ij}$ & Route length going from node $i\in D$ to the center node $p\in C$ and then to customer node $j\in C$ and back to node $i$\\
      $\delta_{ji}$ & Distance of customer node $j\in C$ to the drone node $i\in D$ \\
      \hline
       $y_p$ & Binary decision variable that is 1 if node $p\in C$ is served by the truck and 0 otherwise \\
       $x^p_{ij}$ & Binary decision variable that is $1$ when a drone travels from node $i\in D$ to the center node $p\in C$ and then to customer node $j\in C$ and back to node $i$ and is $0$ otherwise \\
       $z_{ip}$ & Binary decision variable that is $1$ when drone $i\in D$ is used for delivery at center node $p\in C$ and is $0$ otherwise \\
       $q_{ip}$ & Total travel time of the drone with base at node $i\in D$ that is assigned to center $p\in C$\\
       $Q_p$ & Maximum time spent by all drones assigned to center $p\in C$\\
       $\Delta_{ip}$ & Maximum length of the last tour segment among all tours formed by drone $i \in D$ through center $p \in C$\\
       $\eta^{p}_{ij}$ & Binary decision variable equal to 1 if $d_{ji}x^p_{ij} = \Delta_{ip}$ and 0 otherwise \\
       $\gamma_{pp'}$ & Binary decision variable equal to 1 if the path from $p$ to $p'$ has been used by the truck \\
       \hline
   \end{tabular}
\end{center}
\vspace{-10pt}
\end{table}

An MILP model for this problem can be written as following:  
\begin{eqnarray}
    \minimize \;\;\; \frac{1}{v_T}\sum_{p\in C}\sum_{p'\in C} d_{pp'}\gamma_{pp'} & + & \sum_{p} Q_p \qquad \qquad \qquad \st \nonumber \\
    \sum_{p\in C} y_p & \geq & 1 \; , \label{eq:Problem3CentersOnCustomersMinOneTruckNode} \\
    \sum_{i\in D}\sum_{p\in C} x^p_{ij} & = & 1-y_j \; , \quad \qquad \qquad \forall j\in C \label{eq:Problem3CentersOnCustomersAssignment} \\
    \sum_{i\in D} z_{ip} & \leq & M y_p \; , \quad \qquad \qquad  \forall p\in C  \label{eq:Problem3CentersOnCustomersTruckThenDrone} \\
    \sum_{p\in C} z_{ip} & \leq & 1 \; , \quad \qquad \qquad \qquad \forall i\in D  \label{eq:Problem3CentersOnCustomersEachDroneOneCenter} \\
    \sum_{j \in C} x^p_{ij} & \leq & M z_{ip} \; , \quad \qquad \qquad \forall i\in D, p\in C  \label{eq:Problem3CentersOnCustomersIfDroneAtCenter} \\
    d^p_{ij} x^p_{ij} & \leq & L \; , \quad \qquad \qquad \qquad \forall i\in D, j,p\in C  \label{eq:Problem3CentersOnCustomersDroneRange} \\ 
        \frac{1}{v_D}\sum_{j\in C} d^p_{ij}x^p_{ij} - M(1-z_{ip}) & \leq & q_{ip} \; ,  \qquad \qquad \qquad \forall i\in D, p\in C \label{eq:Problem3CentersOnCustomersCenterDroneTime} \\
    q_{ip} -  (1/v_D) \max_{j} \delta_{ji} x^p_{ij} & \leq & Q_p \; , \qquad \qquad \qquad \forall i\in D, p\in C \label{eq:Problem3CentersOnCustomersCenterMakespan} \\
    2\gamma_{pp'} & \leq & y_{p} + y_{p'} \; , \qquad \qquad \forall p\neq p' \in C 
    \label{eq:Problem3CentersOnCustomersTSPtruckNodesMakePath} \\
    \sum_{p\in C} \gamma_{pp'} & = & y_{p'} \; , \qquad \qquad \qquad \forall p'\in C \label{eq:Problem3CentersOnCustomersTSPassignment1} \\
    \sum_{p'\in C} \gamma_{pp'} & = & y_p \; , \qquad \qquad \qquad \forall p \in C \label{eq:Problem3CentersOnCustomersTSPassignment2} \\
    \sum_{p\in S} \; \sum_{p'\in S,\, p'\neq p} \gamma_{pp'} & \leq & |S| -1 \; , \qquad \qquad \forall S \subset C\,, \; 2 \leq |S| \leq |C|-2 \label{eq:Problem3CentersOnCustomersTSPsubtour} \\
    M \gamma_{pp} + \sum_{j\in C} y_{j} -1 & \leq & M \; , \qquad \qquad \qquad \forall p \in C \label{eq:Problem3CentersOnCustomersTSPorStar} \\
    x^p_{ij}, \; \gamma_{pp'}, \; y_p & \in & \{0,1\} \;,  \qquad \qquad \forall i\in D, j,p,p' \in C  \nonumber \\
    q_{ip}, \; z_{ip}, \; Q_p & \geq & 0 \;,  \qquad \qquad \qquad \forall i\in D, p\in C  \nonumber %\\
\end{eqnarray}
The objective function minimizes the sum of the time that truck takes to travel between the pickup centers on its route and the total time that it takes to serve customers with crowdsourced drones at these stopping points. 
Constraint (\ref{eq:Problem3CentersOnCustomersMinOneTruckNode}) ensures that at least one customer node is set to be truck node. 
Constraint (\ref{eq:Problem3CentersOnCustomersAssignment}) ensures all customers are visited once, either by the truck or by one of the drones. Constraint (\ref{eq:Problem3CentersOnCustomersTruckThenDrone}) ensures that a node can serve as pickup center for drones if it is visited by the truck. 
Constraint (\ref{eq:Problem3CentersOnCustomersEachDroneOneCenter}) enforces the assumption that each drone can only be used at one truck center.
Constraint \eqref{eq:Problem3CentersOnCustomersIfDroneAtCenter} ensures that customers are served by a drone through a center node if that drone is used at that center.
Constraint (\ref{eq:Problem3CentersOnCustomersDroneRange}) ensures that a drone cannot fly more than its range (battery limit). 
Constraints (\ref{eq:Problem3CentersOnCustomersCenterDroneTime}) and (\ref{eq:Problem3CentersOnCustomersCenterMakespan}) are used to find the time spent at each stop of the truck to serve a group of customers. 
Note that truck does not have to wait at a pickup center after all packages in that neighborhood are delivered so the time spent on the last segment of the last route for each drone is deducted. This implies that we may prefer to assign the drone to a farther customer.
Constraint \eqref{eq:Problem3CentersOnCustomersTSPtruckNodesMakePath} is to make sure that the truck tour is formed through truck nodes.
Constraints (\ref{eq:Problem3CentersOnCustomersTSPassignment1})--(\ref{eq:Problem3CentersOnCustomersTSPsubtour}) are TSP constraints for the truck. Note that we treat the set of constraints in \eqref{eq:Problem3CentersOnCustomersTSPsubtour} as lazy constraints and add them to our model in the implementation in a lazy fashion only when they are needed.
Finally, constraint \eqref{eq:Problem3CentersOnCustomersTSPorStar} allows the model to switch from having a TSP tour to a one center star shape (hub-and-spoke) model when the range and speed of drones are large enough.

To linearize the 
 $\max_{j} \delta_{ji}x^p_{ij}$ term in \eqref{eq:Problem3CentersOnCustomersCenterMakespan}, we replace this term with $\Delta_{ip}$ and constraint \eqref{eq:Problem3CentersOnCustomersCenterMakespan} with
\begin{eqnarray}
 \hspace{0.6in}   q_{ip} -  \frac{1}{v_D} \Delta_{ip} & \leq & Q_p \; , \qquad \qquad \qquad \forall i\in D, p\in C \label{eq:Problem3CentersOnCustomersCenterMakespanLinearized}
\end{eqnarray} 
and add the following four constraints
\begin{eqnarray} 
 \hspace{1.7in} \delta_{ji}x^p_{ij} & \leq & \Delta_{ip} \;, \quad \qquad \qquad \qquad \forall i\in D, p\in C, j\in C  \label{eq:Problem3CentersOnCustomersCenterMakespanLinearized1} \\
 \Delta_{ip} & \leq & \delta_{ji}x^p_{ij} + M(1-\eta^p_{ij}) \;, \quad \forall i\in D, p\in C, j\in C  \label{eq:Problem3CentersOnCustomersCenterMakespanLinearized2} \\   
 \sum_{j\in C} \eta^p_{ij} & \geq & z_{ip} \;, \qquad \qquad \qquad \qquad \forall i\in D, p\in C \label{eq:Problem3CentersOnCustomersCenterMakespanLinearized3} \\
  \eta^p_{ij} & \leq & x^p_{ij} \qquad \qquad \qquad \qquad \forall i\in D, p\in C, j\in C  \label{eq:Problem3CentersOnCustomersCenterMakespanLinearized4} \\
 \eta^p_{ij} & \in & \{0,1\} \quad \qquad \qquad \qquad \forall i\in D, p\in C, j\in C \nonumber 
\end{eqnarray} 

to form an MILP model.
We solved this model for small instances but for larger problems we rely on a heuristic algorithm.

\subsection{Problem IV--(Multiple Centers with Revisiting)} 	
\label{subsec:ProblemIV}
In this problem we consider multiple pickup centers and assume that after each delivery, drones could revisit the pickup center for another delivery before going back to their home bases. Similar to the latter part of the previous section, we assume the pickup centers (truck nodes) are a subset of customer nodes. 
Also, similar to Problem II, we define a feasible drone route as a route that starts from a drone node $i\in D$, satisfies the demand of one or multiple customer nodes after picking their packages up at a center node, and comes back to the starting drone node while not violating the drone range $L$. We assume that the revisits happen to the same center and the drone could start serving through another center after going back to its base. Relaxing this assumption is doable but comes at huge computational costs. We consider no restriction on the number of revisits for each drone.
We pre-calculate all such routes.
The sets, parameters, and variables used in our model are defined in Table \ref{tab:Sets-Parameters-Variables-ProblemIII-CentersOnCustomers}.

\begin{table}[!h]
\protect\caption{Sets, parameters, and variables used in the general (multi-center) model.}
\label{tab:Sets-Parameters-Variables-ProblemIV-CentersOnCustomers}
\begin{center}
\begin{tabular}{c|p{0.72\textwidth}}  
      \hline
      \hline
      \textbf{Sets/Parameters/Variables} & \textbf{Description} \\
      \hline
      $\Pi^p_{i}$ & Set of all feasible routes that start from drone node $i\in D$ and go through pickup center $p \in C$\\
      $\Pi^p_{i}(j)$ & Set of all feasible routes that start from drone node $i\in D$, go through pickup center $p \in C$, and traveled through customer node $j \in C$\\
      \hline
      $d_{pp'}$ & Distance between node $p\in C$ and $p'\in C$ \\
      $d_{\pi^p_{ik}}$ & Length of route $\pi^p_{ik} \in \Pi^p_i$ for drone node $i \in D$\\
      $\delta_{\pi^p_{ik}}$ & Length of the last line segment of route $\pi^p_{ik}$ \\
      \hline
       $y_p$ & Binary decision variable that is 1 if node $p\in C$ is served by the truck and 0 otherwise \\
       $x_{\pi^p_{ik}}$ & Binary decision variable that is $1$ when a drone travels along the route $\pi^p_{ik}\in \Pi^p_{i}$ and is $0$ otherwise \\
       $z_{ip}$ & Binary decision variable that is $1$ when drone $i\in D$ is used for delivery at center node $p\in C$ and is $0$ otherwise \\
       $q_{ip}$ & Total travel time of the drone with base at node $i\in D$ that is assigned to center $p\in C$\\
       $Q_p$ & Maximum time spent by all drones assigned to center $p\in C$\\
       $\Delta_{ip}$ & Maximum length of the last tour segment among all routes $\pi^p_{ik}\in \Pi^p_{i}$ \\
       $\eta_{\pi^p_{ik}}$ & Binary decision variable equal to 1 if $\delta_{\pi^p_{ik}}x_{\pi^p_{ik}} = \Delta_{ip}$ and 0 otherwise \\
       $\gamma_{pp'}$ & Binary decision variable equal to 1 if the path from $p$ to $p'$ has been used by the truck \\
       \hline
   \end{tabular}
\end{center}
\vspace{-16pt} 
\end{table}
An MILP model for this problem can be written as following:  
\begin{eqnarray}
    \minimize \;\;\; \frac{1}{v_T}\sum_{p\in C}\sum_{p'\in C} d_{pp'}\gamma_{pp'} & + & \sum_{p} Q_p \qquad \qquad \st \nonumber \\
    \sum_{p\in C} y_p & \geq & 1 \; , \label{eq:Problem4CentersOnCustomersMinOneTruckNode} \\
    \sum_{p\in C}\sum_{i\in D} \sum_{\pi^p_{ik}\in \Pi^p_i(j)}x_{\pi^p_{ik}} & = & 1-y_j \; , \quad \qquad \forall j\in C \label{eq:Problem4CentersOnCustomersAssignment} \\
    \sum_{i\in D} z_{ip} & \leq & M y_p \; ,  \qquad \qquad  \forall p\in C  \label{eq:Problem4CentersOnCustomersTruckThenDrone} \\
    \sum_{p\in C} z_{ip} & \leq & 1 \; , \quad \qquad \qquad \qquad \forall i\in D  \label{eq:Problem4CentersOnCustomersEachDroneOneCenter} \\
        \sum_{\pi^p_{ik}\in \Pi^p_i} x_{\pi^p_{ik}}  & \leq & M z_{ip} \; ,  \qquad \qquad \forall i\in D, p\in C  \label{eq:Problem4CentersOnCustomersIfDroneAtCenter} \\
   \frac{1}{v_D} \sum_{\pi^p_{ik}\in \Pi^p_i} d_{\pi^p_{ik}} x_{\pi^p_{ik}} - M(1-z_{ip}) & \leq & q_{ip} \; , \quad \qquad \qquad \forall i\in D, \;  p\in C  \label{eq:Problem4CentersOnCustomersCenterDroneTime} \\
    q_{ip} - (1/v_D) \max_{\pi^p_{ik}\in \Pi^p_i} \delta_{\pi^p_{ik}}x_{\pi^p_{ik}}  & \leq & Q_p \; , \quad \qquad \qquad \forall i\in D, \; p\in C \label{eq:Problem4CentersOnCustomersCenterMakespan} \\
        2\gamma_{pp'} & \leq & y_{p} + y_{p'} \; , \qquad \qquad \forall p\neq p'\in C \label{eq:Problem4CentersOnCustomersTSPtruckNodesMakePath} \\
    \sum_{p\in C} \gamma_{pp'} & = & y_{p'} \; , \qquad \qquad \qquad \forall p'\in C \label{eq:Problem4CentersOnCustomersTSPassignment1} \\
    \sum_{p'\in C} \gamma_{pp'} & = & y_p \; , \qquad \qquad \qquad \forall p \in C \label{eq:Problem4CentersOnCustomersTSPassignment2} \\
    \sum_{p\in S} \; \sum_{p'\in S,\, p'\neq p} \gamma_{pp'} & \leq & |S| -1 \; , \quad \qquad \forall S \subset C\,, \; 2 \leq |S| \leq |C|-2 \label{eq:Problem4CentersOnCustomersTSPsubtour} \\
    M \gamma_{pp} + \sum_{j\in C} y_{j} -1 & \leq & M \; , \qquad \qquad \qquad \forall p \in C \label{eq:Problem4CentersOnCustomersTSPorStar} \\
    x_{\pi^p_{ik}}, \; \gamma_{pp'},\; y_p & \in & \{0,1\} \; , \; \qquad \qquad \forall i\in D,\; p,p'\in C,\; \pi^p_{ik} \in \Pi^p_i \nonumber \\
    q_{ip}, \; Q_p & \geq & 0 \;, \qquad \qquad \qquad \forall i\in D,\; p\in C \nonumber %\\
\end{eqnarray}
The objective function minimizes the sum of the time that spent on the truck route and the drone routes. 
Constraint (\ref{eq:Problem4CentersOnCustomersMinOneTruckNode}) ensures that at least one customer node is set to be truck node (pickup center). 
Constraint (\ref{eq:Problem4CentersOnCustomersAssignment}) ensures either the truck or one of the drones visits each of the customers. 
Constraint (\ref{eq:Problem4CentersOnCustomersTruckThenDrone}) ensures that a customer nodes can only be used as pickup center for drones if that node is visited by the truck. 
Constraint (\ref{eq:Problem4CentersOnCustomersEachDroneOneCenter}) enforces the assumption that each drone can only be used at one truck center.
Constraint \eqref{eq:Problem4CentersOnCustomersIfDroneAtCenter} enforces drone routes to be formed at pickup centers.
Constraints (\ref{eq:Problem4CentersOnCustomersCenterDroneTime}) and (\ref{eq:Problem4CentersOnCustomersCenterMakespan}) are used to find the time spent at each stop of the truck to serve a group of customers. Note that at each center the truck may not wait until the last drone route is completed, which may, at times, prioritize a farther customer over a closer customer for a drone.
Constraint \eqref{eq:Problem4CentersOnCustomersTSPtruckNodesMakePath} is to make sure that the truck tour is formed through truck nodes.
Constraints (\ref{eq:Problem4CentersOnCustomersTSPassignment1})--(\ref{eq:Problem4CentersOnCustomersTSPsubtour}) are TSP constraints for the truck. Note that we treat the set of constraints in \eqref{eq:Problem4CentersOnCustomersTSPsubtour} as lazy constraints and add them to our model in the implementation in a lazy fashion only when they are needed.
Finally, constraint \eqref{eq:Problem4CentersOnCustomersTSPorStar} allows the model to switch from having a TSP tour to a one center star shape (hub-and-spoke) model when the range and speed of drones are large enough.

To linearize the $\max_{\pi^p_{ik}\in \Pi^p_i} \delta_{\pi^p_{ik}}x_{\pi^p_{ik}}$ term in \eqref{eq:Problem4CentersOnCustomersCenterMakespan}, we replace this term with $\Delta_{ip}$ and constraint \eqref{eq:Problem4CentersOnCustomersCenterMakespan} with
\begin{eqnarray}
 \hspace{0.1in}   q_{ip} -  \frac{1}{v_D} \Delta_{ip} & \leq & Q_p \; , \qquad \qquad \forall i\in D, p\in C \label{eq:Problem4CentersOnCustomersCenterMakespanLinearized}
\end{eqnarray} 
and add the following four constraints
\begin{eqnarray}
 \hspace{1.7in}  \delta_{\pi^p_{ik}}x_{\pi^p_{ik}} & \leq & \Delta_{ip} \;, \qquad \qquad \qquad \forall i\in D, p\in C, \, \pi^p_{ik}\in \Pi^p_i \label{eq:Problem4CentersOnCustomersCenterMakespanLinearized1} \\
 \Delta_{ip} & \leq & \delta_{\pi^p_{ik}}x_{\pi^p_{ik}} + M(1-\eta_{\pi^p_{ik}})  \;, \;\; \forall i\in D, p\in C, \pi^p_{ik}\in \Pi^p_i  \label{eq:Problem4CentersOnCustomersCenterMakespanLinearized2} \\  
 \sum_{\pi^p_{ik}\in \Pi^p_i} \eta_{\pi^p_{ik}} & \geq & z_{ip} \;, \label{eq:Problem4CentersOnCustomersCenterMakespanLinearized3} \\
   \eta_{\pi^p_{ik}} & \leq & x_{\pi^p_{ik}} \;, \;\; \qquad \qquad \forall i\in D, p\in C, \pi^p_{ik}\in \Pi^p_i \label{eq:Problem4CentersOnCustomersCenterMakespanLinearized4} \\
 \eta_{\pi^p_{ik}} & \in & \{0,1\} \;, \qquad \qquad \forall i\in D, p\in C, \, \pi^p_{ik}\in \Pi^p_i \nonumber 
\end{eqnarray} 

to form an MILP model. We also solved this model for small instances but for larger problems we rely on our heuristic algorithm.

\section{Algorithm}
\label{sec:algorithm}
It is clear that the problem is NP-hard and we have to take a sub-optimal approach to solve the problem. 
Our algorithm consists of several sub-routines as explained in the following. The high-level steps of our algorithm are as follows:
\begin{description}
    \item [Step 1:] Running a \emph{binary search} and the \emph{$k$-means clustering algorithm} to find the minimum feasible $k$ to cluster all customers into $k$ clusters with radius $L/4$ and finding the cluster centers.
    \item [Step 2:] Partitioning the region with the \emph{Voronoi tessellation} generated by the $k$ cluster center points found in Step 1 and assigning customers and drone bases to their nearest centers.    
    \item [Step 3:]  Running a \emph{Tabu Search algorithm} to solve the sub-problem in each cluster where parcels are assigned to drones to be delivered to customers in a way that minimizes total delivery time.     
    \item [Step 4:]  Solving the truck tour by the \emph{Lin-Kernighan-Helsgaun (LKH) algorithm} to go through all cluster centers as well as all customer nodes that do not have any drone node in that cluster.
\end{description}

\begin{figure}[hbt]
    \centering
    \begin{subfigure}[b]{0.25\linewidth}
        \includegraphics[width=\linewidth]{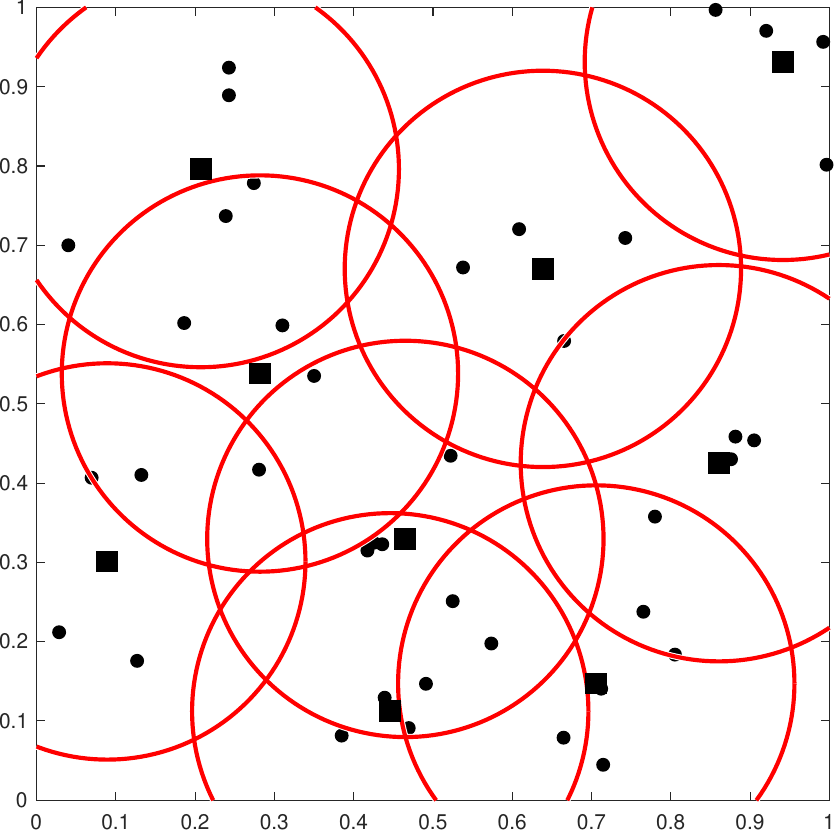}
        \caption{$k$-means clustering}
        \label{fig:ClusteringStep}
    \end{subfigure}
    \qquad
    \begin{subfigure}[b]{0.25\linewidth}
        \includegraphics[width=\linewidth]{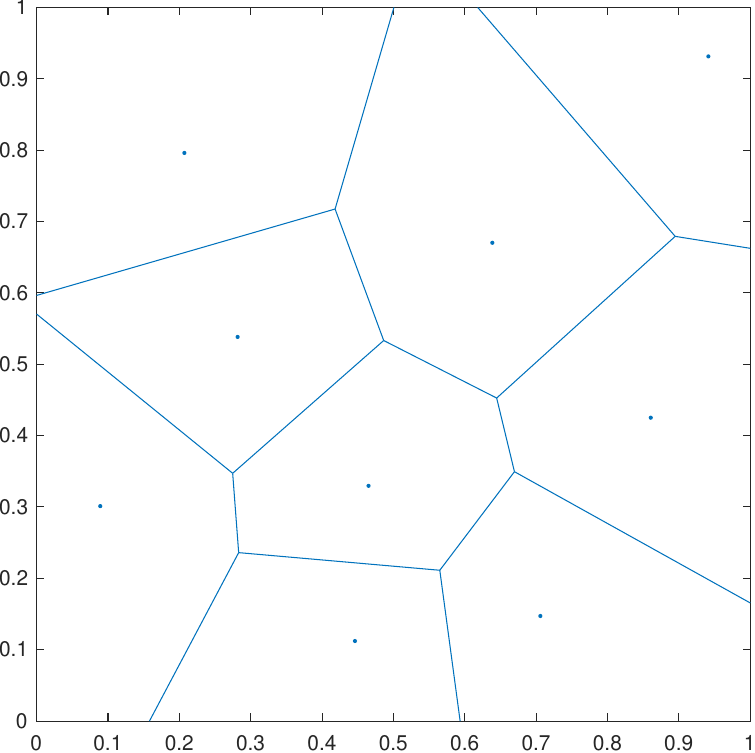}
        \caption{Voronoi Partition}
        \label{fig:VoronoiStep}
    \end{subfigure}
    \\
    \begin{subfigure}[b]{0.27\linewidth}
        \includegraphics[width=\linewidth]{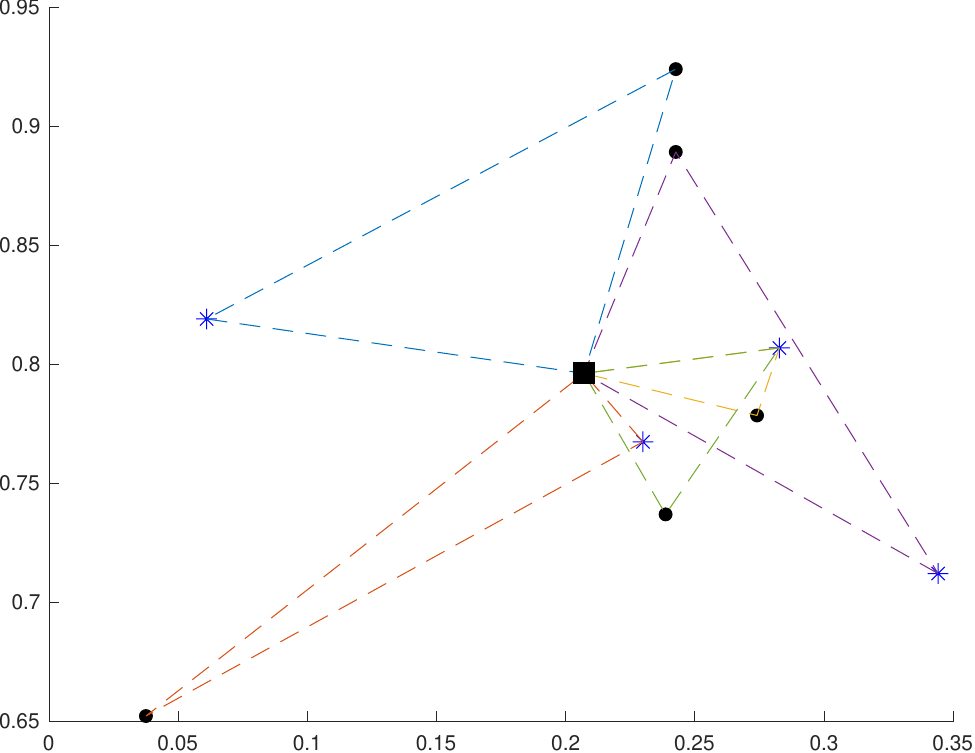}
        \caption{One cluster result}
        \label{fig:OneClusterTS}
    \end{subfigure}
    \quad 
    \begin{subfigure}[b]{0.25\linewidth}
        \includegraphics[width=\linewidth]{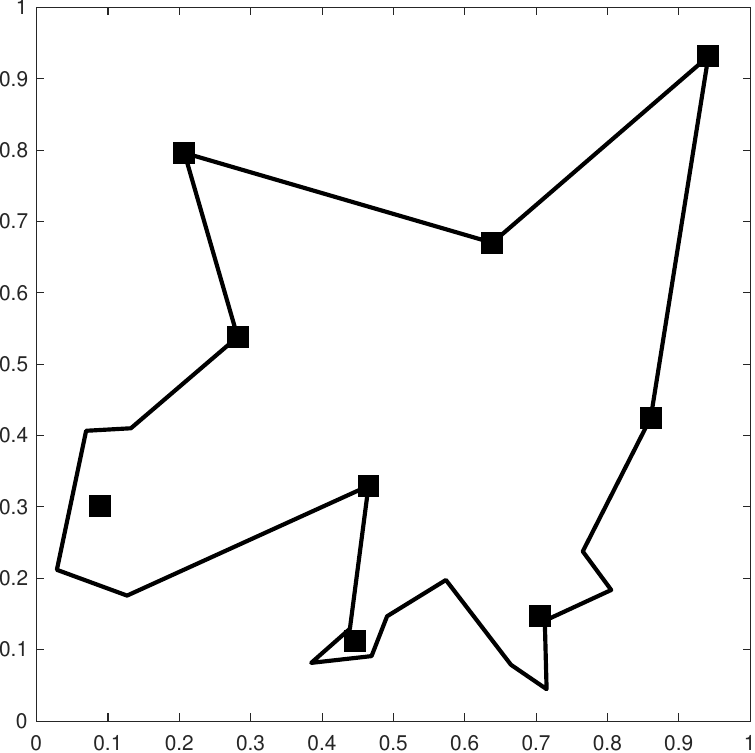}
        \caption{Truck route}
        \label{fig:TruckRouteTSP}
    \end{subfigure}
        \\
    \begin{subfigure}[b]{0.3\linewidth}
       \includegraphics[width=\linewidth]{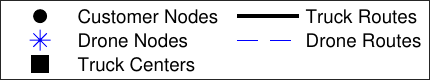}
    \end{subfigure}
    \caption{Illustration of the main steps of the Algorithm for a problem with 40 customers and $k=9$ clusters. The $k$-means clustering step is shown in (\subref{fig:ClusteringStep}) and the Voronoi partitioning step is illustrated in (\subref{fig:VoronoiStep}). (\subref{fig:OneClusterTS}) shows the routes for the drones in the upper left cluster and (\subref{fig:TruckRouteTSP}) shows the truck route going through cluster centers and customer nodes that are not served by the drones.}
    \label{fig:ocrTR}
\end{figure}

We first cluster all customers into $k$ groups, where each group will have a center that will serve as the truck's stoping point, by solving a \emph{$k$-means clustering problem} to ensure all customers are within radius $L/4$ of their closest centers. A similar idea for finding truck's stopping points for launching drones in the centralized system of combined truck and drone has been previously shown to be efficient \cite{Ferrandez2016optimization}.
To solve the $k$-means clustering problem we use \emph{Lloyd's algorithm} \cite{lloyd1982least}, which is described in \cref{sec:LloydsAlg} of the Online Supplement of this paper. This minimizes the total squared Euclidean distances between customer nodes and their closest cluster centers. The choice of radius $L/4$ is to make sure that with one full battery the drones in each cluster can finish the delivery and return to their base. The number $k$ is the smallest integer that makes this feasible and will be found using a binary search. The feasibility here means that all customers are located in at least one circle that is centered at a cluster center with radius $L/4$. The clustering step is shown in \cref{fig:ClusteringStep}. The final geographic partitioning, i.e., identifying the boundaries of each cluster, is done using a \emph{Voronoi partitioning scheme} \cite{voronoi1908nouvelles_a,voronoi1908nouvelles_b,de2000computational} with the cluster centers as the generators. This determines the assignment of customer nodes and drone nodes to their closest cluster centers. 
These nodes would be assigned to the nearest center only if their distance to the center is less than $L/4$. 
Clusters with no drone will be served by the truck. The partitioning step is illustrated in \cref{fig:VoronoiStep}.

In the next step, in each cluster (Voronoi cell), a \emph{Tabu Search algorithm} \cite{glover1986future} will be used to solve the one center problem. The steps of a generic Tabu Search algorithm are explained in in \cref{sec:TSAlg} of the Online Supplement of this paper. In each cluster, several drones will fly from their origin  (base), go to the truck center, deliver the packages to the customer and go back to their base. We represent a solution of the Tabu Search with a list of all nodes. As an example consider the following solution
\begin{center}
          \begin{tabular}{|c|c|c|c|c|c|c|c|c|c|c|c|c|c|c|}
          \hline
          8 & 1 & 5& 7& 1& 2& 9& 1& 3& 10& 1& 4& 9& 1& 6\\
          \hline
          \end{tabular}
\end{center}
for a cluster (one center problem) with 5 customer nodes and 4 drone nodes. In this solution, $1$ represents the truck node, $2$ to $6$ represents the customer nodes and $7$ to $10$ represents the drone nodes. Therefore, a solution representation like this indicates that a drone from location 8 goes to the cluster center (truck node) at location 1 and then goes to customer node 5, while drones from locations 7 and 10 also go to the cluster center and deliver the packages to customer nodes 2 and 4, respectively. Drone from location 9 first serves the customer node 3 and comes back to its home base to change the battery, then it goes to serve the customer node 6. In this solution, four drones will deliver packages for five customers, and the solution for one cluster is shown in \cref{fig:OneClusterTS}, where the black square is the cluster center, blue stars are the customer-owned drones' location and small black circles are the customers' location. It is clear that each such list represents a unique solution of the problem but each solution can be represented with multiple such lists.

The Tabu Search algorithm is applied to this solution to find a new solution. During each iteration, an action will be applied to the current solution. An action is either to randomly switch the position of customer nodes (switching 3 and 5 in the following example)
\begin{center}
          \begin{tabular}{|c|c|c|c|c|c|c|c|c|c|c|c|c|c|c|}
          \hline
          8 & 1 & {\color{red}3} & 7& 1& 2& 9& 1&{\color{red}5} & 10& 1& 4& 9& 1& 6\\
          \hline
          \end{tabular}
\end{center}
or randomly change a drone node among all the drone nodes in that cluster (changing 7 to 9 in the following example). \vspace{-8pt}
\begin{center}
          \begin{tabular}{|c|c|c|c|c|c|c|c|c|c|c|c|c|c|c|}
          \hline
              8 & 1 & 5& {\color{red}9} & 1& 2& 9& 1& 3& 10& 1& 4& 9& 1& 6\\
          \hline
          \end{tabular}
\end{center}

This gives us an action list that has a size of $|C|(|C|-1)/2+|D|$.
     
In the final step of the algorithm, we solve a travelling salesman problem (TSP) to find the truck tour among all cluster centers and all customer nodes that do not have any drone node nearby. All customers in any cluster that has no drones in it, will be served by the truck. \emph{Lin-Kernighan-Helsgaun (LKH) algorithm} \cite{lin1973effective,helsgaun2000effective} is used for this purpose. The LKH algorithm is a local search heuristic to solve TSP with swapping operations inside the tours to make a new tour. It is an iterative algorithm with two set of links, with each link connecting two nodes, and in each step it will decide pairs of links that need to be switched by deleting links in the first set from the tour and adding the links from the second set to the tour (a generalization of the well-known 2-opt and 3-opt exchanges). \cref{fig:TruckRouteTSP} shows the TSP tour in our example. It can be seen that there are some cluster centers that are not on the truck route. This is because in those clusters all of the customers are served by the truck and thus there is no reason for the truck to go to the cluster center and wait there for the drones anymore. 

\section{Improvements of the Algorithm}
\label{sec:Improvements}
Next, we will show that some further improvements can be applied to the proposed model with sharing drones. In this paper, we analyze four such improvements. In the first three, we try to find the best location and the best number of truck centers by moving these centers to one of the customer nodes and merging these centers, if we could still cover the all customers using the new centers. In the last improvement, we try to utilize the battery of drones, that is if the battery remaining in the drones can deliver another package, instead of going back to the base, we let them to go back to the cluster center to pick up the next package immediately after delivering the previous one. 

\subsection{Moving the Cluster Centers to Their Nearest Feasible Customer Nodes}
\label{subsec:MovingCentersNearestCustomerNode}
Both of our MILP models developed for Problems III and IV assume pickup centers (truck nodes) as a subset of customer nodes. We will also show in \cref{sec:MovingCentersImpact} that in our algorithm moving the cluster centers to some well-selected customer nodes, provided that the new pickup center still could serve all of its assigned customers via the available assigned drones to that center, improves the performance of the algorithm. Therefore, to make the comparison between the optimization models and the algorithm valid, we incorporate this modification into our algorithm.

\begin{figure}[t]
        \centering
        \includegraphics[width=0.8\textwidth]{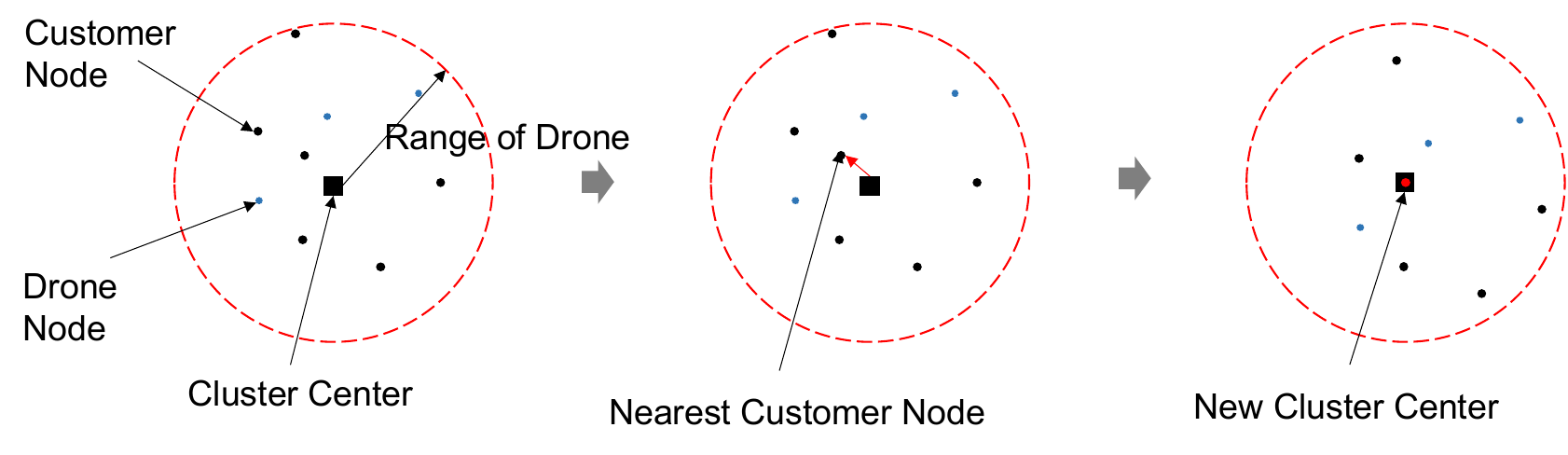}
        \caption{Moving cluster centers}
        \label{fig:mc1}
\end{figure}

In the first modification we move the cluster center to the nearest customer node to that center provided that the center could still cover all of its assigned customers using its assigned drones as illustrated in \cref{fig:mc1}. While we can save one drone trip to the customer that is now served by the truck, the length of the truck route or the routes of the other drones could either decrease or increase. We will analyze this trade-off in Section \ref{sec:MovingCentersImpact}. \cref{alg:RouteFinderMC1} summarizes the steps of the algorithm with this modification. One could take an alternative approach by incorporating the \emph{$k$-medoids algorithm} that automatically selects one of the central customers in each cluster as the cluster center \cite{kaufman1990partitioning}.

\subsection{Moving the Cluster Centers to Nearest Feasible Customer Node to Area Center}
\label{subsec:MovingCentersNearestCustomerNode2AreaCenter}

The modification in Algorithm 1 may provide a slight improvement. However, we can do better by moving the cluster centers strategically having the goal of making the truck route shorter. The intuition behind this is due to the general understanding that the drones are faster than trucks especially in urban delivery operations. Also, in the centralized alternative (TSPD) it is known that for having any significant gain by combining a truck and a drone in delivery operations we must have the speed of the drone at least twice as much as the speed of the truck \cite{Ferrandez2016optimization,agatz2018optimization} and it is most often assumed to have drones faster than the truck. Thus, since drone has speed advantage, we expect to gain in delivery time by decreasing the length of truck route and increasing the length of drone routes. We will analyze this trade-off in Section \ref{sec:MovingCentersImpact}. Therefore, we modify the algorithm as follows. 

We first find the centroid of the cluster centers (truck nodes). Let's call it $\bm{g}$. Then, for each cluster center we find one of its feasible assigned customer nodes closest to $\bm{g}$ among all other assigned customers to this center and closer to $\bm{g}$ than the cluster center itself. We choose this customer node, if existing, as the new cluster center (truck node). The steps of the algorithm with this modification are summarized in \cref{alg:RouteFinderMC2}.

\subsection{Merging Centers}
\label{subsec:MergingCenters}
We can still do a significant improvement to Algorithm \ref{alg:RouteFinderMC2} by trying to merge any two centers if all customers assigned to one could be served via the other by using drones assigned to either of these two centers. This comes from our assumption that the drone is faster than the truck and making the truck route shorter by reducing the number of centers could benefit us in the total delivery time. We first sort centers found in Algorithm \ref{alg:RouteFinderMC2} with respect to their distance from $\bm{g}$ in a descending order. We then start from the the first center in the list, which is the farthest from $\bm{g}$, and check the possibility of its merge into its closest center. If it is feasible, we assign this center, which itself is a customer node, and all customer nodes assigned to this center, to its closest center and we serve them using combined set of drones assigned to these two centers. The feasibility here is defined as whether drones in this combined set could serve all customers assigned to these two centers via the latter center while staying within their flying range $L$. This also limits the impact of the very restricting constraint of considering a range of $L/4$ in the initial clustering step. Once we have checked the possibility of a merge for all $k$ centers and if at least one merge is happened, we repeat this step with the remaining number of centers. This process continues until no more merge is possible. The pseudocode of the algorithm with this modification is presented in \cref{alg:RouteFinderMC3}.

\twocolumn
\begin{algorithm}[h]
\protect\caption{\label{alg:RouteFinderMC1}$\noun{RechargingRouteFinderMC1}\:(R,C,D,L)$  \hfill 
Generates the truck route and all the drone routes after finding the best location for the truck centers in the continuous space and then moving those centers to the nearest customer's place.}
\SetAlgoLined
\BlankLine
\KwIn{Input region $R$, set of customer nodes $C$, set of drone nodes $D$, and drone range $L$.}
\KwOut{Truck route $T_0$ and routes of drones $T_j\,, \; j\in D$ used for serving the customers.}
\BlankLine
\tcc{**********************************************}
Let $T_0$ denote the truck route and $T_j$ denote the route for drone $j\in D$\;
Find a feasible $k$ and cluster all customer nodes in $C$ into $k$ clusters within distance $L/4$ using a binary search and $k$-means clustering algorithms\;
Let $\bm{m}_1,...,\bm{m}_k$ be the cluster centers and generate their Voronoi tessellation\;
Let $V_1,...,V_k$ be the $k$ Voronoi cells and use them to assign all customer and drone nodes to their corresponding centers if their distance to the center is less than $L/4$\;
\For{$\ell=1,...,k$}{
  Let $\bm{z}_{\ell} = \argmin_{\bm{z}_j \in V_{\ell}} \|\bm{z}_j - \bm{m}_{\ell}\|\,, \; \forall j\in C$\;
  Let $\bm{m}_{\ell}^{\prime} = \bm{z}_{\ell}$\;
  Let $C_{\ell}=\{\bm{z}_j \in V_{\ell} \,, \; j\in C\}$\;
  \For{$\bm{z}_j \in C_{\ell} \; \&\& \; \bm{z}_j  \neq \bm{m}_{\ell}^{\prime}$}{
       \eIf{$\bm{z}_j$ can be served via $\bm{m}_{\ell}^{\prime}$ using drones in $V_{\ell}$}{
           \eIf{$|C_{\ell}\backslash \{\bm{m}_{\ell}^{\prime}\}| = 1$ \&\& 
            $2*\|\bm{z}_j  - \bm{m}_{\ell}^{\prime}\| / v_{T}  \leq (\|\bm{y}_i - \bm{m}_{\ell}^{\prime}\| + \|\bm{m}_{\ell}^{\prime} - \bm{z}_j\| + \|\bm{z}_j - \bm{y}_i\|)/v_{D} \,, \; \bm{y}_i \in V_{\ell}\,, \; i\in D$}{
            Let $\bm{z}_j$ be a truck node, i.e., to be served by the truck\;
           }{
           Assign $\bm{z}_j$ to $\bm{m}_{\ell}^{\prime}$, i.e., to be served by a drone\;
           }
           }{
           Let $\bm{z}_j$ be a truck node, i.e., to be served by the truck\;
       }
  }  
}
Update $k$ to be the updated number of truck nodes and reindex them as $\bm{m}_1^{\prime},...,\bm{m}_k^{\prime}$\;
Run the Tabu Search algorithm to find the drone routes in all Voronoi cells, i.e., $T_j\,, \; j \in D$\;
\If{$k=2$}{
 Set $T_0 = \{\bm{m}_1^{\prime}, \bm{m}_2^{\prime}\}$\;
 }
\If{$k\geq 3$}{
   Run the Lin-Kernighan-Helsgaun algorithm to find the truck route $T_0$ through $\bm{m}_1^{\prime},...,\bm{m}_k^{\prime}$\;
}
\KwRet{$T_0$ \textrm{and} $T_j\,, \; j\in D$}\;
\end{algorithm}

\begin{algorithm}[thb]
\protect\caption{\label{alg:RouteFinderMC2}$\noun{RechargingRouteFinderMC2}\:(R,C,D,L)$  \hfill 
Generates the truck route and all the drone routes. It first finds the best location for the truck centers and then moves those centers to the nearest customer's place. Finally, it finds the centroid of these centers and moves these centers to one of their assigned customers that is closest to this centroid and keeps the delivery feasible.}
\SetAlgoLined
\BlankLine
\KwIn{Input region $R$, set of customer nodes $C$, set of drone nodes $D$, and drone range $L$.}
\KwOut{Truck route $T_0$ and routes of drones $T_j\,, \; j\in D$ used for serving the customers.}
\BlankLine
\tcc{**********************************************}
Let $T_0$ denote the truck route and $T_j$ denote the route for drone $j\in D$\;
Find a feasible $k$ and cluster all customer nodes in $C$ into $k$ clusters within distance $L/4$ using a binary search and $k$-means clustering algorithms\;
Let $\bm{m}_1,...,\bm{m}_k$ be the cluster centers and generate their Voronoi tessellation\;
Let $V_1,...,V_k$ be the $k$ Voronoi cells and use them to assign all customer and drone nodes to their corresponding centers if their distance to the center is less than $L/4$\;
\For{$\ell=1,...,k$}{
  Let $\bm{z}_{\ell} = \argmin_{\bm{z}_j \in V_{\ell}} \|\bm{z}_j - \bm{m}_{\ell}\|\,, \; \forall j\in C$\;
  Let $\bm{m}_{\ell}^{\prime} = \bm{z}_{\ell}$\;
  Let $C_{\ell}=\{\bm{z}_j \in V_{\ell} \,, \; j\in C\}$\;
  \For{$\bm{z}_j \in C_{\ell} \; \&\& \; \bm{z}_j  \neq \bm{m}_{\ell}^{\prime}$}{
       \eIf{$\bm{z}_j$ can be served via $\bm{m}_{\ell}^{\prime}$ using drones in $V_{\ell}$}{
           \eIf{$|C_{\ell}\backslash \{\bm{m}_{\ell}^{\prime}\}| = 1$ \&\& 
            $2*\|\bm{z}_j  - \bm{m}_{\ell}^{\prime}\| / v_{T}  \leq (\|\bm{y}_i - \bm{m}_{\ell}^{\prime}\| + \|\bm{m}_{\ell}^{\prime} - \bm{z}_j\| + \|\bm{z}_j - \bm{y}_i\|)/v_{D} \,, \; \bm{y}_i \in V_{\ell}\,, \; i\in D$}{
            Let $\bm{z}_j$ be a truck node, i.e., to be served by the truck\;
           }{
           Assign $\bm{z}_j$ to $\bm{m}_{\ell}^{\prime}$, i.e., to be served by a drone\;
           }
           }{
           Let $\bm{z}_j$ be a truck node, i.e., to be served by the truck\;
       }
  }  
}
Update $k$ to be the updated number of truck nodes and reindex them as $\bm{m}_1^{\prime},...,\bm{m}_k^{\prime}$\;
Let $\bm{g}$ be the centroid of $\bm{m}_1^{\prime},...,\bm{m}_k^{\prime}$\;
\For{$\ell=1,...,k$}{
    Set $d_j = \|\bm{z}_j - \bm{g}\|\,, \; \forall j\in C$ and $d_{\bm{m}_{\ell}^{\prime}} = \|\bm{m}_{\ell}^{\prime} - \bm{g}\|$\;
    Let $A_{\ell}$ be an array of customer nodes in $V_{\ell}$ sorted in ascending order of $d_j\,, \; \forall j\in C$\;
    Find the customer node $\bm{z}_{\ell}$ with the smallest index in $A_{\ell}$ and $d_{\bm{z}_{\ell}} \leq d_{\bm{m}_{\ell}^{\prime}}$ that has at least one drone within distance $L/4$ and if chosen as cluster center we can still serve all other customers in $V_{\ell}$ using drones in $V_{\ell}$ and let $\bm{m}_{\ell}^{\prime\prime} = \bm{z}_{\ell}$. If no such point exists let $\bm{m}_{\ell}^{\prime\prime} = \bm{m}_{\ell}^{\prime}$\;
}
Run the Tabu Search algorithm to find the drone routes in all Voronoi cells, i.e., $T_j\,, \; j \in D$\;
\If{$k=2$}{
 Set $T_0 = \{\bm{m}_1^{\prime\prime}, \bm{m}_2^{\prime\prime}\}$\;
 }
\If{$k\geq 3$}{
   Run the Lin-Kernighan-Helsgaun algorithm to find the truck route $T_0$ through $\bm{m}_1^{\prime\prime},...,\bm{m}_k^{\prime\prime}$\;
}
\KwRet{$T_0$ \textrm{and} $T_j\,, \; j\in D$}\;
\end{algorithm}

\onecolumn

\begin{center}
 \scalebox{.80}{
  \begin{minipage}{\linewidth}
\begin{algorithm}[H]
\protect\caption{\label{alg:RouteFinderMC3}$\noun{RechargingRouteFinderMC3}\:(R,C,D,L)$ -- Generates the truck route and all the drone routes. It first finds the best location for the truck centers and then moves those centers to the nearest customer's place. It then finds the centroid of these centers and moves these centers to one of their assigned customers that is closest to this centroid and keeps the delivery feasible. Finally, it tries to merge these centers to shorten the truck tour.}
\SetAlgoLined
\BlankLine

\KwIn{Input region $R$, set of customer nodes $C$, set of drone nodes $D$, and drone range $L$.}
\KwOut{Truck route $T_0$ and routes of drones $T_j\,, \; j\in D$ used for serving the customers.}
\BlankLine
\tcc{***************************************************************************************************************}
Let $T_0=\emptyset$ denote the truck route and $T_j=\emptyset$ denote the route for drone $j\in D$\;
Find a feasible $k$ and cluster all customer nodes in $C$ into $k$ clusters within distance $L/4$ using a binary search and $k$-means clustering\; 
Let $\bm{m}_1,...,\bm{m}_k$ be the cluster centers and generate their Voronoi tessellation\;
Let $V_1,...,V_k$ be the $k$ Voronoi cells and use them to assign all customer and drone nodes to their corresponding centers within distance $L/4$\;
\For{$\ell=1,...,k$}{
  Let $\bm{z}_{\ell} = \argmin_{\bm{z}_j \in V_{\ell}} \|\bm{z}_j - \bm{m}_{\ell}\|\,, \; \forall j\in C$\;
  Let $\bm{m}_{\ell}^{\prime} = \bm{z}_{\ell}$\;
  Let $C_{\ell}=\{\bm{z}_j \in V_{\ell} \,, \; j\in C\}$\;
  \For{$\bm{z}_j \in C_{\ell} \; \&\& \; \bm{z}_j  \neq \bm{m}_{\ell}^{\prime}$}{
       \eIf{$\bm{z}_j$ can be served via $\bm{m}_{\ell}^{\prime}$ using drones in $V_{\ell}$}{
           \eIf{$|C_{\ell}\backslash \{\bm{m}_{\ell}^{\prime}\}| = 1$ \&\& 
            $2*\|\bm{z}_j  - \bm{m}_{\ell}^{\prime}\| / v_{T}  \leq (\|\bm{y}_i - \bm{m}_{\ell}^{\prime}\| + \|\bm{m}_{\ell}^{\prime} - \bm{z}_j\| + \|\bm{z}_j - \bm{y}_i\|)/v_{D} \,, \; \bm{y}_i \in V_{\ell}\,, \; i\in D$}{
            Let $\bm{z}_j$ be a truck node, i.e., to be served by the truck\;
           }{
           Assign $\bm{z}_j$ to $\bm{m}_{\ell}^{\prime}$, i.e., to be served by a drone\;
           }
           }{
           Let $\bm{z}_j$ be a truck node, i.e., to be served by the truck\;
       }
  }  
}
Update $k$ to be the updated number of truck nodes and reindex them as $\bm{m}_1^{\prime},...,\bm{m}_k^{\prime}$\;
Let $\bm{g}$ be the centroid of $\bm{m}_1^{\prime},...,\bm{m}_k^{\prime}$\;
\For{$\ell=1,...,k$}{
    Set $d_j = \|\bm{z}_j - \bm{g}\|\,, \; \forall j\in C$ and $d_{\bm{m}_{\ell}^{\prime}} = \|\bm{m}_{\ell}^{\prime} - \bm{g}\|$\;
    Let $A_{\ell}$ be an array of customer nodes in $V_{\ell}$ sorted in ascending order of $d_j\,, \; \forall j\in C$\;
    Find the customer node $\bm{z}_{\ell}$ with the smallest index in $A_{\ell}$ and $d_{\bm{z}_{\ell}} \leq d_{\bm{m}_{\ell}^{\prime}}$ that has at least one drone within distance $L/4$ and if chosen as cluster center we can still serve all other customers in $V_{\ell}$ using drones in $V_{\ell}$ and let $\bm{m}_{\ell}^{\prime\prime} = \bm{z}_{\ell}$. If no such point exists let $\bm{m}_{\ell}^{\prime\prime} = \bm{m}_{\ell}^{\prime}$\;
}
Let $\Omega = [\bm{m}_1^{\prime\prime},...,\bm{m}_{k}^{\prime\prime}]$, $\Omega^{\prime} = \emptyset$, and $count=k$\;
\While{$count > 0$}{
     Set $count=0$ \;
     \If{$k = 1$}{
         break\;
     }
     Sort centers $\bm{m}_{\ell}^{\prime\prime} \in \Omega, \ell=1,...,k$ in a descending order of their distances to $\bm{g}$ and reindex them accordingly\;
     Reallocate customer and drone nodes to their nearest centers. For customer nodes pick the closest center that could serve it\;
     \For{$r=1:k$}{
          Let $\bm{\mu} = \argmin_{\bm{m}_{\ell}^{\prime\prime} \in \Omega \backslash \{\Omega^{\prime}\cup \{\bm{m}_{r}^{\prime\prime}\}\}} \|\bm{m}_{\ell}^{\prime\prime} - \bm{m}_{r}^{\prime\prime}\|\,, \; \ell=1,...,k$ and let $s$ be its index in $\Omega$, i.e., $\bm{\mu} = \bm{m}_{s}^{\prime\prime}$. Break a tie with closeness to $\bm{g}$\; 
          Let $feasibility = true$\;
          \For{$\bm{z}_j \in V_r \cup V_s\,, \; j\in C$}{
              Find $\bm{y}_i \in V_{r} \cup V_{s}\,, \; i\in D$ such that $\|\bm{y}_i - \bm{\mu}\| + \|\bm{\mu} - \bm{z}_j\| + \|\bm{z}_j - \bm{y}_i\| \leq L$. If no such drone exists, let $feasibility = false$\;
              \If{$feasibility = false$}{
                 break\;
              }
          }
          \If{$feasibility = true$}{
             \eIf{$k = 2$}{
              Merge the center with fewer customers (say $\bm{m}_{2}^{\prime\prime}$) into the center with more customers (say $\bm{m}_{1}^{\prime\prime}$) and reallocate its assigned customers to that center\;
              Set $\Omega^{\prime} = \Omega^{\prime} \cup \{\bm{m}_{2}^{\prime\prime}\}$\;
              break\;
            }{
             Merge the center $\bm{m}_{r}^{\prime\prime}$ and $V_r$ into $\bm{m}_{s}^{\prime\prime}$ and $V_s$ and allocate its assigned (customer and drone) nodes to $\bm{m}_{s}^{\prime\prime}$\;
             Set $\Omega^{\prime} = \Omega^{\prime} \cup \{\bm{m}_{r}^{\prime\prime}\}$\;
             Set $count = count +1$\;
           }
        }  
     }
     Set $\Omega = \Omega \backslash \Omega^{\prime}$, $k = length(\Omega)$, and $\Omega^{\prime} = \emptyset$\;
}
Run the Tabu Search algorithm to find the drone routes in all Voronoi cells, i.e., $T_j\,, \; j \in D$\;
\If{$k=2$}{
 Set $T_0 = \{\bm{m}_1^{\prime\prime}, \bm{m}_2^{\prime\prime}\}$\;
 }
\If{$k\geq 3$}{
   Run the Lin-Kernighan-Helsgaun algorithm to find the truck route $T_0$ through $\bm{m}_1^{\prime\prime},...,\bm{m}_k^{\prime\prime}$\;
}
\KwRet{$T_0$ \textrm{and} $T_j\,, \; j\in D$}\;
\end{algorithm}
\end{minipage}
} 
\end{center}

\subsection{Utilization of Battery}
\label{sec:batteryUtilization}

In Problem IV, we assumed that drones can revisit the pickup center (truck's stopping place) as long as they do not violate their range (battery duration) limit. Therefore, we need to modify the algorithm to utilize the battery of the drones and allow drones to revisit their assigned pickup center before returning to their base. We assume no limit on the number of times a drone could go back to a pickup center to deliver another package as long as it remains within the range. By implementing this modification, obviously, we may expect some savings since if one drone can deliver another package without going back to the base, it can save time on that trip. However, this adds to the waiting time of the truck and makes the truck route costlier. This trade-off is analyzed in Section \ref{sec:batteryUtilizationImpact}. We apply this modification to all three Algorithms \ref{alg:RouteFinderMC1}, \ref{alg:RouteFinderMC2}, and \ref{alg:RouteFinderMC3}. 
We skip presenting the pseudocode here for brevity (the only difference is allowing revisiting in range conditions and in the TS step). The modification in the tabu search step is applied by changing the solution representation. Our original example of a solution in a cluster for a problem with 5 customer nodes and 4 drone nodes can change as following 
\begin{center}
          \begin{tabular}{|c|c|c|c|c|c|c|c|c|c|c|c|c|c|c|}
          \hline
          8 & 1 & 5& 7& 1& 2& \textcolor{blue}{9} & \textcolor{blue}{1} & \textcolor{blue}{3} & \textcolor{blue}{1} & \textcolor{blue}{6} & 10& 1& 4\\
          \hline
          \end{tabular}
\end{center}
In this solution, drone from location 9 first serves the customer node 3 and then revisits the truck node $1$ (without going to its base first) to serve the customer node 6.

\section{Computational Results}
\label{sec:ComputationalResults}
As problems III and IV are the generalized versions of Problems I and II, we only do our computational results for Problems III and IV. We will see that for certain values of input parameters we actually get solutions for the special cases of Problems I and II. We tested the MILP model of Problems III and IV on a set of two synthetic problem instances with 16 and 12 ($P1$) and 16 and 16 ($P2$) customer and drone nodes, respectively, distributed uniformly at random in a unit square. The speed of drones is set to be 2 times the speed of the truck in both models. For the drone flying range $L$ we took the maximum length of a drone in a unit box, i.e., $2+\sqrt{2} \simeq 3.4$, and with increments of 0.2 we generated 17 instance for each problem. In the tabu search step of the algorithms we have set |Tabu List| = 0.7 * |Action List| with a termination condition of 1000 iterations. We have also set the precision to be $\epsilon=10^{-5}$ for all solutions.

We implemented our models in Julia 1.9.2 using Gurobi Optimizer 10.0.2 on  
an Intel Core i7-5600U @2.6GHz, 8GB DDR3 RAM computer. We ran our algorithms for the same problem instances on the same computer using MATLAB R2022b. 
 We also ran our algorithms on larger problem instances for comparison with the results of traditional centralized delivery system (pure TSP model) and the coordinated truck and drone delivery (TSP-D model).  

\subsection{Results for Multiple Centers with Recharging}
\label{subsec:resultsProblemIII}
\begin{figure}[p]
    \centering
    \begin{subfigure}[b]{0.16\linewidth}
        \includegraphics[width=\linewidth]{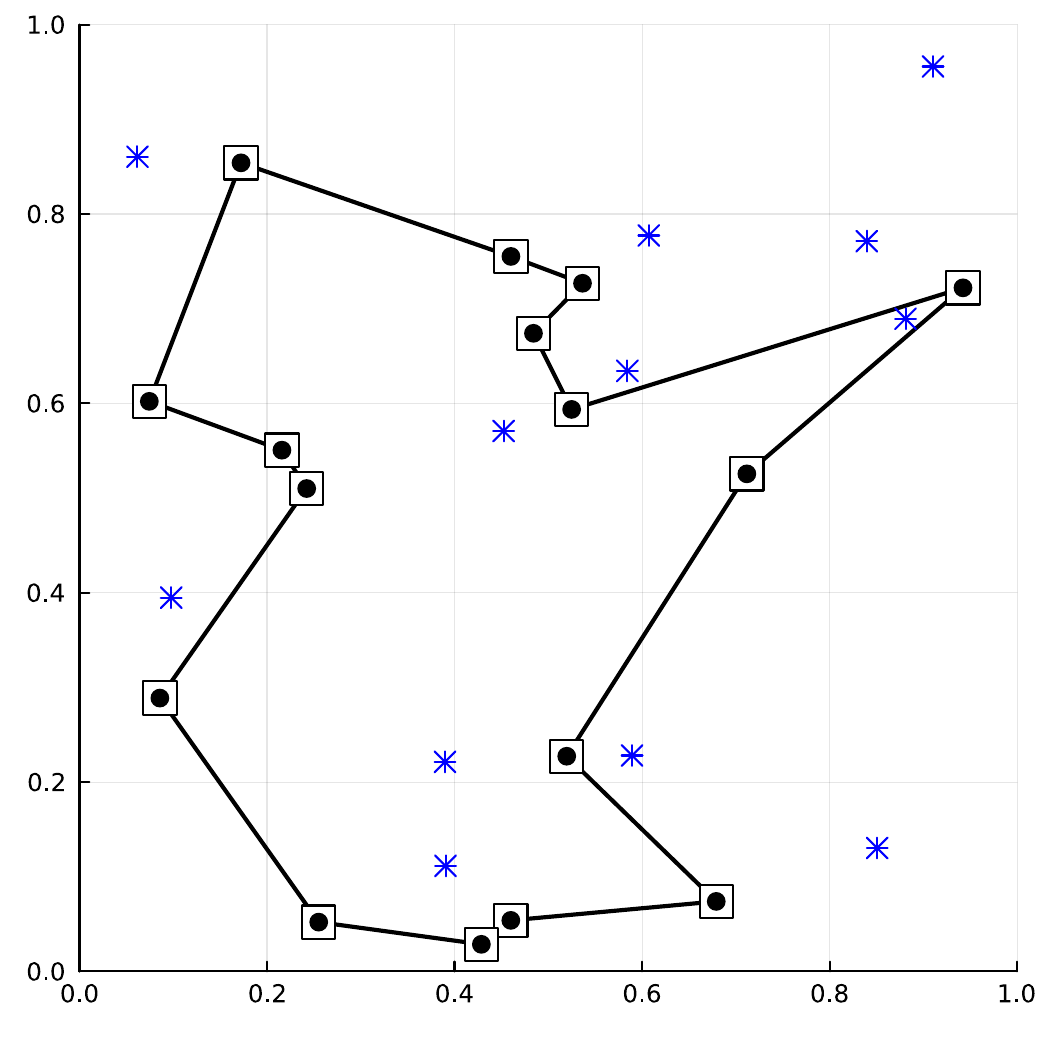}
        \caption{$P1,\;L = 0.2$}
    \end{subfigure}
    \begin{subfigure}[b]{0.16\linewidth}
        \includegraphics[width=\linewidth]{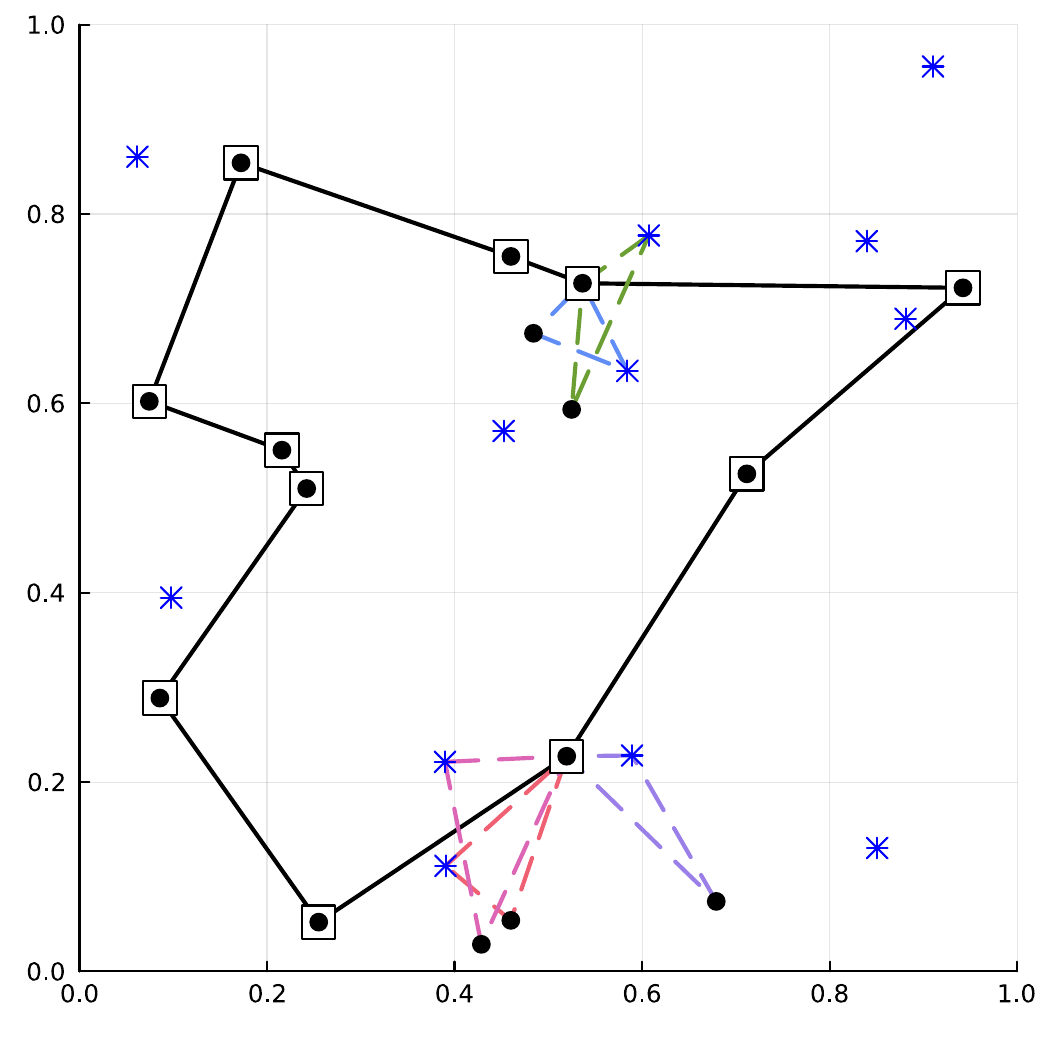}
        \caption{$P1,\;L = 0.6$}
    \end{subfigure}
    \begin{subfigure}[b]{0.16\linewidth}
        \includegraphics[width=\linewidth]{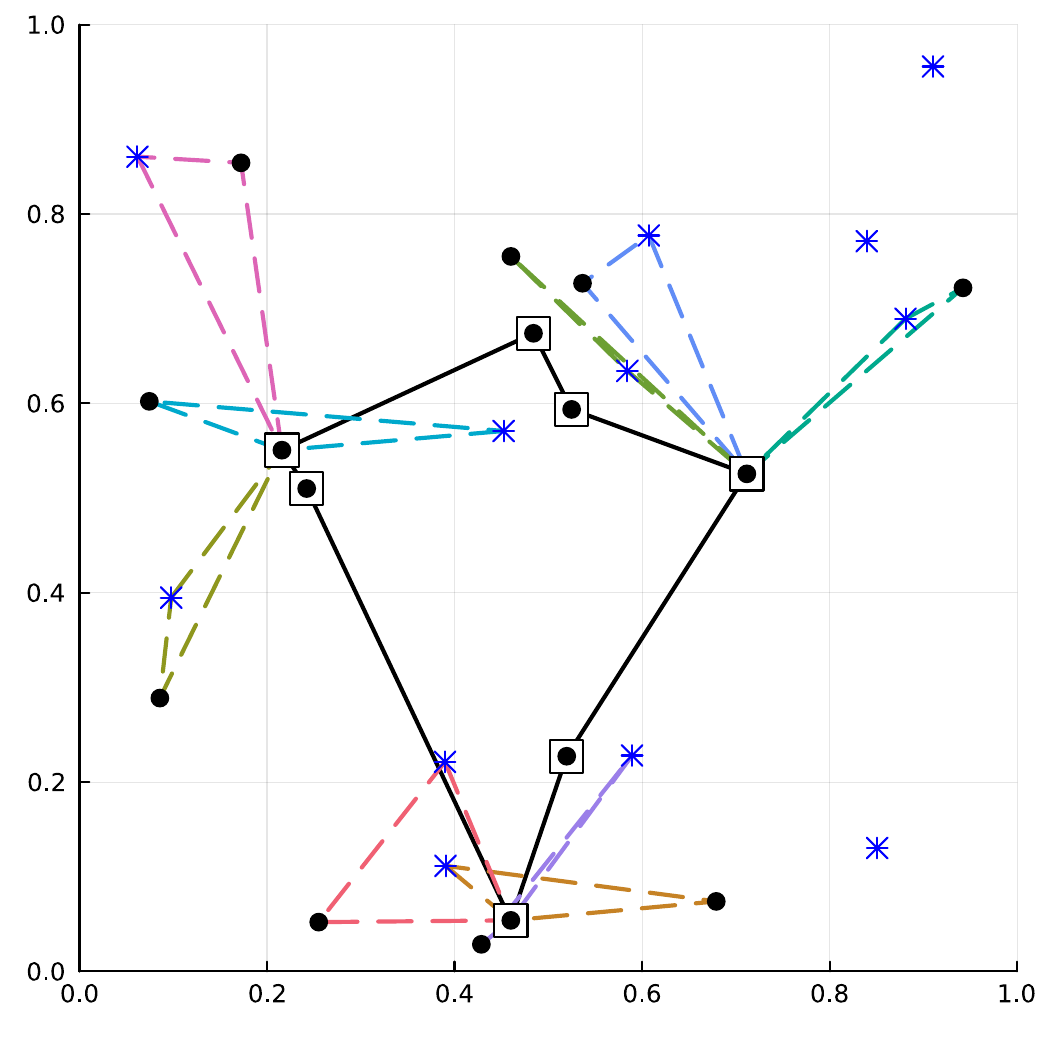}
        \caption{$P1,\;L = 0.8$}
    \end{subfigure} 
    \quad
        \begin{subfigure}[b]{0.16\linewidth}
        \includegraphics[width=\linewidth]{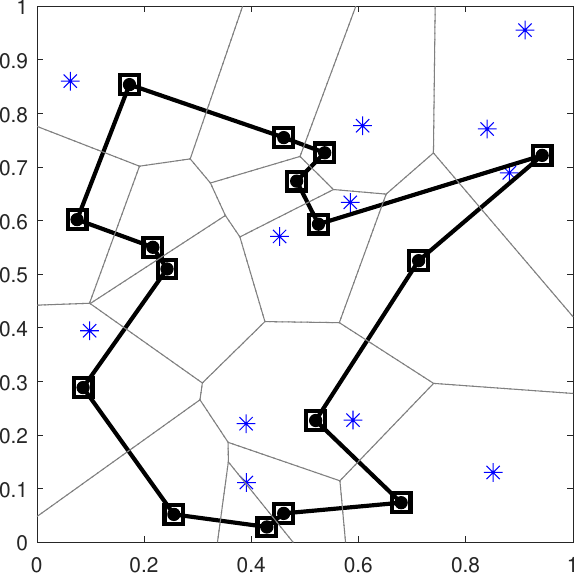}
        \caption{$P1,\;L = 0.2$}
    \end{subfigure}
    \begin{subfigure}[b]{0.16\linewidth}
        \includegraphics[width=\linewidth]{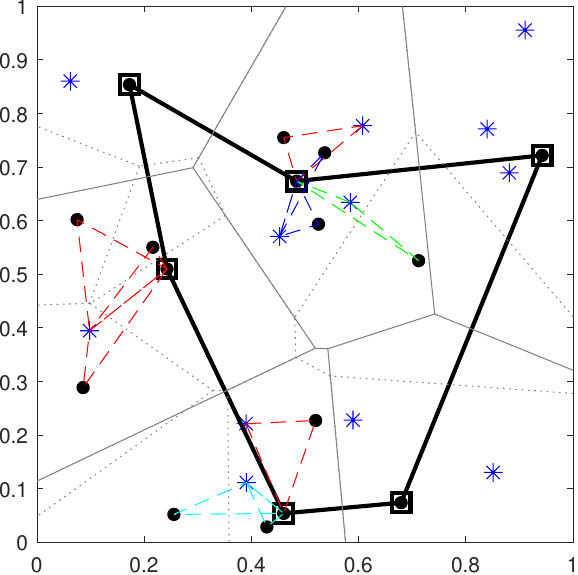}
        \caption{$P1,\;L = 0.6$}
    \end{subfigure}
    \begin{subfigure}[b]{0.16\linewidth}
        \includegraphics[width=\linewidth]{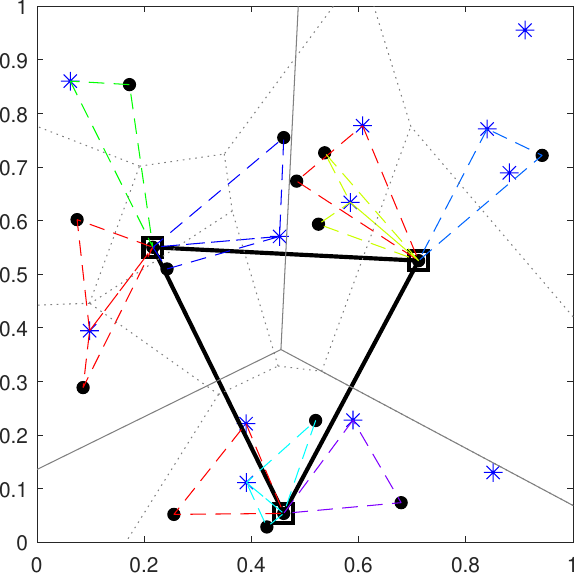}
        \caption{$P1,\;L = 0.8$}
    \end{subfigure} 
    \\
    \begin{subfigure}[b]{0.16\linewidth}
        \includegraphics[width=\linewidth]{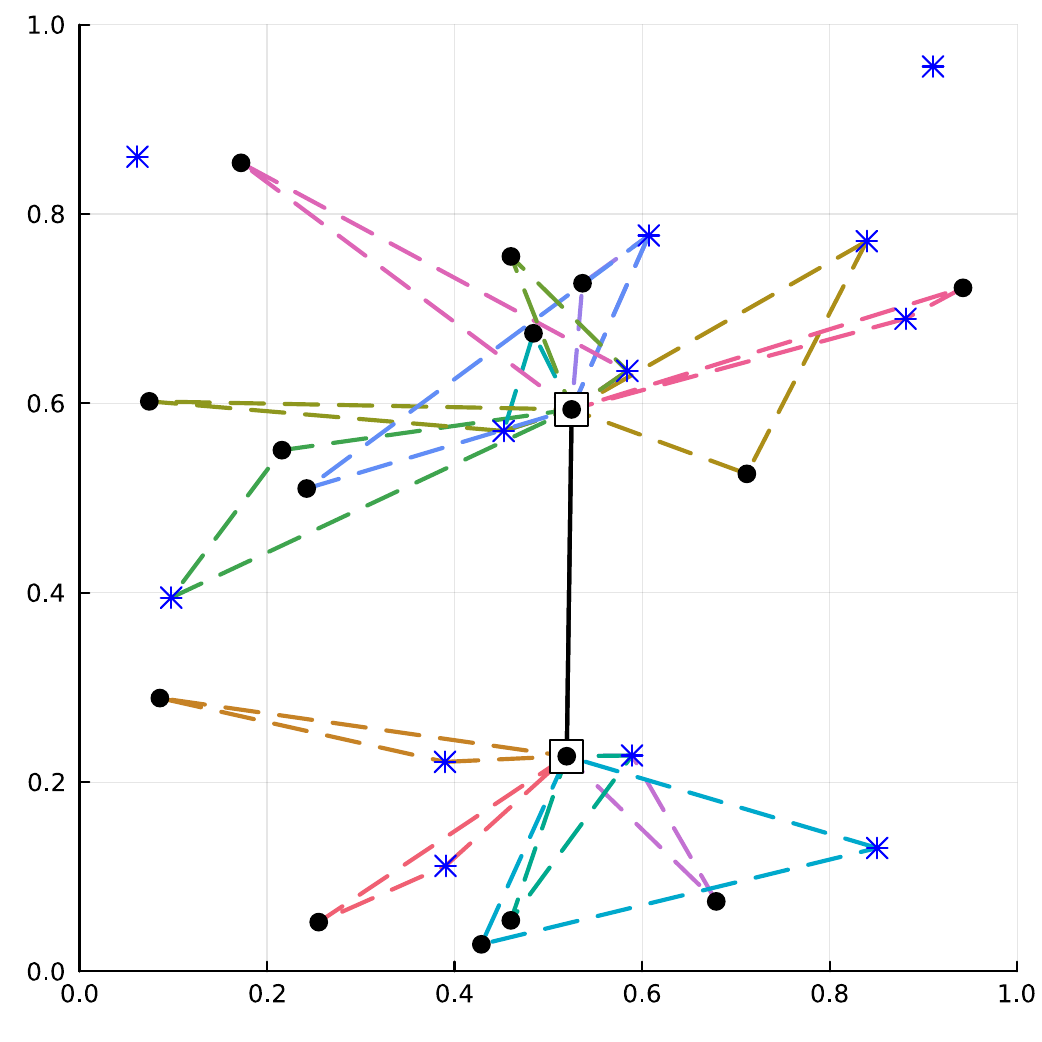}
        \caption{$P1,\;L = 1.0$}
    \end{subfigure}
    \begin{subfigure}[b]{0.16\linewidth}
        \includegraphics[width=\linewidth]{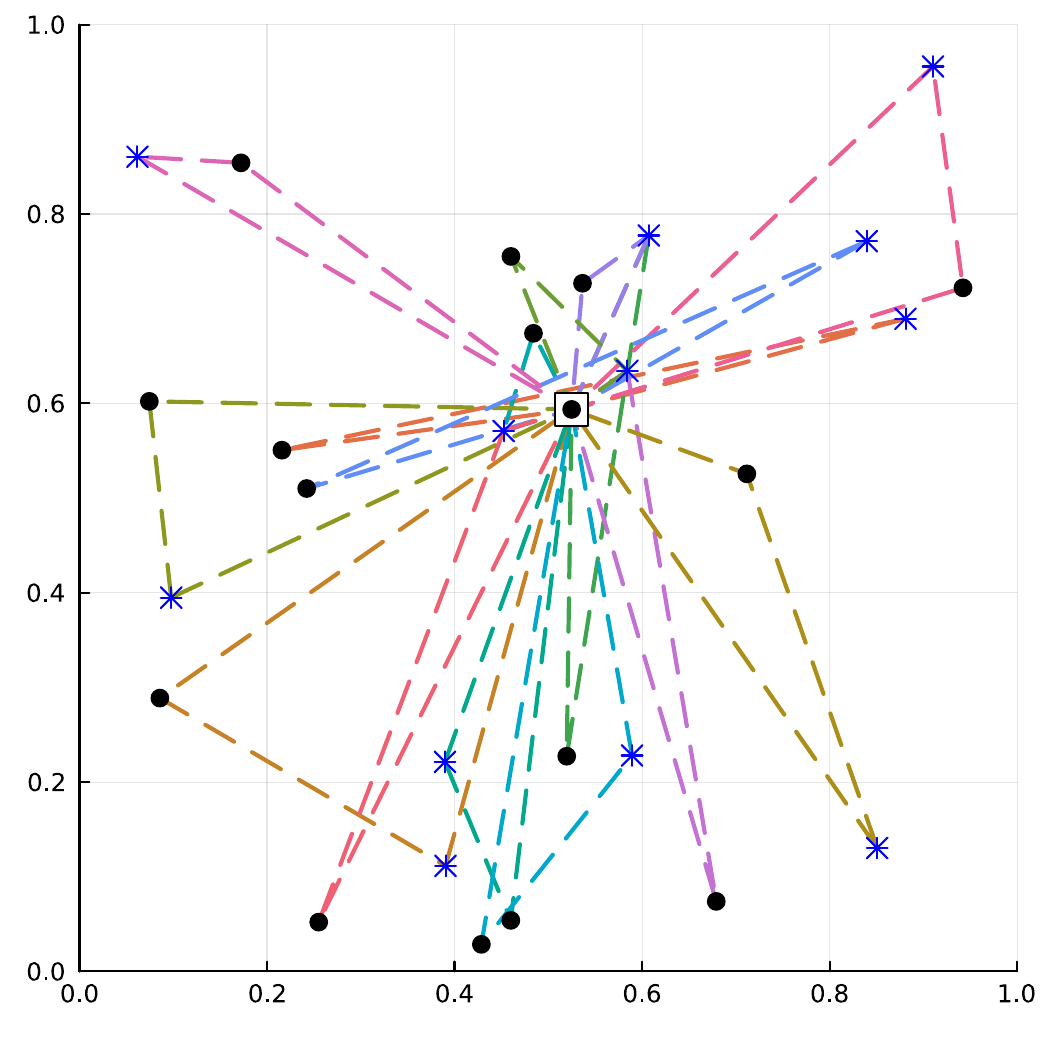}
        \caption{$P1,\;L = 1.4$}
    \end{subfigure}
    \begin{subfigure}[b]{0.16\linewidth}
        \includegraphics[width=\linewidth]{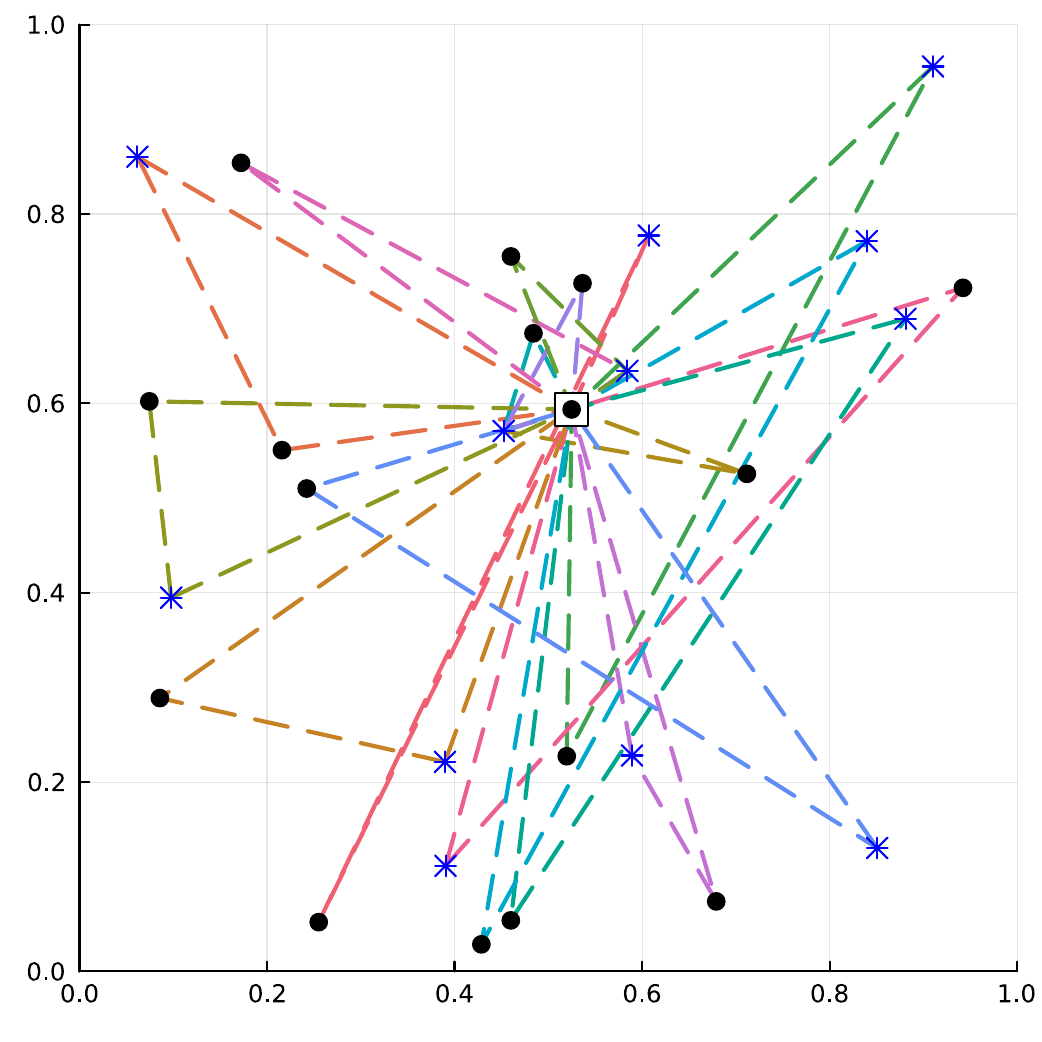}
        \caption{$P1,\;L = 1.8$}
    \end{subfigure}
        \quad
        \begin{subfigure}[b]{0.16\linewidth}
        \includegraphics[width=\linewidth]{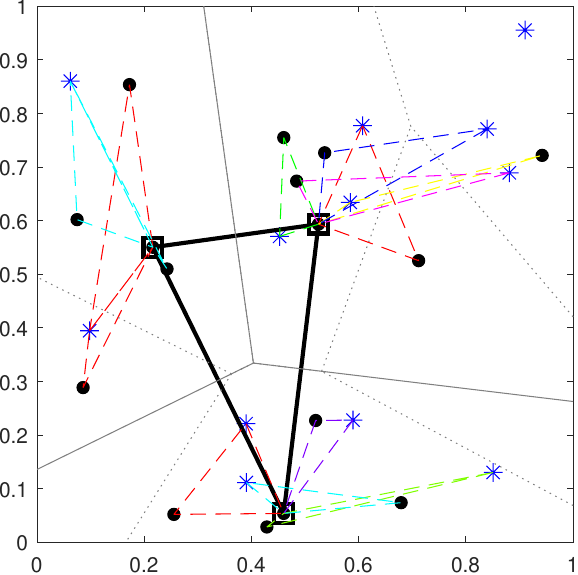}
        \caption{$P1,\;L = 1.0$}
    \end{subfigure}
    \begin{subfigure}[b]{0.16\linewidth}
        \includegraphics[width=\linewidth]{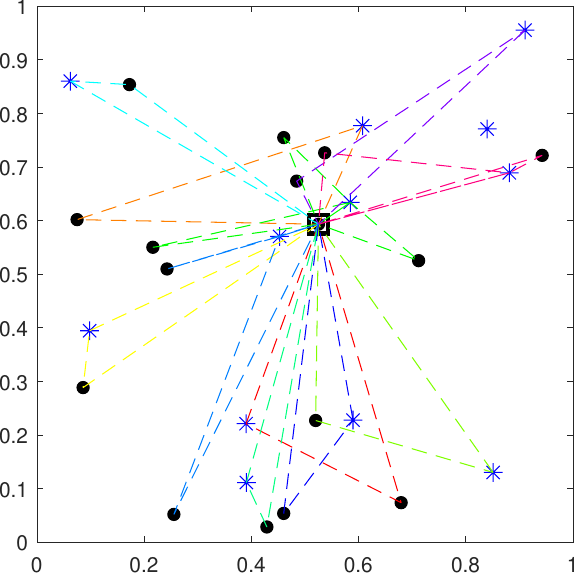}
        \caption{$P1,\;L = 1.4$}
    \end{subfigure}
    \begin{subfigure}[b]{0.16\linewidth}
        \includegraphics[width=\linewidth]{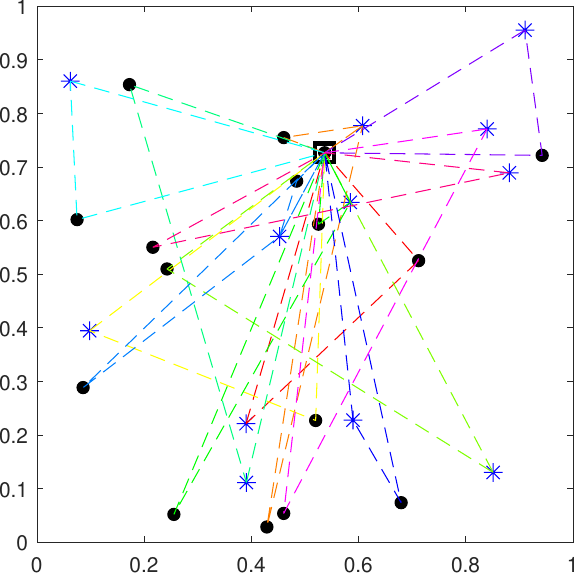}
        \caption{$P1,\;L = 1.8$}
    \end{subfigure} 
    \\
        \begin{subfigure}[b]{0.16\linewidth}
        \includegraphics[width=\linewidth]{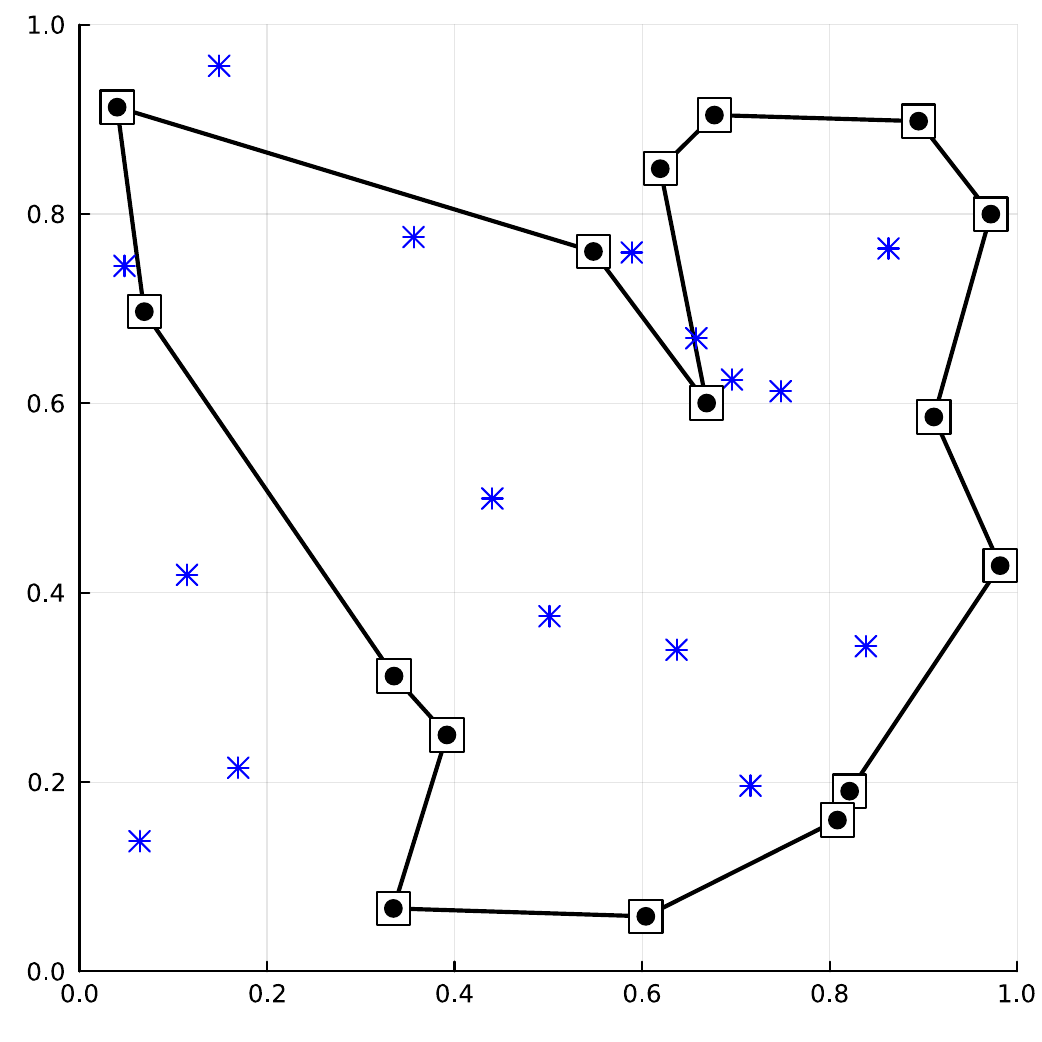}
        \caption{$P2,\;L = 0.2$}
    \end{subfigure}
    \begin{subfigure}[b]{0.16\linewidth}
        \includegraphics[width=\linewidth]{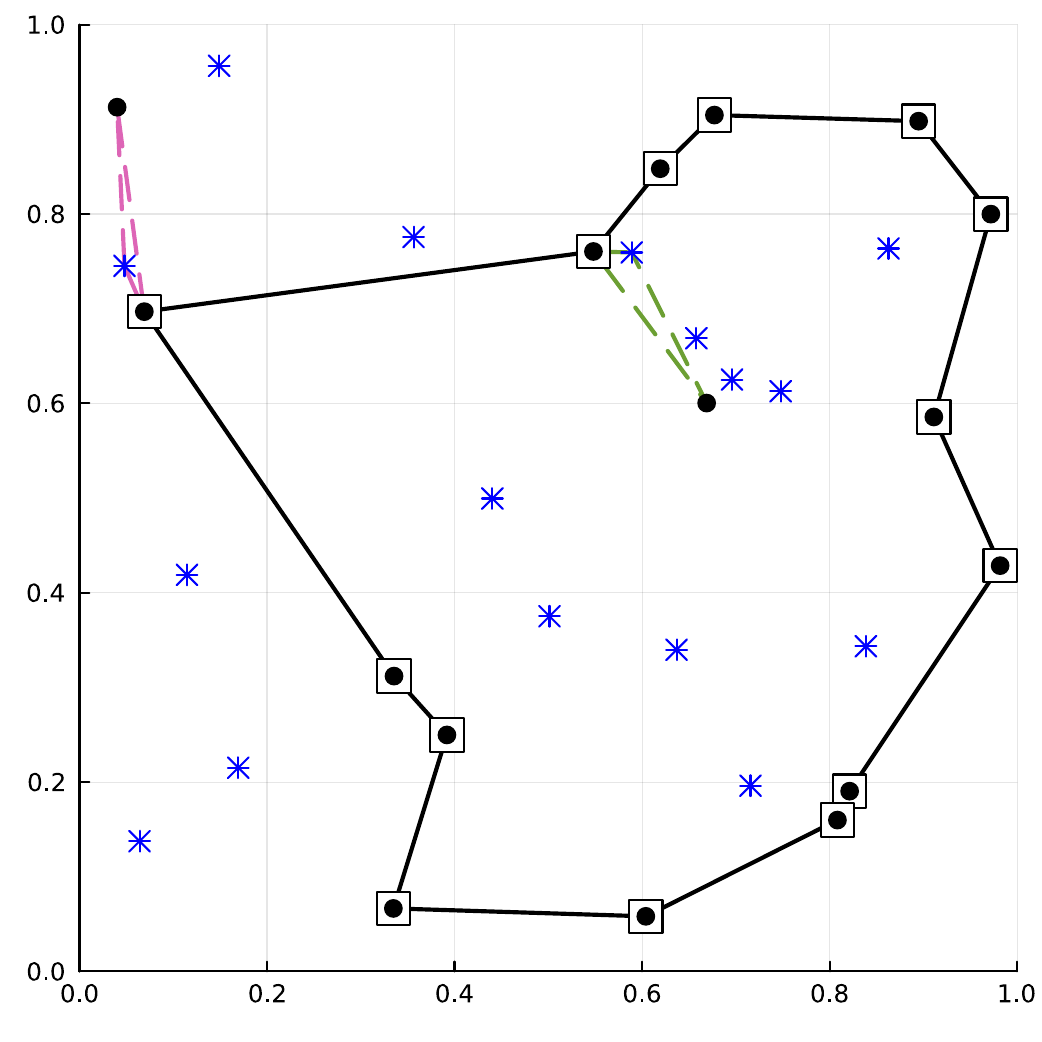}
        \caption{$P2,\;L = 0.6$}
    \end{subfigure}
    \begin{subfigure}[b]{0.16\linewidth}
        \includegraphics[width=\linewidth]{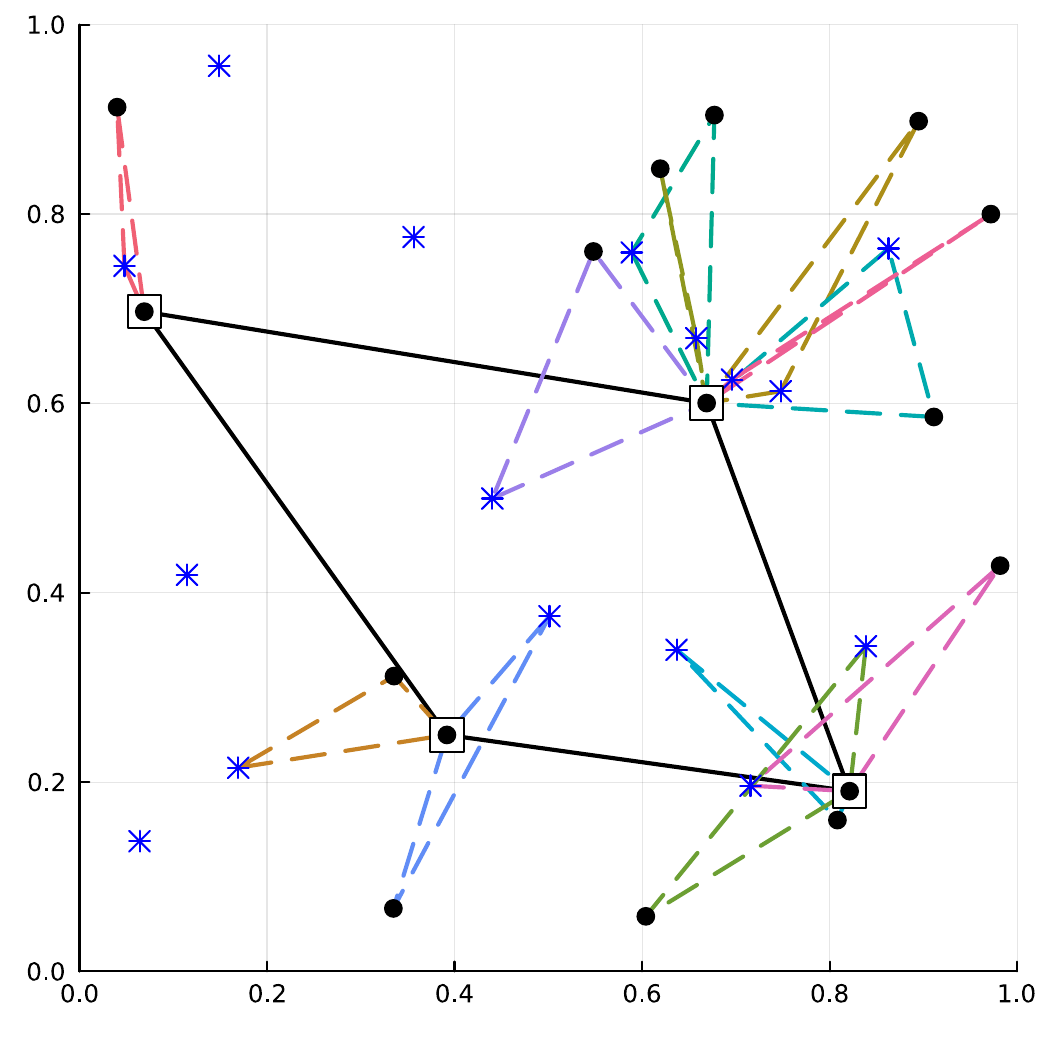}
        \caption{$P2,\;L = 0.8$}
    \end{subfigure} 
     \quad
        \begin{subfigure}[b]{0.16\linewidth}
        \includegraphics[width=\linewidth]{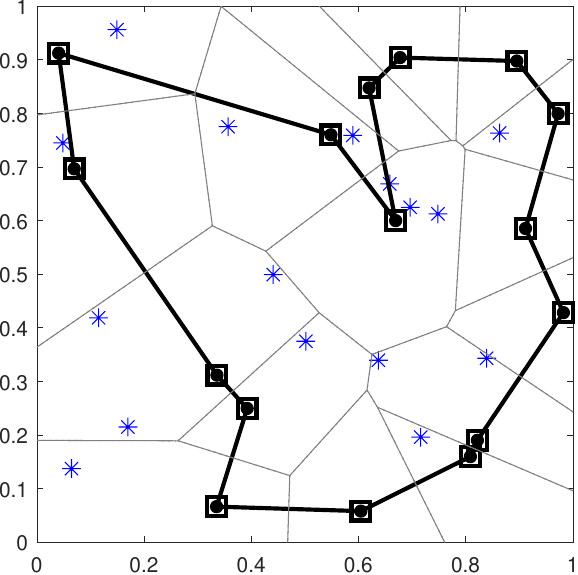}
        \caption{$P2,\;L = 0.2$}
    \end{subfigure}
    \begin{subfigure}[b]{0.16\linewidth}
        \includegraphics[width=\linewidth]{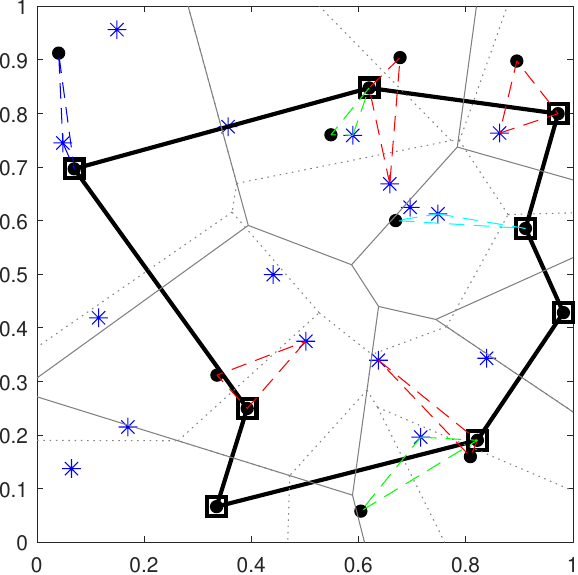}
        \caption{$P2,\;L = 0.6$}
    \end{subfigure}
    \begin{subfigure}[b]{0.16\linewidth}
        \includegraphics[width=\linewidth]{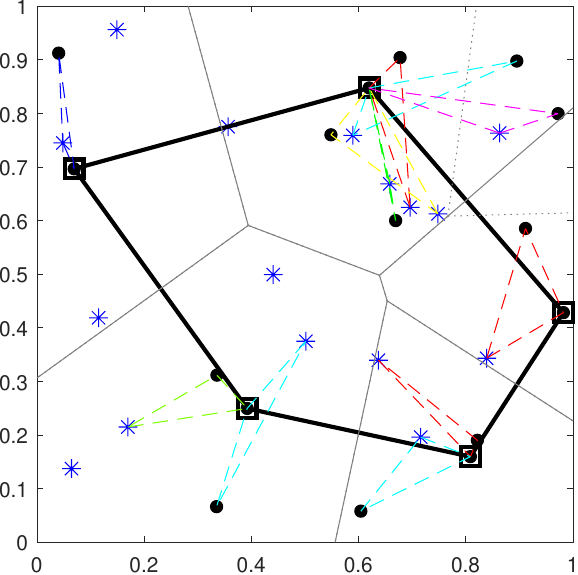}
        \caption{$P2,\;L = 0.8$}
    \end{subfigure} 
    \\
    \begin{subfigure}[b]{0.16\linewidth}
        \includegraphics[width=\linewidth]{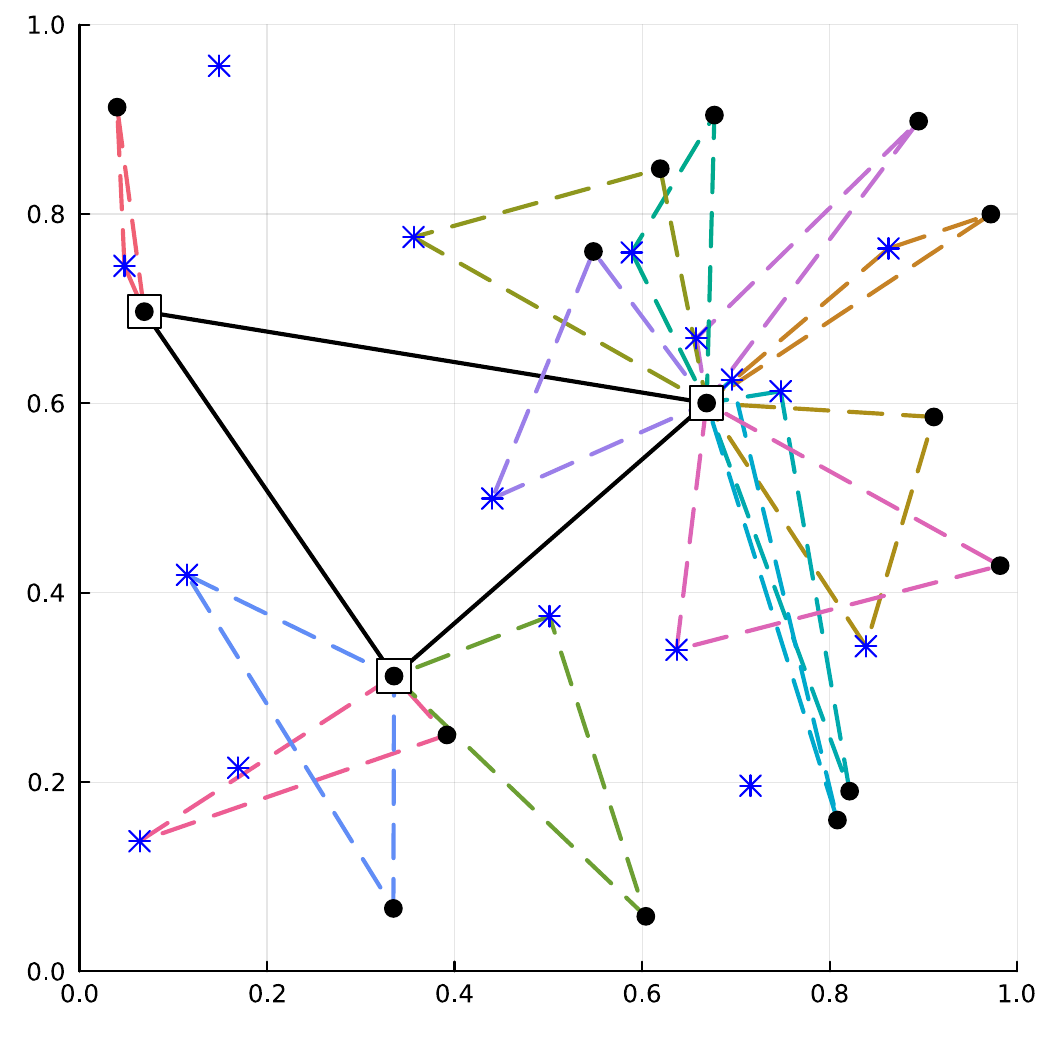}
        \caption{$P2,\;L = 1.0$}
    \end{subfigure}
    \begin{subfigure}[b]{0.16\linewidth}
        \includegraphics[width=\linewidth]{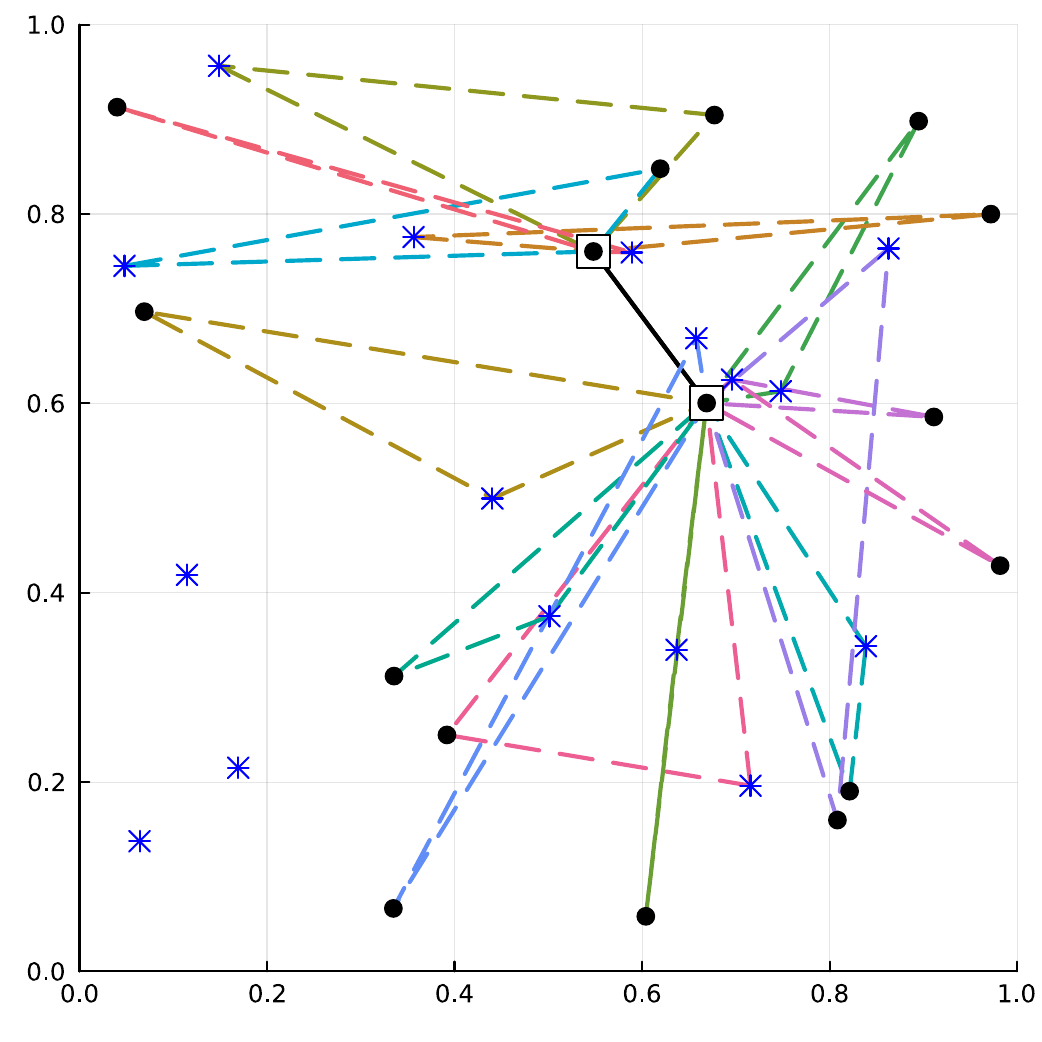}
        \caption{$P2,\;L = 1.4$}
    \end{subfigure}
    \begin{subfigure}[b]{0.16\linewidth}
        \includegraphics[width=\linewidth]{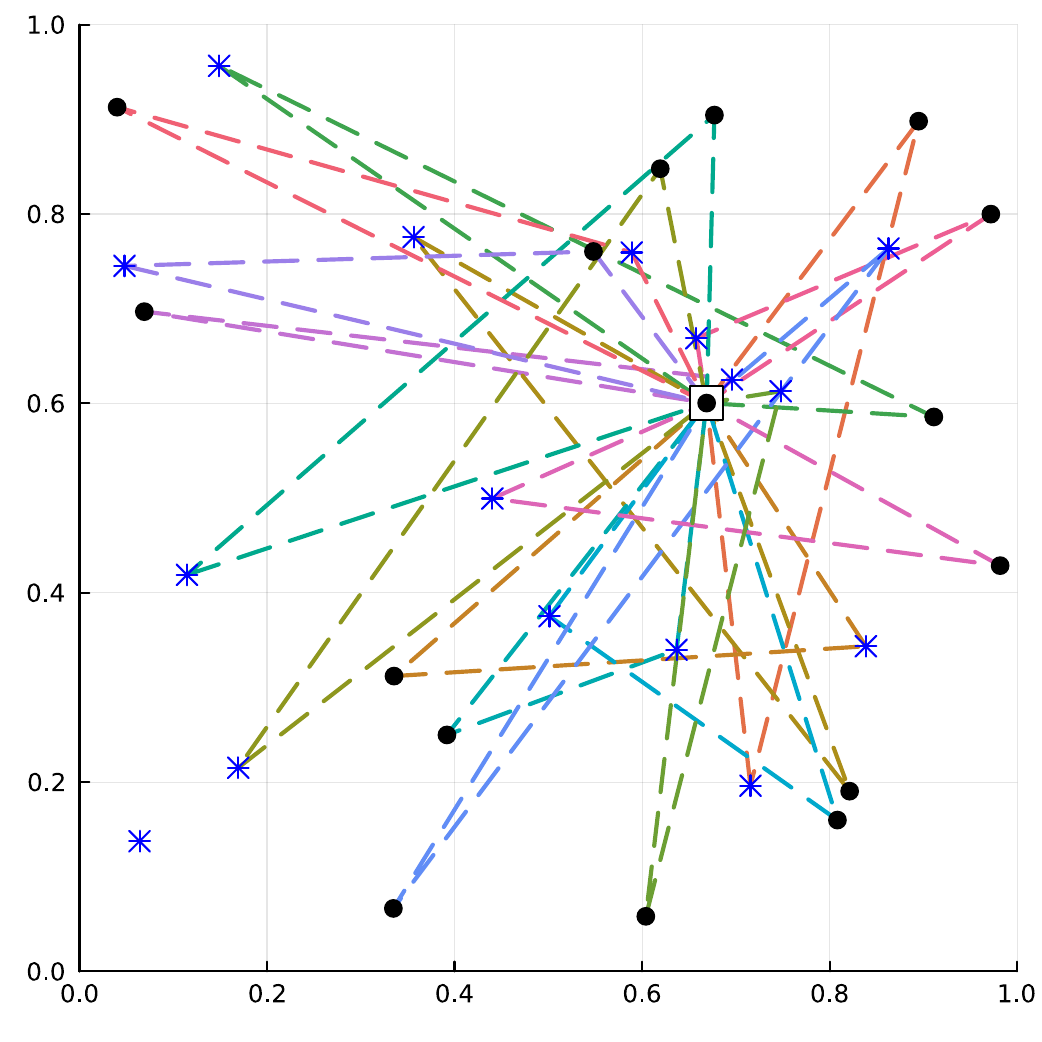}
        \caption{$P2,\;L = 1.8$}
    \end{subfigure}
            \quad
        \begin{subfigure}[b]{0.16\linewidth}
        \includegraphics[width=\linewidth]{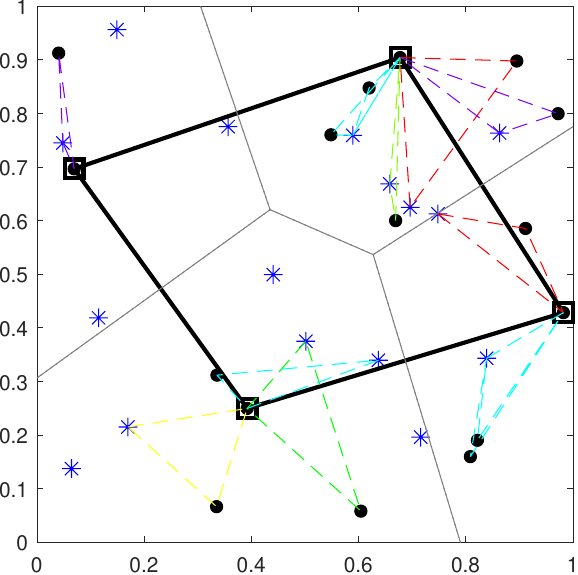}
        \caption{$P2,\;L = 1.0$}
    \end{subfigure}
    \begin{subfigure}[b]{0.16\linewidth}
        \includegraphics[width=\linewidth]{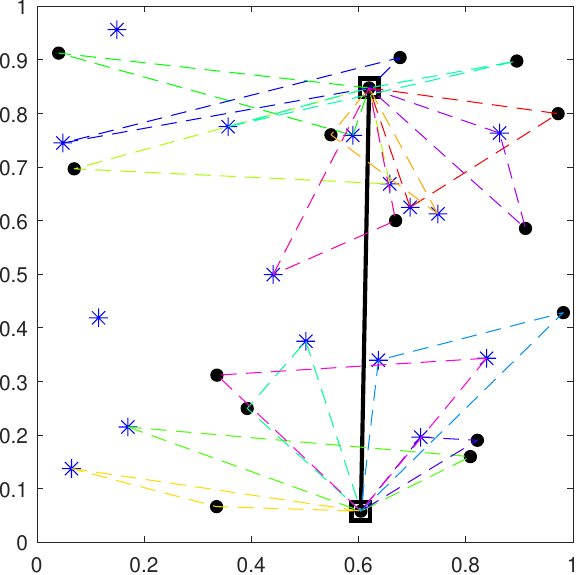}
        \caption{$P2,\;L = 1.4$}
    \end{subfigure}
    \begin{subfigure}[b]{0.16\linewidth}
        \includegraphics[width=\linewidth]{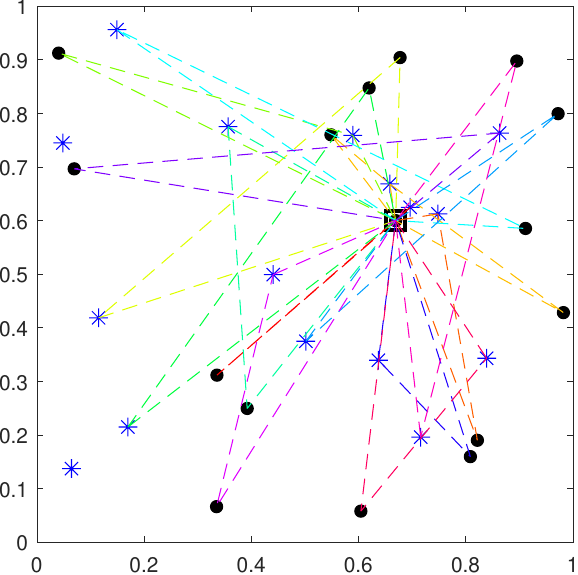}
        \caption{$P2,\;L = 1.8$}
    \end{subfigure} 
    \\
    \begin{subfigure}[b]{0.3\linewidth}
       \includegraphics[width=\linewidth]{legend.pdf}
    \end{subfigure}
    \caption{Examples of optimal solutions of Problem III (multiple centers with recharging) obtained by the optimization model for problems $P1$ with 16 customer nodes and 12 drone nodes and $P2$ with 16 customer nodes and 16 drone nodes for select values of drone range $L$ are provided on the left three columns. Their corresponding (optimal/suboptimal) solutions obtained by Algorithm \ref{alg:RouteFinderMC3} is shown on the three columns on the right. 
    It can be seen that drones are not necessarily assigned to their closest center or closest customer. This is due to the fact that the truck is not required to wait at a center while its assigned drones are traversing the last edge of their route (after delivery) to go back to their base.
    For small drone range such as $L = 0.2$ (with a 2-to-1 relative drone speed assumption) we are better off to not use the drones and for very large drone ranges such as $L\geq 1.8$ the problem reduces to a one center (drone-only) delivery problem.}
    \label{fig:resultsModelRecharging}
\end{figure}

Table \ref{tab:RechargingResultsModelAlgComp} compares the computational results of the our optimization model for Problem III and the proposed algorithms on the 34 instances of $P1$ and $P2$. As it is evident from the table, our algorithms are much faster than the optimization model. The average running time of the final version of our algorithm (Algorithm \ref{alg:RouteFinderMC3}) is 8.94 seconds versus 885.13 seconds of the optimization model for $P1$ and 9.56 seconds versus 172.23 seconds for $P2$. The average (maximum) optimality gap of the solutions provided by Algorithm \ref{alg:RouteFinderMC3} for $P1$ and $P2$ is 11.92\% (57.66\%) and 16.61\% (102.71\%), respectively. 
Figure \ref{fig:resultsModelRecharging} illustrates optimal solutions obtained by our optimization model and optimal/suboptimal solutions obtained by Algorithm \ref{alg:RouteFinderMC3} for $P1$ and $P2$ for select instances of $L$.

\begin{table}[htbp]
\scriptsize
  \centering
  \caption{\label{tab:RechargingResultsModelAlgComp} Comparison between the results of the optimization model for Problem III (Multiple Centers with Recharging) and the algorithms on two synthetic examples with 16 and 12 ($P1$) and 16 and 16 ($P2$) customers nodes and drone nodes, respectively, distributed randomly in a unit box and for different values of $L$. The column ``Time (s)'' shows the computational time in seconds. The column ``Gap'' presents the optimality gap of the solutions found by the algorithms. The instances for which the optimal solution is found by our algorithms are shown in bold (\textbf{0.00\% gap}). \\}
    \begin{tabular}{cccccccccccc}
    \bottomrule \hline
          & \multicolumn{11}{c}{\textbf{P1 (16-12) -- Recharging}} \\
\cmidrule{2-12}    \multirow{2}[4]{*}{\textbf{L }} & \multicolumn{2}{c}{\textbf{Opt. Model}} & \multicolumn{3}{c}{\textbf{Algorithm 1}} & \multicolumn{3}{c}{\textbf{Algorithm 2}} & \multicolumn{3}{c}{\textbf{Algorithm 3}} \\
\cmidrule{2-12}          & Time (s) & Obj. Value  & Time (s) & Obj. Value  & Gap (\%)  & Time (s) & Obj. Value  & Gap (\%)  & Time (s) & Obj. Value  & Gap (\%)  \\
    \midrule
    \multicolumn{1}{c|}{0.2} & 2.80  & 3.3325 & 0.36  & 3.3325 & \textbf{0.00\%} & 0.36  & 3.3325 & \textbf{0.00\%} & 0.42  & 3.3325 & \textbf{0.00\%} \\
    \multicolumn{1}{c|}{0.4} & 21.93 & 3.2561 & 0.63  & 3.2837 & 0.85\% & 0.61  & 3.2837 & 0.85\% & 0.62  & 3.3385 & 2.53\% \\
    \multicolumn{1}{c|}{0.6} & 384.27 & 3.0865 & 0.70  & 3.2659 & 5.81\% & 0.67  & 3.2868 & 6.49\% & 1.12  & 3.7672 & 22.05\% \\
    \multicolumn{1}{c|}{0.8} & 7772.69 & 2.4655 & 0.89  & 3.2056 & 30.02\% & 0.81  & 3.2881 & 33.36\% & 1.65  & 2.7818 & 12.83\% \\
    \multicolumn{1}{c|}{1} & 483.70 & 1.5588 & 0.96  & 3.6451 & 133.85\% & 0.92  & 3.6421 & 133.66\% & 1.62  & 2.4575 & 57.66\% \\
    \multicolumn{1}{c|}{1.2} & 5445.16 & 1.4339 & 1.14  & 2.8359 & 97.78\% & 1.14  & 2.8080 & 95.84\% & 3.56  & 1.7979 & 25.39\% \\
    \multicolumn{1}{c|}{1.4} & 11.45 & 0.5176 & 1.19  & 3.1534 & 509.28\% & 1.35  & 2.9181 & 463.82\% & 13.04 & 0.6887 & 33.07\% \\
    \multicolumn{1}{c|}{1.6} & 97.43 & 0.4828 & 1.05  & 3.1534 & 553.13\% & 1.17  & 2.9218 & 505.17\% & 12.88 & 0.5029 & 4.16\% \\
    \multicolumn{1}{c|}{1.8} & 43.86 & 0.4686 & 3.12  & 1.7029 & 263.40\% & 3.18  & 1.7107 & 265.07\% & 13.28 & 0.5849 & 24.82\% \\
    \multicolumn{1}{c|}{2} & 100.87 & 0.4686 & 3.03  & 1.7029 & 263.40\% & 3.06  & 1.4503 & 209.51\% & 12.84 & 0.4772 & 1.84\% \\
    \multicolumn{1}{c|}{2.2} & 92.24 & 0.4686 & 2.96  & 1.7029 & 263.40\% & 2.95  & 1.4503 & 209.51\% & 12.83 & 0.4869 & 3.90\% \\
    \multicolumn{1}{c|}{2.4} & 97.30 & 0.4686 & 12.63 & 0.4737 & 1.09\% & 12.82 & 0.4869 & 3.90\% & 13.11 & 0.4832 & 3.12\% \\
    \multicolumn{1}{c|}{2.6} & 100.90 & 0.4686 & 12.64 & 0.4694 & 0.18\% & 12.62 & 0.4755 & 1.47\% & 14.06 & 0.4734 & 1.03\% \\
    \multicolumn{1}{c|}{2.8} & 96.93 & 0.4686 & 12.66 & 0.4836 & 3.21\% & 12.85 & 0.4832 & 3.12\% & 12.69 & 0.4775 & 1.91\% \\
    \multicolumn{1}{c|}{3} & 98.45 & 0.4686 & 12.91 & 0.4832 & 3.12\% & 12.96 & 0.5029 & 7.32\% & 12.86 & 0.4828 & 3.03\% \\
    \multicolumn{1}{c|}{3.2} & 98.74 & 0.4686 & 12.78 & 0.4836 & 3.21\% & 12.71 & 0.4734 & 1.03\% & 12.86 & 0.4846 & 3.41\% \\
    \multicolumn{1}{c|}{3.4} & 98.49 & 0.4686 & 12.78 & 0.4755 & 1.47\% & 13.15 & 0.4832 & 3.12\% & 12.59 & 0.4775 & 1.91\% \\
          &       &       &       &       &       &       &       &       &       &       &  \\
	\rowcolor{mygray}   \textbf{Average} & 885.13 & -     & 5.44  & -     & 125.48\% & 5.49  & -     & 114.31\% & 8.94  & -     & 11.92\% \\
    \midrule
          &       &       &       &       &       &       &       &       &       &       &  \\
          &       &       &       &       &       &       &       &       &       &       &  \\
    \midrule
          & \multicolumn{11}{c}{\textbf{P2 (16-16) -- Recharging}} \\
\cmidrule{2-12}    \multirow{2}[4]{*}{\textbf{L }} & \multicolumn{2}{c}{\textbf{Opt. Model}} & \multicolumn{3}{c}{\textbf{Algorithm 1}} & \multicolumn{3}{c}{\textbf{Algorithm 2}} & \multicolumn{3}{c}{\textbf{Algorithm 3}} \\
\cmidrule{2-12}          & Time (s) & Obj. Value  & Time (s) & Obj. Value  & Gap (\%)  & Time (s) & Obj. Value  & Gap (\%)  & Time (s) & Obj. Value  & Gap (\%)  \\
    \midrule
    \multicolumn{1}{c|}{0.2} & 2.63  & 3.5817 & 0.55  & 3.5817 & \textbf{0.00\%} & 0.44  & 3.5817 & \textbf{0.00\%} & 0.47  & 3.5817 & \textbf{0.00\%} \\
    \multicolumn{1}{c|}{0.4} & 82.08 & 3.5817 & 0.46  & 3.6717 & 2.51\% & 0.40  & 3.6717 & 2.51\% & 0.49  & 3.7663 & 5.15\% \\
    \multicolumn{1}{c|}{0.6} & 19.84 & 3.2325 & 0.35  & 3.4684 & 7.30\% & 0.33  & 3.7302 & 15.40\% & 0.55  & 3.7508 & 16.03\% \\
    \multicolumn{1}{c|}{0.8} & 90.52 & 2.7967 & 0.90  & 3.5360 & 26.43\% & 0.66  & 3.5360 & 26.43\% & 1.11  & 3.3767 & 20.74\% \\
    \multicolumn{1}{c|}{1} & 209.00 & 2.2349 & 1.25  & 3.6115 & 61.59\% & 1.22  & 3.6115 & 61.59\% & 1.22  & 3.6115 & 61.59\% \\
    \multicolumn{1}{c|}{1.2} & 631.79 & 1.7369 & 1.28  & 3.3403 & 92.31\% & 1.31  & 3.3403 & 92.31\% & 3.16  & 2.5199 & 45.08\% \\
    \multicolumn{1}{c|}{1.4} & 1645.00 & 1.1662 & 1.17  & 3.3403 & 186.43\% & 1.20  & 3.3403 & 186.43\% & 3.26  & 2.3639 & 102.71\% \\
    \multicolumn{1}{c|}{1.6} & 69.15 & 0.4673 & 1.38  & 3.2021 & 585.27\% & 1.36  & 2.9233 & 525.62\% & 15.25 & 0.4673 & \textbf{0.00\%} \\
    \multicolumn{1}{c|}{1.8} & 16.60 & 0.4438 & 1.19  & 3.5377 & 697.22\% & 1.25  & 3.3313 & 650.72\% & 15.41 & 0.4438 & \textbf{0.00\%} \\
    \multicolumn{1}{c|}{2} & 6.51  & 0.4438 & 1.13  & 3.0457 & 586.34\% & 1.18  & 3.0333 & 583.56\% & 15.05 & 0.5803 & 30.76\% \\
    \multicolumn{1}{c|}{2.2} & 37.22 & 0.4438 & 3.78  & 1.6636 & 274.90\% & 3.44  & 1.6636 & 274.90\% & 15.66 & 0.4438 & \textbf{0.00\%} \\
    \multicolumn{1}{c|}{2.4} & 21.44 & 0.4438 & 3.56  & 1.6636 & 274.90\% & 3.40  & 1.6636 & 274.90\% & 15.02 & 0.4438 & \textbf{0.00\%} \\
    \multicolumn{1}{c|}{2.6} & 13.89 & 0.4438 & 3.41  & 1.6636 & 274.90\% & 3.43  & 1.6636 & 274.90\% & 14.99 & 0.4438 & \textbf{0.00\%} \\
    \multicolumn{1}{c|}{2.8} & 11.76 & 0.4438 & 14.77 & 0.4438 & \textbf{0.00\%} & 15.72 & 0.4438 & \textbf{0.00\%} & 15.50 & 0.4450 & 0.28\% \\
    \multicolumn{1}{c|}{3} & 23.71 & 0.4438 & 15.00 & 0.4438 & \textbf{0.00\%} & 14.89 & 0.4450 & 0.28\% & 15.15 & 0.4438 & \textbf{0.00\%} \\
    \multicolumn{1}{c|}{3.2} & 23.79 & 0.4438 & 14.71 & 0.4598 & 3.61\% & 14.75 & 0.4438 & \textbf{0.00\%} & 15.26 & 0.4438 & \textbf{0.00\%} \\
    \multicolumn{1}{c|}{3.4} & 22.98 & 0.4438 & 14.87 & 0.4438 & \textbf{0.00\%} & 14.88 & 0.4438 & \textbf{0.00\%} & 14.94 & 0.4438 & \textbf{0.00\%} \\
          &       &       &       &       &       &       &       &       &       &       &  \\
\rowcolor{mygray}  \textbf{Average} & 172.23 & -     & 4.69  & -     & 180.81\% & 4.70  & -     & 174.68\% & 9.56  & -     & 16.61\% \\
   \hline \toprule 
    \end{tabular}%
  \label{tab:addlabel}%
\end{table}%

\subsection{Results for Multiple Centers with Revisiting}
\label{subsec:resultsProblemIV}

\begin{figure}[p]
    \centering
    \begin{subfigure}[b]{0.16\linewidth}
        \includegraphics[width=\linewidth]{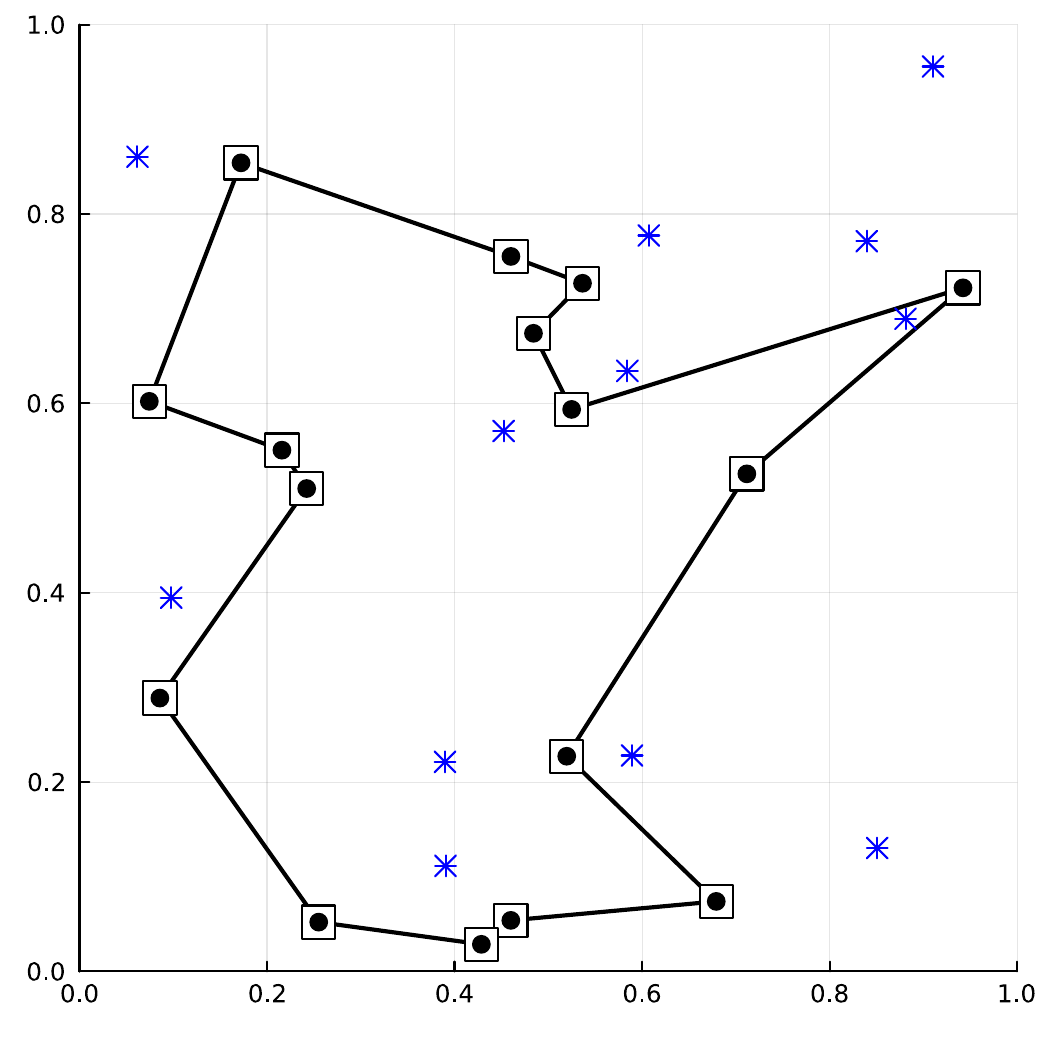}
        \caption{$P1,\;L = 0.2$}
    \end{subfigure}
    \begin{subfigure}[b]{0.16\linewidth}
        \includegraphics[width=\linewidth]{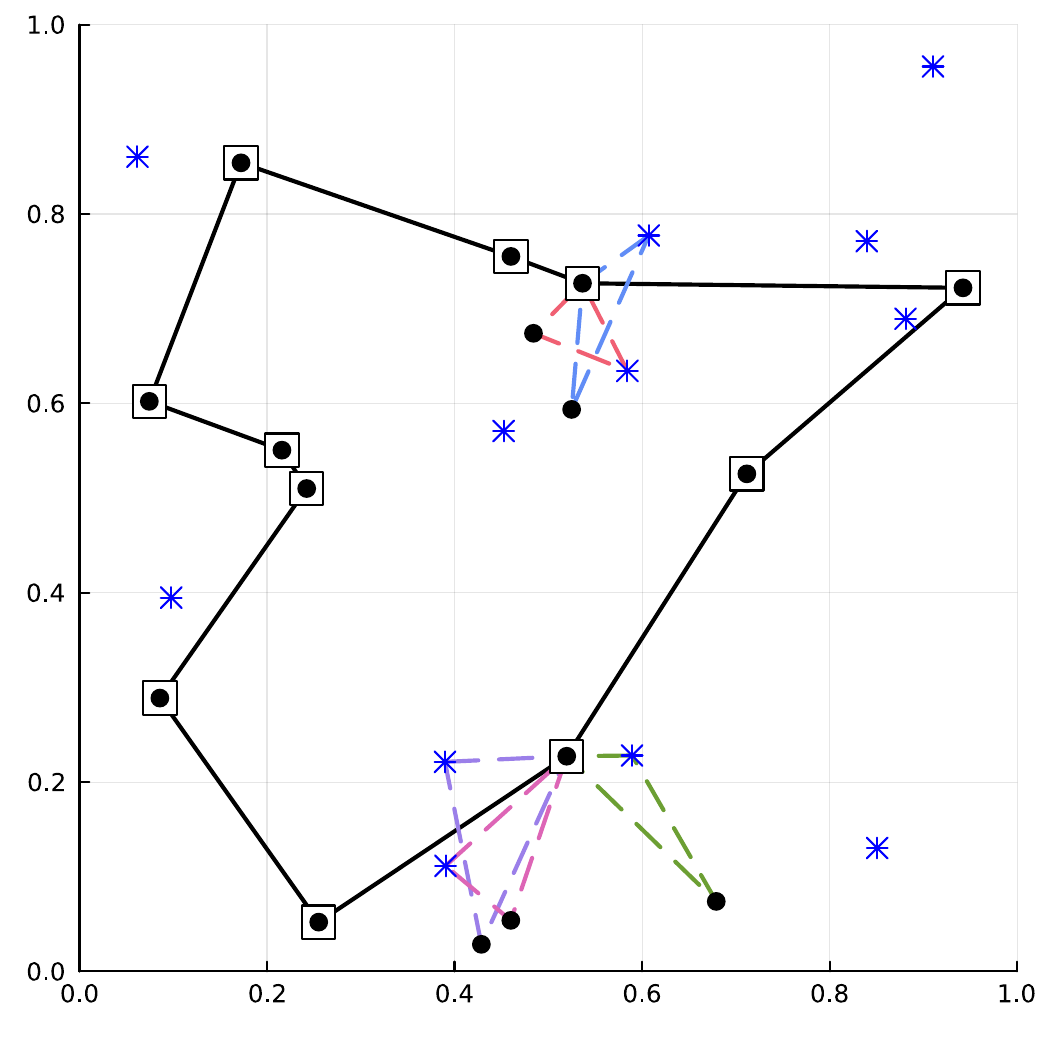}
        \caption{$P1,\;L = 0.6$}
    \end{subfigure}
    \begin{subfigure}[b]{0.16\linewidth}
        \includegraphics[width=\linewidth]{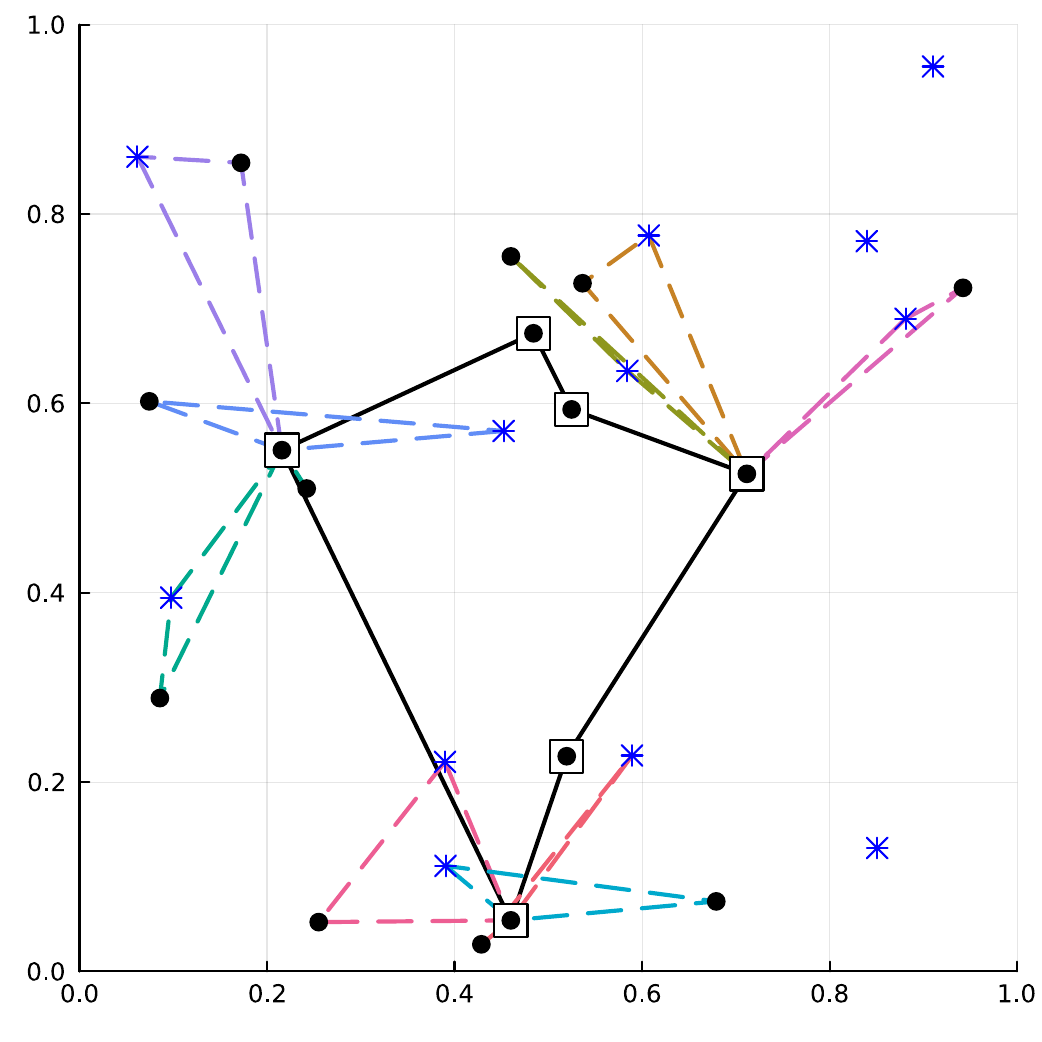}
        \caption{$P1,\;L = 0.8$}
    \end{subfigure} 
    \quad
        \begin{subfigure}[b]{0.16\linewidth}
        \includegraphics[width=\linewidth]{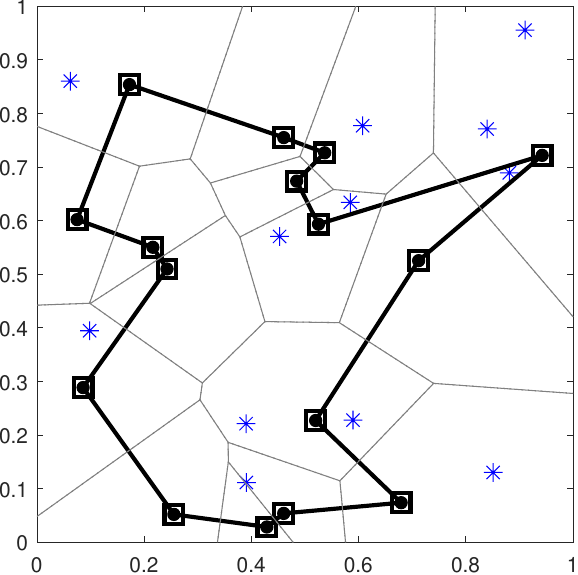}
        \caption{$P1,\;L = 0.2$}
    \end{subfigure}
    \begin{subfigure}[b]{0.16\linewidth}
        \includegraphics[width=\linewidth]{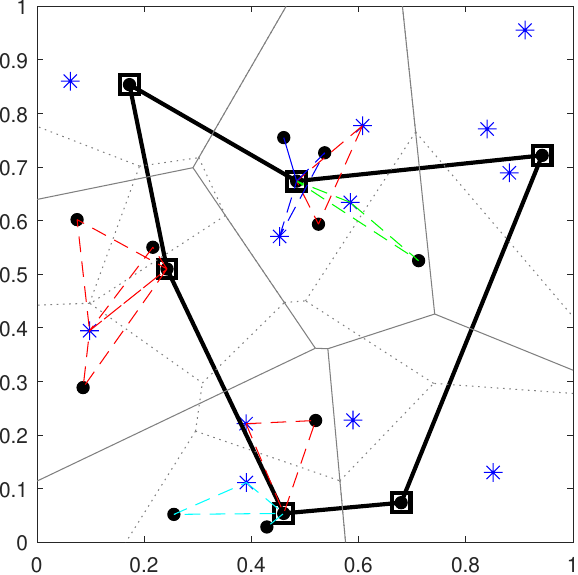}
        \caption{$P1,\;L = 0.6$}
    \end{subfigure}
    \begin{subfigure}[b]{0.16\linewidth}
        \includegraphics[width=\linewidth]{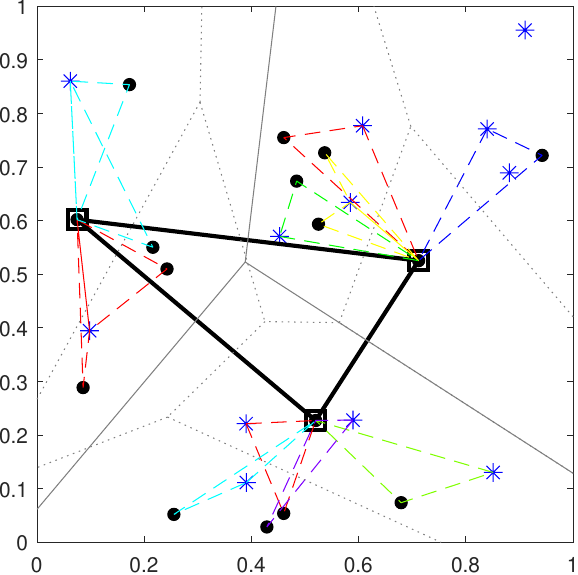}
        \caption{$P1,\;L = 0.8$}
    \end{subfigure} 
    \\
    \begin{subfigure}[b]{0.16\linewidth}
        \includegraphics[width=\linewidth]{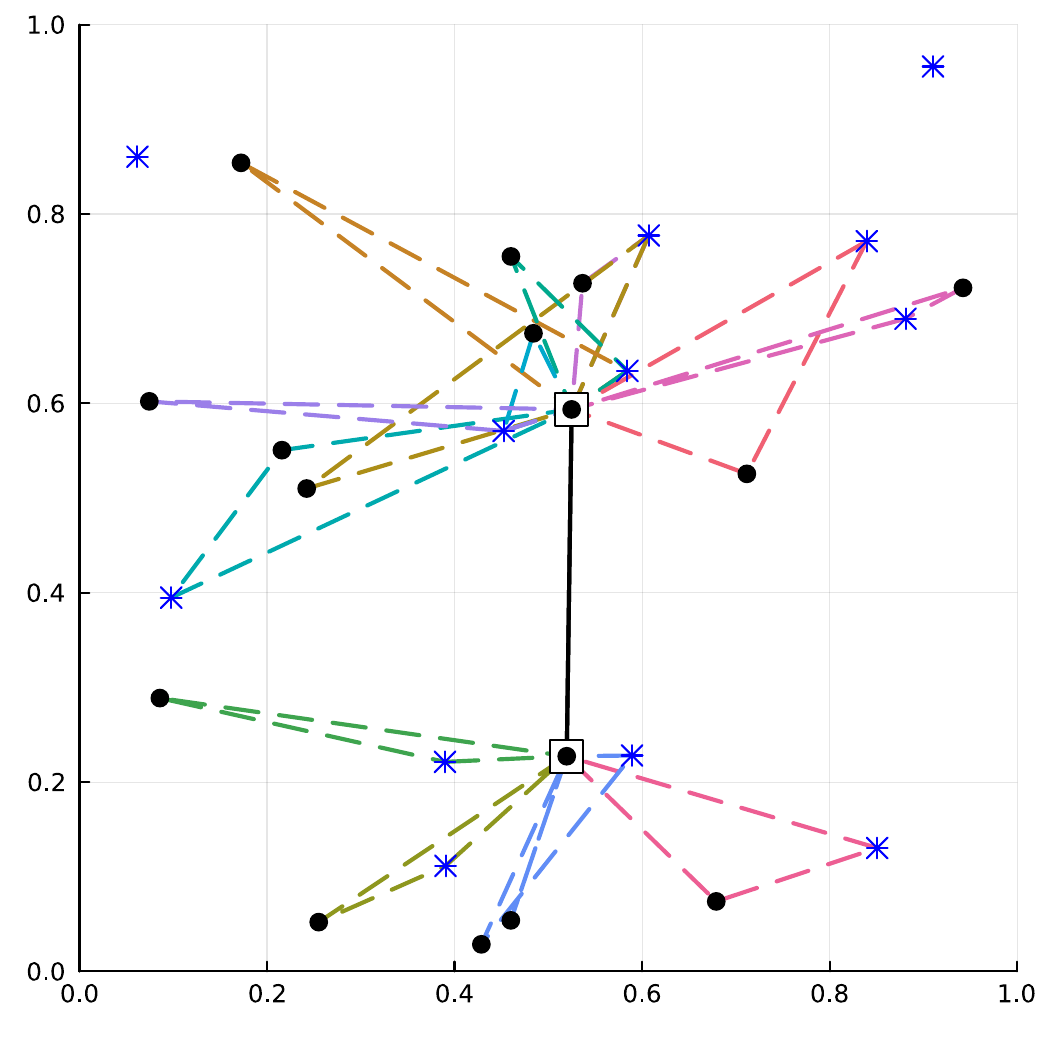}
        \caption{$P1,\;L = 1.0$}
    \end{subfigure}
    \begin{subfigure}[b]{0.16\linewidth}
        \includegraphics[width=\linewidth]{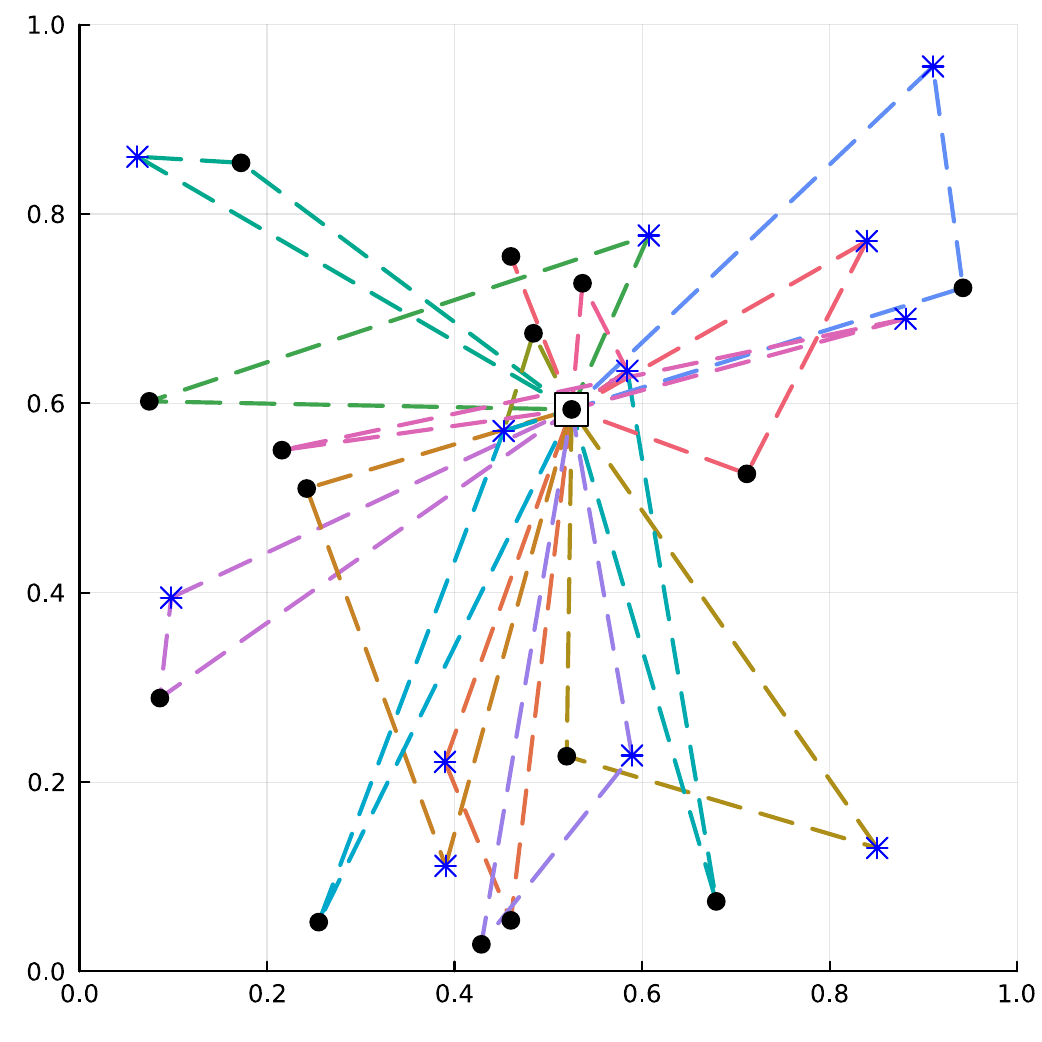}
        \caption{$P1,\;L = 1.4$}
    \end{subfigure}
    \begin{subfigure}[b]{0.16\linewidth}
        \includegraphics[width=\linewidth]{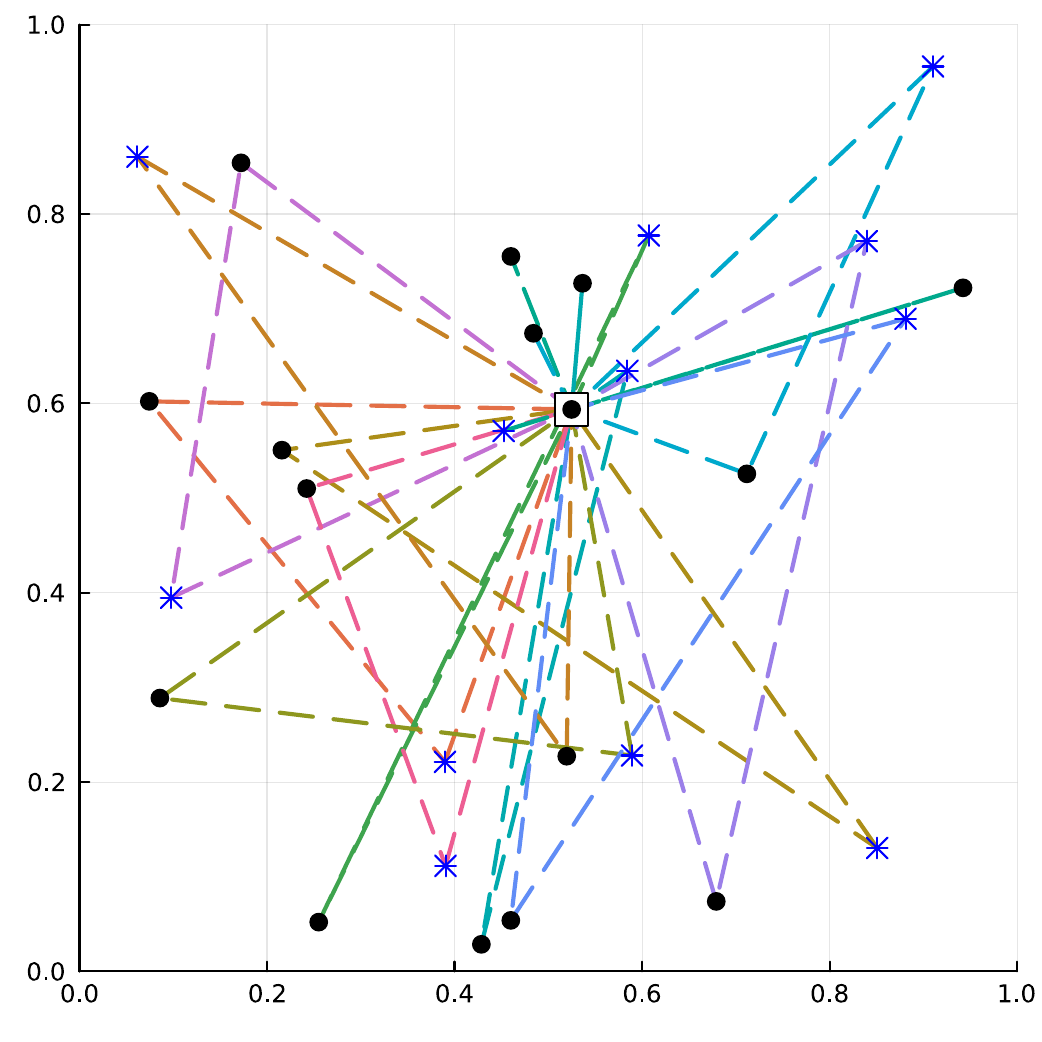}
        \caption{$P1,\;L = 1.8$}
        \label{fig:model-P1-L1.8}
    \end{subfigure}
    \quad
        \begin{subfigure}[b]{0.16\linewidth}
        \includegraphics[width=\linewidth]{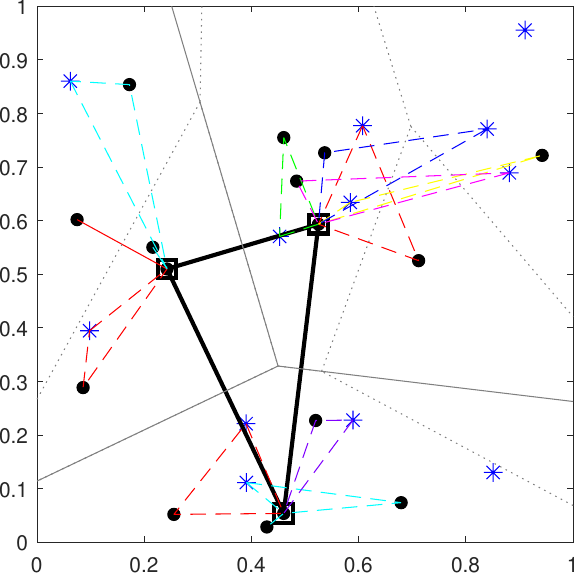}
        \caption{$P1,\;L = 1.0$}
    \end{subfigure}
    \begin{subfigure}[b]{0.16\linewidth}
        \includegraphics[width=\linewidth]{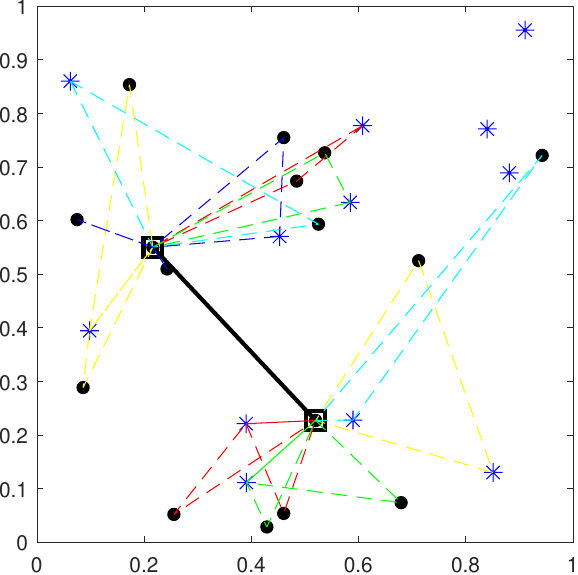}
        \caption{$P1,\;L = 1.4$}
    \end{subfigure}
    \begin{subfigure}[b]{0.16\linewidth}
        \includegraphics[width=\linewidth]{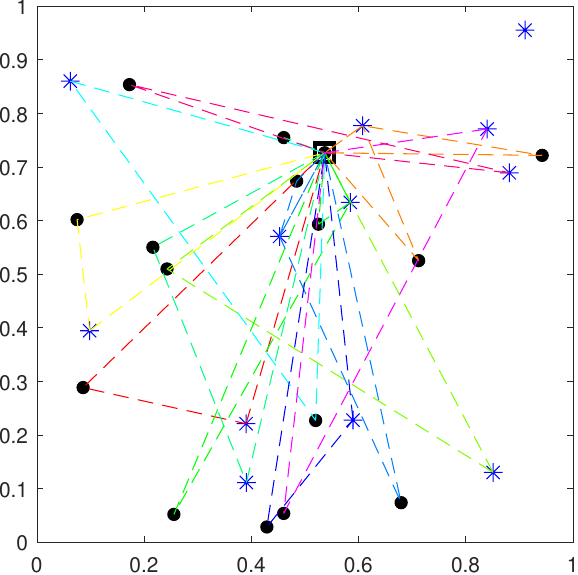}
        \caption{$P1,\;L = 1.8$}
    \end{subfigure} 
    \\
        \begin{subfigure}[b]{0.16\linewidth}
        \includegraphics[width=\linewidth]{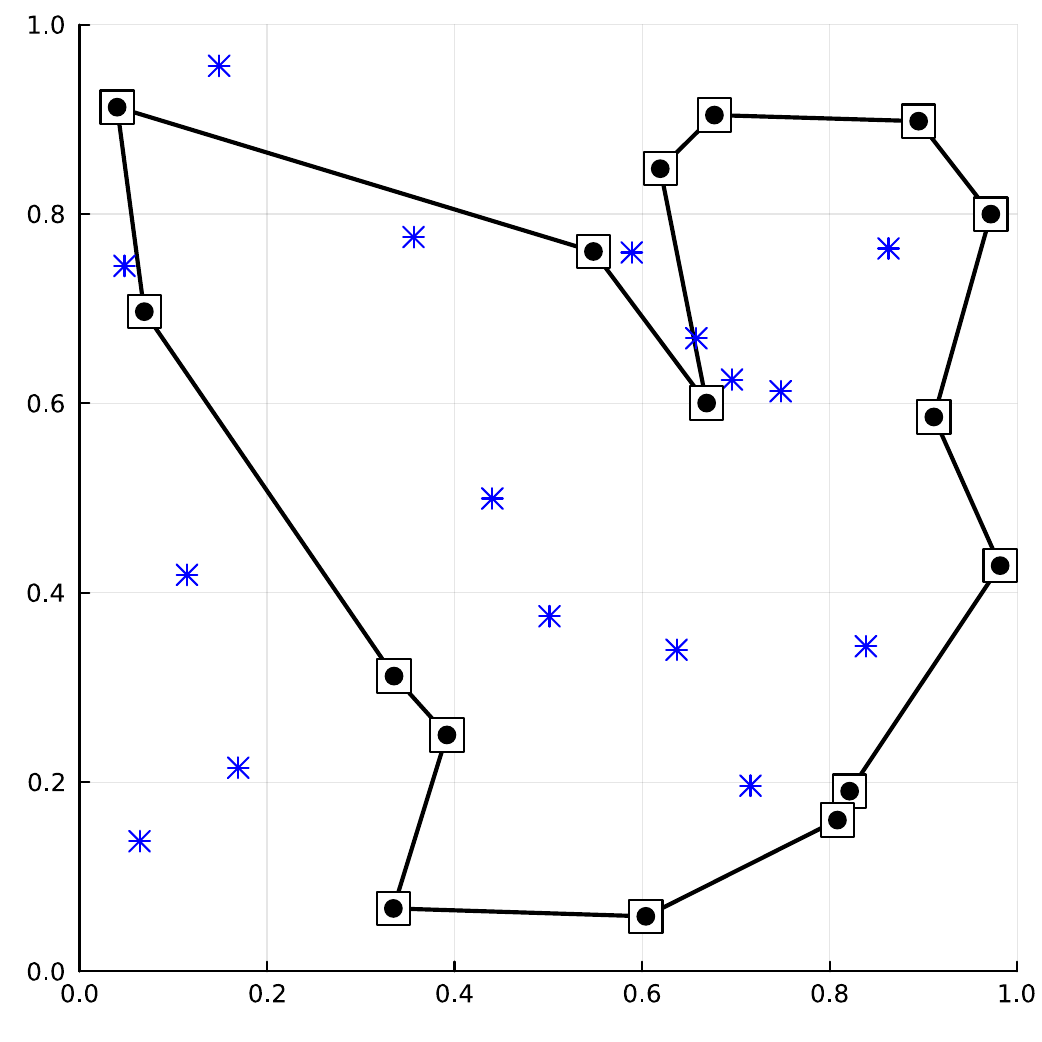}
        \caption{$P2,\;L = 0.2$}
    \end{subfigure}
    \begin{subfigure}[b]{0.16\linewidth}
        \includegraphics[width=\linewidth]{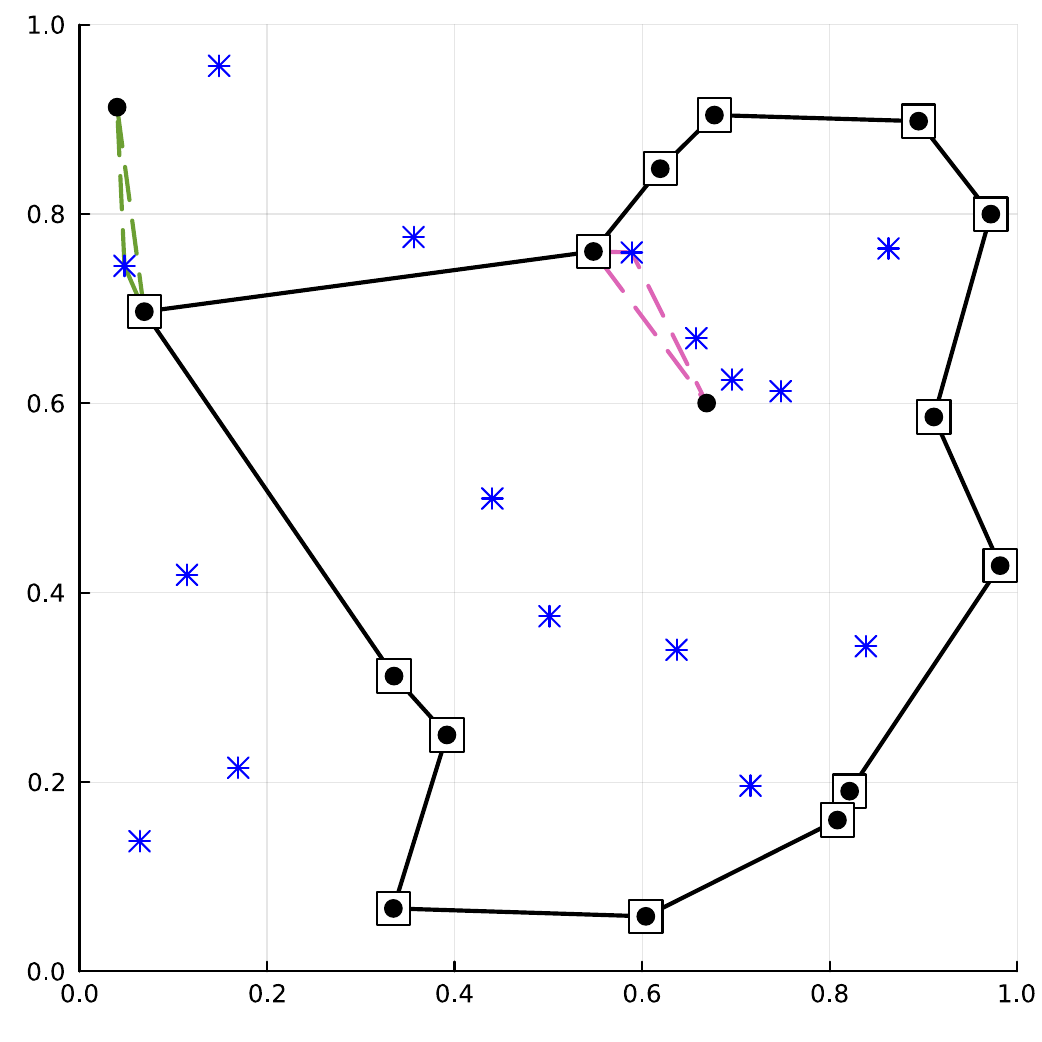}
        \caption{$P2,\;L = 0.6$}
    \end{subfigure}
    \begin{subfigure}[b]{0.16\linewidth}
        \includegraphics[width=\linewidth]{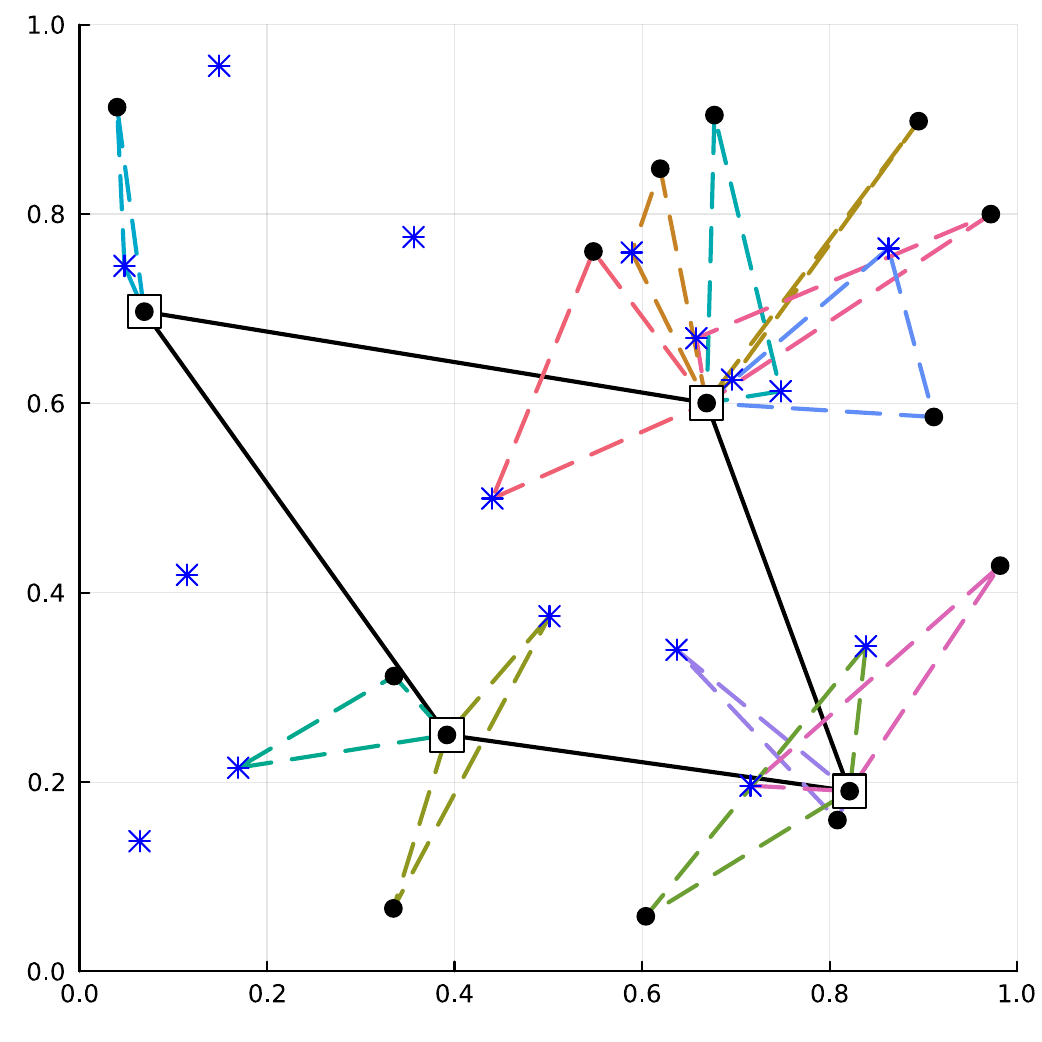}
        \caption{$P2,\;L = 0.8$}
    \end{subfigure} 
    \quad
        \begin{subfigure}[b]{0.16\linewidth}
        \includegraphics[width=\linewidth]{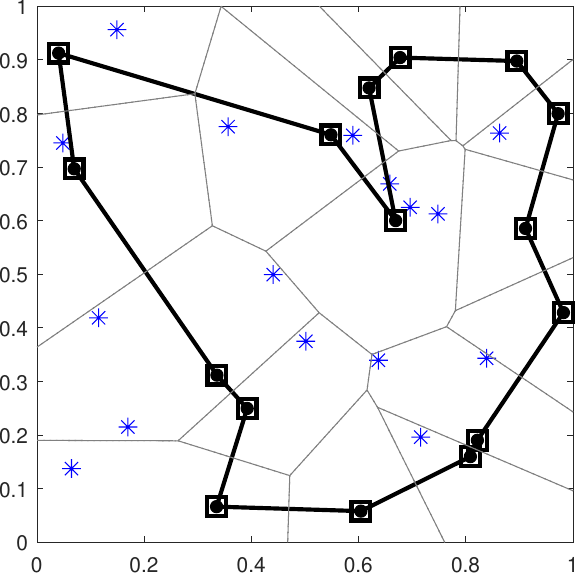}
        \caption{$P2,\;L = 0.2$}
    \end{subfigure}
    \begin{subfigure}[b]{0.16\linewidth}
        \includegraphics[width=\linewidth]{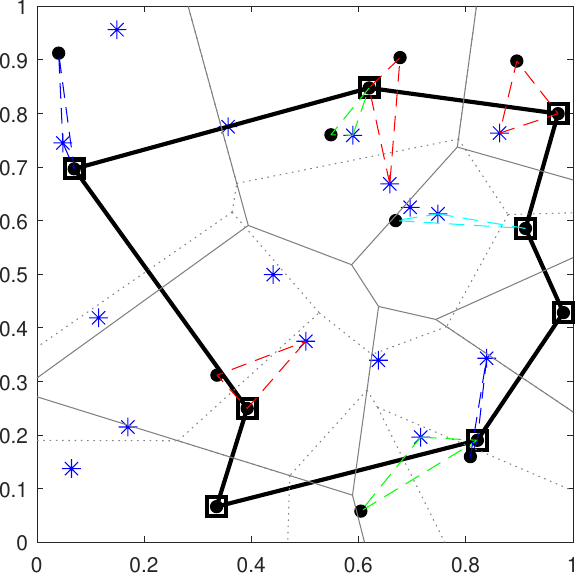}
        \caption{$P2,\;L = 0.6$}
    \end{subfigure}
    \begin{subfigure}[b]{0.16\linewidth}
        \includegraphics[width=\linewidth]{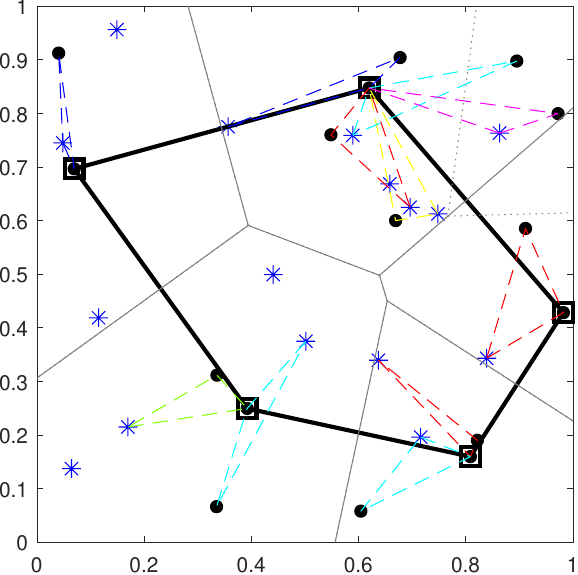}
        \caption{$P2,\;L = 0.8$}
    \end{subfigure} 
    \\
    \begin{subfigure}[b]{0.16\linewidth}
        \includegraphics[width=\linewidth]{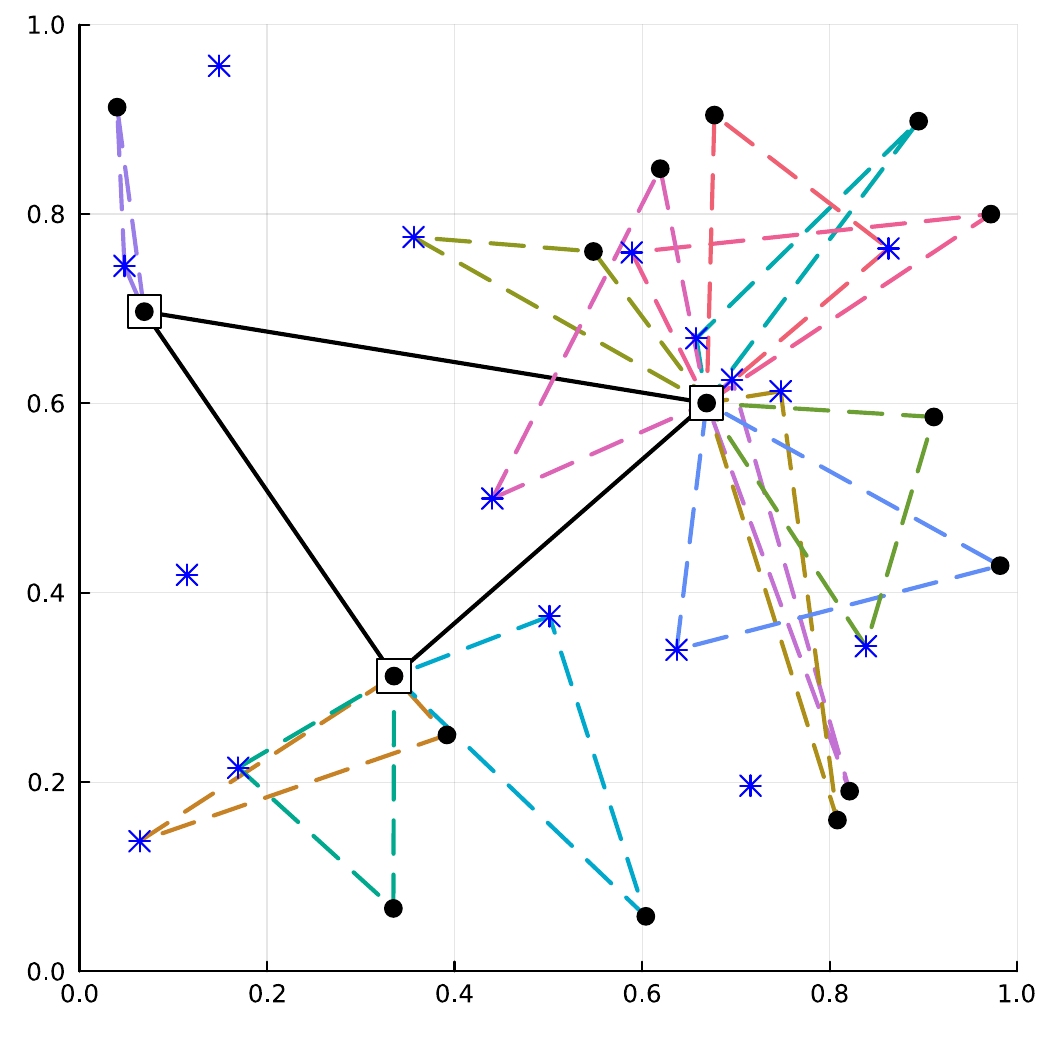}
        \caption{$P2,\;L = 1.0$}
    \end{subfigure}
    \begin{subfigure}[b]{0.16\linewidth}
        \includegraphics[width=\linewidth]{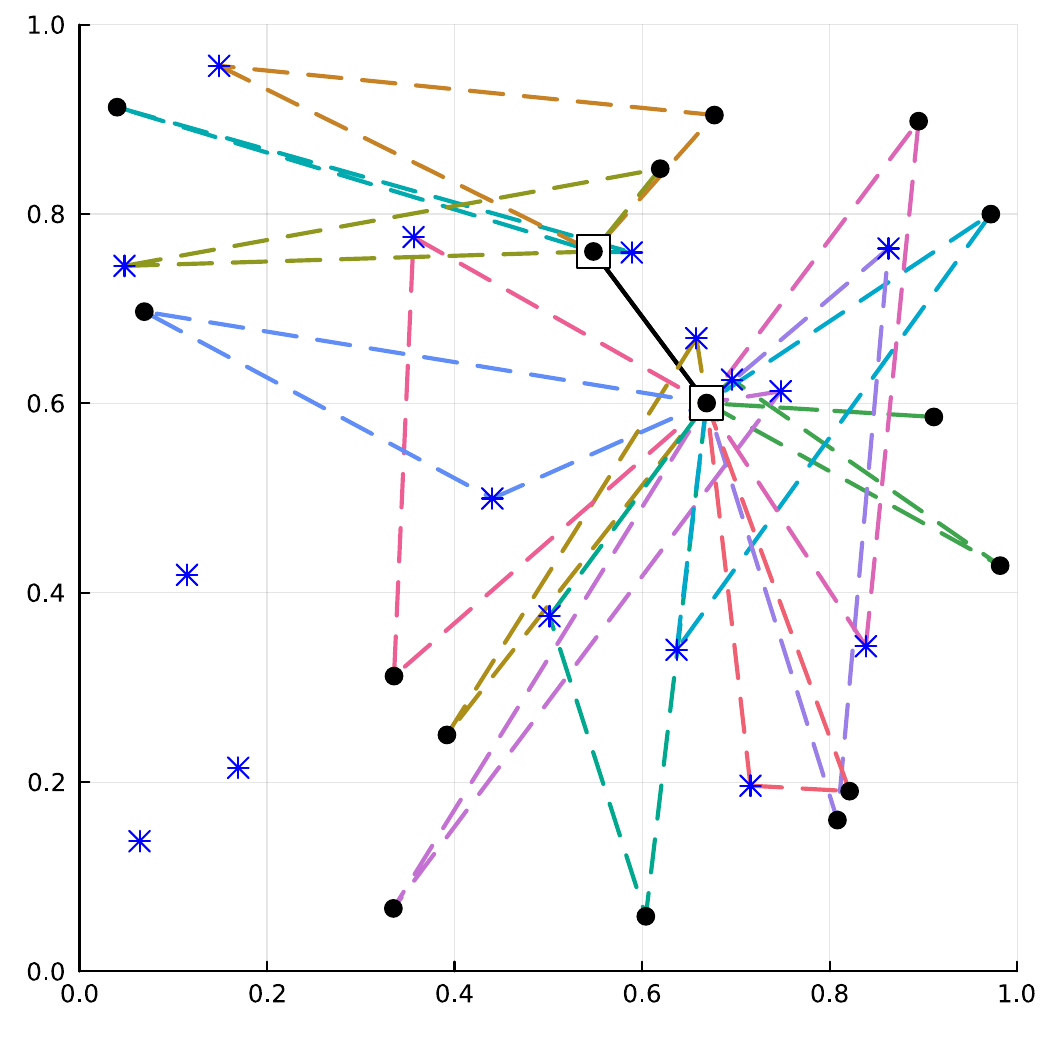}
        \caption{$P2,\;L = 1.4$}
    \end{subfigure}
    \begin{subfigure}[b]{0.16\linewidth}
        \includegraphics[width=\linewidth]{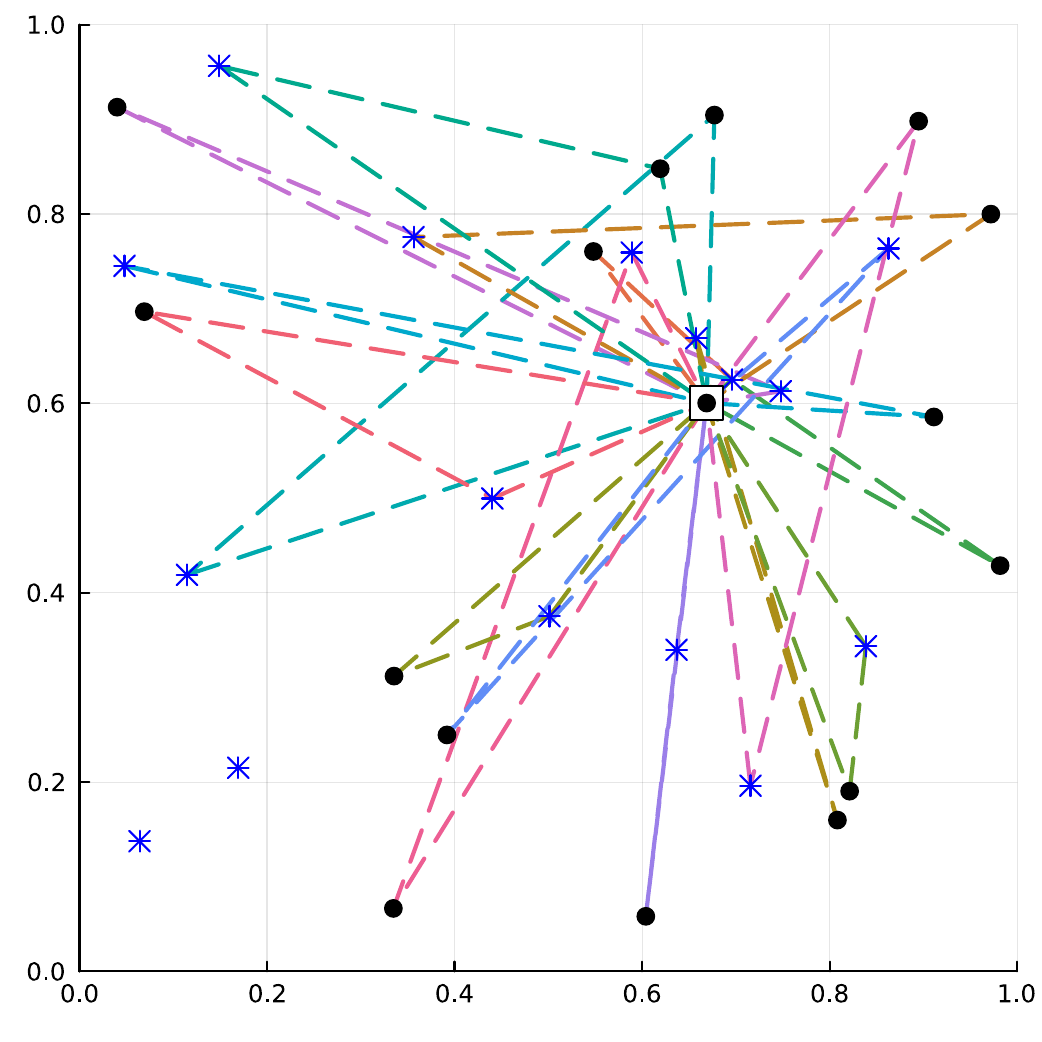}
        \caption{$P2,\;L = 1.8$}
    \end{subfigure}
    \quad
        \begin{subfigure}[b]{0.16\linewidth}
        \includegraphics[width=\linewidth]{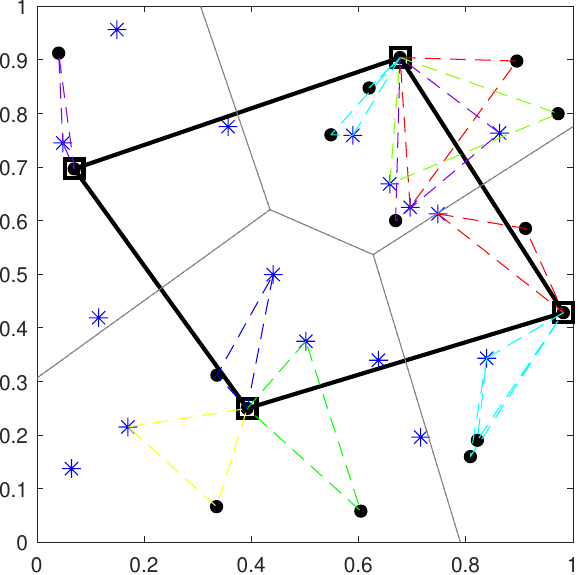}
        \caption{$P2,\;L = 1.0$}
    \end{subfigure}
    \begin{subfigure}[b]{0.16\linewidth}
        \includegraphics[width=\linewidth]{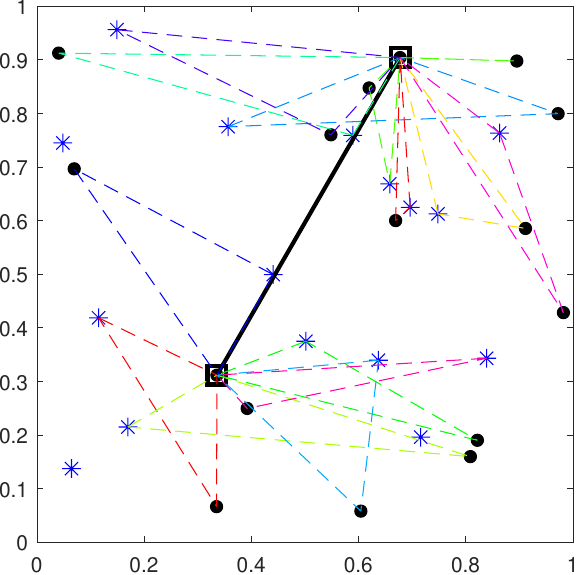}
        \caption{$P2,\;L = 1.4$}
        \label{fig:MC3-P2-L1.4}
    \end{subfigure}
    \begin{subfigure}[b]{0.16\linewidth}
        \includegraphics[width=\linewidth]{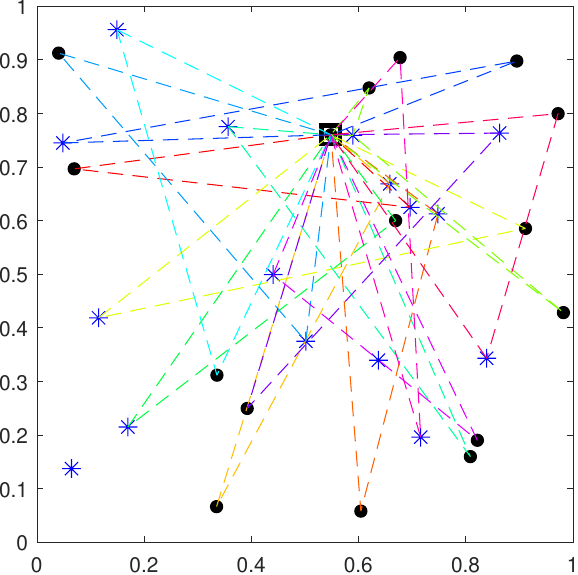}
        \caption{$P2,\;L = 1.8$}
    \end{subfigure} 
    \\
    \begin{subfigure}[b]{0.3\linewidth}
       \includegraphics[width=\linewidth]{legend.pdf}
    \end{subfigure}
    \caption{Examples of solutions for Problem IV (multiple centers with revisiting) on $P1$ and $P2$ instances for select values of maximum drone range $L$ obtained by the optimization model (optimal solutions) on the three columns on the left and by the revisiting version of Algorithm \ref{alg:RouteFinderMC3} (optimal/suboptimal solutions) on the three columns on the right. 
Similar to the results of Problem III, small drone range makes the problem closer to a pure TSP (truck-only) model and for reasonably large $L$ the problem gets closer and closer to the one center (drone-only) delivery problem. 
It can be seen that for some large enough values of $L$ in both model and algorithm solutions, e.g., in (\subref{fig:model-P1-L1.8}) and (\subref{fig:MC3-P2-L1.4}), some of the drones were able to revisit the pickup center after delivering a package at a customer node. Such incidents are distinguishable by hub-and-spoke type (non-triangular) drone routes.}
    \label{fig:resultsModelRevisiting}
\end{figure}

We also tested the MILP model of Problem IV and our algorithms 
on the same problem instances with the same relative drone speed of 2.  
\cref{tab:RevisitingResultsModelAlgComp} summarizes the results of this experiment. We observe that the average computational time of the optimization model, when compared to that of Problem III, is significantly increased. However, this is mostly due to one to two values of $L$ and for the rest the run time is closer to, although still higher than, what we observed for Problem III. This increase is reasonable because the solution space of Problem IV is substantially larger than the solution space of Problem III. The revisiting model achieved better solutions than the recharging model in 13 instances of $P1$ and one instance of $P2$. 
The average running time of the final version of our algorithm (Algorithm \ref{alg:RouteFinderMC3}) is 6.11 seconds versus 3397.84 seconds of the optimization model for $P1$ and 6.66 seconds versus 1202.10 seconds for $P2$. The average (maximum) optimality gap of the solutions provided by Algorithm \ref{alg:RouteFinderMC3} for $P1$ and $P2$ is 28.01\% (279.14\%) and 17.82\% (83.38\%), respectively. 
\cref{fig:resultsModelRevisiting} illustrates solutions of the considered instances of Problem IV solved by our optimization model for select values of $L$ as well as the corresponding solutions obtained by Algorithm \ref{alg:RouteFinderMC3}.

\begin{table}[htbp]
\scriptsize
  \centering
  \caption{\label{tab:RevisitingResultsModelAlgComp} Comparison between the results of the optimization model for Problem IV (Multiple Centers with Revisiting) and the algorithms on two synthetic examples with 16 and 12 and 16 and 16 customers nodes and drone nodes, respectively, distributed randomly in a unit box and for different values of $L$. The column ``Time (s)'' shows the computational time in seconds. The column ``Gap'' presents the optimality gap of the solutions found by the algorithms. The instances for which the optimal solution is found by our algorithms are shown in bold (\textbf{0.00\% gap}). The \underline{underlined} optimal objective values show the instances that the revisiting model achieved better solution than the recharging model. \\}
    \begin{tabular}{cccccccccccc}
    \bottomrule  \hline
          & \multicolumn{11}{c}{\textbf{P1 (16-12) -- Revisiting}} \\
\cmidrule{2-12}    \multirow{2}[4]{*}{\textbf{L }} & \multicolumn{2}{c}{\textbf{Opt. Model}} & \multicolumn{3}{c}{\textbf{Algorithm 1}} & \multicolumn{3}{c}{\textbf{Algorithm 2}} & \multicolumn{3}{c}{\textbf{Algorithm 3}} \\
\cmidrule{2-12}          & Time (s) & Obj. Value  & Time (s) & Obj. Value  & Gap (\%)  & Time (s) & Obj. Value  & Gap (\%)  & Time (s) & Obj. Value  & Gap (\%)  \\
    \midrule
    \multicolumn{1}{c|}{0.2} & 2.99  & 3.3325 & 0.40  & 3.3325 & \textbf{0.00\%} & 0.37  & 3.3325 & \textbf{0.00\%} & 0.48  & 3.3325 & \textbf{0.00\%} \\
    \multicolumn{1}{c|}{0.4} & 18.87 & 3.2561 & 0.63  & 3.6008 & 10.59\% & 0.51  & 3.5956 & 10.42\% & 0.61  & 3.3176 & 1.89\% \\
    \multicolumn{1}{c|}{0.6} & 954.06 & 3.0865 & 0.48  & 3.4071 & 10.39\% & 0.46  & 3.4071 & 10.39\% & 0.75  & 3.6581 & 18.52\% \\
    \multicolumn{1}{c|}{0.8} & 12324.75 & 2.4651 & 0.52  & 3.4377 & 39.45\% & 0.61  & 3.4765 & 41.03\% & 1.25  & 2.8226 & 14.50\% \\
    \multicolumn{1}{c|}{1} & 1758.83 & \underline{1.5250} & 0.68  & 3.2693 & 114.38\% & 0.83  & 3.3081 & 116.92\% & 1.09  & 2.2191 & 45.51\% \\
    \multicolumn{1}{c|}{1.2} & 40669.77 & \underline{1.3796} & 0.95  & 3.2953 & 138.86\% & 0.87  & 2.7749 & 101.14\% & 2.39  & 1.9162 & 38.89\% \\
    \multicolumn{1}{c|}{1.4} & 34.30 & \underline{0.5029} & 0.99  & 2.7771 & 452.25\% & 0.91  & 2.7862 & 454.06\% & 2.03  & 1.9066 & 279.14\% \\
    \multicolumn{1}{c|}{1.6} & 56.83 & \underline{0.4742} & 1.11  & 2.9175 & 515.21\% & 1.09  & 2.6562 & 460.12\% & 9.84  & 0.5181 & 9.25\% \\
    \multicolumn{1}{c|}{1.8} & 87.52 & \underline{0.4563} & 2.23  & 1.6961 & 271.72\% & 2.21  & 1.6961 & 271.72\% & 9.48  & 0.6040 & 32.37\% \\
    \multicolumn{1}{c|}{2} & 120.90 & \underline{0.4563} & 2.35  & 1.6742 & 266.91\% & 2.40  & 1.4503 & 217.86\% & 9.44  & 0.4690 & 2.79\% \\
    \multicolumn{1}{c|}{2.2} & 142.73 & \underline{0.4563} & 2.42  & 1.6742 & 266.91\% & 2.31  & 1.4447 & 216.61\% & 9.55  & 0.4775 & 4.66\% \\
    \multicolumn{1}{c|}{2.4} & 178.49 & \underline{0.4563} & 9.54  & 0.4709 & 3.21\% & 9.74  & 0.4727 & 3.60\% & 9.61  & 0.4775 & 4.66\% \\
    \multicolumn{1}{c|}{2.6} & 228.82 & \underline{0.4563} & 9.47  & 0.4805 & 5.31\% & 9.32  & 0.4856 & 6.43\% & 9.55  & 0.4846 & 6.20\% \\
    \multicolumn{1}{c|}{2.8} & 226.00 & \underline{0.4563} & 9.60  & 0.4832 & 5.90\% & 9.59  & 0.4686 & 2.70\% & 9.53  & 0.4694 & 2.88\% \\
    \multicolumn{1}{c|}{3} & 272.74 & \underline{0.4563} & 9.39  & 0.4828 & 5.81\% & 9.49  & 0.4828 & 5.81\% & 9.38  & 0.4710 & 3.23\% \\
    \multicolumn{1}{c|}{3.2} & 335.45 & \underline{0.4563} & 9.61  & 0.4828 & 5.81\% & 9.40  & 0.4767 & 4.47\% & 9.42  & 0.4828 & 5.81\% \\
    \multicolumn{1}{c|}{3.4} & 350.22 & \underline{0.4563} & 9.39  & 0.4828 & 5.81\% & 9.45  & 0.4742 & 3.93\% & 9.51  & 0.4832 & 5.90\% \\
          &       &       &       &       &       &       &       &       &       &       &  \\
      \rowcolor{mygray}  \textbf{Average} & 3397.84 & -     & 4.10  & -     & 124.62\% & 4.09  & -     & 113.36\% & 6.11  & -     & 28.01\% \\
    \midrule
          &       &       &       &       &       &       &       &       &       &       &  \\
          &       &       &       &       &       &       &       &       &       &       &  \\
    \midrule
          & \multicolumn{11}{c}{\textbf{P2 (16-16) -- Revisiting}} \\
\cmidrule{2-12}    \multirow{2}[4]{*}{\textbf{L }} & \multicolumn{2}{c}{\textbf{Opt. Model}} & \multicolumn{3}{c}{\textbf{Algorithm 1}} & \multicolumn{3}{c}{\textbf{Algorithm 2}} & \multicolumn{3}{c}{\textbf{Algorithm 3}} \\
\cmidrule{2-12}          & Time (s) & Obj. Value  & Time (s) & Obj. Value  & Gap (\%)  & Time (s) & Obj. Value  & Gap (\%)  & Time (s) & Obj. Value  & Gap (\%)  \\
    \midrule
    \multicolumn{1}{c|}{0.2} & 3.44  & 3.5817 & 0.31  & 3.5817 & \textbf{0.00\%} & 0.33  & 3.5817 & \textbf{0.00\%} & 0.33  & 3.5817 & \textbf{0.00\%} \\
    \multicolumn{1}{c|}{0.4} & 55.07 & 3.5817 & 0.28  & 3.5817 & \textbf{0.00\%} & 0.30  & 3.5817 & \textbf{0.00\%} & 0.38  & 3.7663 & 5.15\% \\
    \multicolumn{1}{c|}{0.6} & 16.14 & 3.2325 & 0.38  & 3.7280 & 15.33\% & 0.34  & 3.7280 & 15.33\% & 0.43  & 3.7508 & 16.03\% \\
    \multicolumn{1}{c|}{0.8} & 97.30 & 2.7967 & 0.58  & 3.5360 & 26.43\% & 0.65  & 3.5360 & 26.43\% & 0.85  & 3.3767 & 20.74\% \\
    \multicolumn{1}{c|}{1} & 475.36 & 2.2349 & 0.86  & 3.5589 & 59.24\% & 1.10  & 3.5571 & 59.16\% & 0.93  & 3.5588 & 59.23\% \\
    \multicolumn{1}{c|}{1.2} & 7490.24 & 1.7369 & 1.15  & 3.1890 & 83.60\% & 1.07  & 3.1868 & 83.47\% & 1.06  & 3.1851 & 83.38\% \\
    \multicolumn{1}{c|}{1.4} & 10233.86 & \underline{1.1598} & 1.13  & 3.1689 & 173.23\% & 1.09  & 3.1851 & 174.64\% & 2.16  & 2.1160 & 82.45\% \\
    \multicolumn{1}{c|}{1.6} & 24.63 & 0.4673 & 1.09  & 3.1689 & 578.16\% & 1.04  & 2.9051 & 521.71\% & 10.62 & 0.5055 & 8.18\% \\
    \multicolumn{1}{c|}{1.8} & 57.27 & 0.4438 & 0.95  & 3.0112 & 578.58\% & 1.06  & 2.9436 & 563.34\% & 10.80 & 0.4774 & 7.59\% \\
    \multicolumn{1}{c|}{2} & 54.01 & 0.4438 & 0.98  & 3.4090 & 668.23\% & 0.96  & 3.1674 & 613.78\% & 10.66 & 0.4438 & \textbf{0.00\%} \\
    \multicolumn{1}{c|}{2.2} & 81.60 & 0.4438 & 5.38  & 1.8542 & 317.85\% & 5.40  & 1.8129 & 308.55\% & 10.80 & 0.4675 & 5.36\% \\
    \multicolumn{1}{c|}{2.4} & 124.37 & 0.4438 & 2.16  & 2.1089 & 375.25\% & 2.12  & 2.0368 & 358.99\% & 10.83 & 0.4631 & 4.37\% \\
    \multicolumn{1}{c|}{2.6} & 167.88 & 0.4438 & 2.13  & 2.1235 & 378.54\% & 2.06  & 1.7725 & 299.44\% & 10.55 & 0.4640 & 4.56\% \\
    \multicolumn{1}{c|}{2.8} & 280.61 & 0.4438 & 10.69 & 0.4438 & \textbf{0.00\%} & 10.65 & 0.4438 & \textbf{0.00\%} & 10.78 & 0.4701 & 5.94\% \\
    \multicolumn{1}{c|}{3} & 295.28 & 0.4438 & 10.61 & 0.4461 & 0.53\% & 10.65 & 0.4438 & \textbf{0.00\%} & 10.56 & 0.4438 & \textbf{0.00\%} \\
    \multicolumn{1}{c|}{3.2} & 460.98 & 0.4438 & 10.63 & 0.4673 & 5.30\% & 10.78 & 0.4626 & 4.25\% & 10.71 & 0.4438 & \textbf{0.00\%} \\
    \multicolumn{1}{c|}{3.4} & 517.70 & 0.4438 & 10.79 & 0.4438 & \textbf{0.00\%} & 10.77 & 0.4673 & 5.30\% & 10.77 & 0.4438 & \textbf{0.00\%} \\
          &       &       &       &       &       &       &       &       &       &       &  \\
      \rowcolor{mygray}  \textbf{Average} & 1202.10 & -     & 3.54  & -     & 191.78\% & 3.55  & -     & 178.49\% & 6.66  & -     & 17.82\% \\
    \hline \toprule
    \end{tabular}%
  \label{tab:addlabel}%
\end{table}%

\subsection{Similarities Between the Results of the Two Problems}
\label{subsec:similaritiesProblemsIII-IV}

It is apparent in Tables \ref{tab:RechargingResultsModelAlgComp} and \ref{tab:RevisitingResultsModelAlgComp} that the objective function value for the optimization model and the algorithms decreases in $L$ as expected.
We also see in Tables \ref{tab:RechargingResultsModelAlgComp} and \ref{tab:RevisitingResultsModelAlgComp} that starting from $L=1.8$ the optimal objective function value remains the same. This is because at this level, the range of the drone is long enough to serve all customers with drones and thus the truck would not be used. This means for large enough $L$ and large enough relative speed of drone to truck, the problem essentially reduces to a one center problem with recharging (Problem I). This is also visible in Figures \ref{fig:resultsModelRecharging} and \ref{fig:resultsModelRevisiting}. 

Furthermore, it can be seen in the results from both Problems III and IV that the running time of the optimization model and the optimality gap of our algorithm are both very sensitive to the values of $L$. It is much easier to solve the problem when $L$ is either small or large. Optimality gap for the middle values of $L$ (between 1.0 to 1.4 in these instances) is also larger.
If the drone range is very low, which means drones can only fulfill the demand of nearby customer nodes, the problem is a pure or nearly TSP problem which is much faster to solve since its feasible region is smaller. In contrast, if the drone maximum allowed distance is in a reasonable range (between 1.0 to 1.4 in these instances), drones can be a big help in delivering the packages while still leaving a good portion of the customers to be served by the truck.  Finding both drone routes and the truck route takes considerable time in this case, although the lion share of total time is due to the drone routes. Therefore, the feasibility region is much bigger than a TSP which increases the run time and makes finding optimal or high quality solutions harder. On the other end of the range, if the drone range is large enough (e.g., larger than 1.8), we can do all deliveries with drones without using the truck. In this case, the solving time is reasonable but still a bit higher than pure TSP. This is because larger values of $L$ increase the number of feasible drone routes. Moreover, we still have to determine which node(s) will serve as the truck node(s), and what's the best truck route (if any). This explains why there is a little bit of increase in solving time at the right tail of the drone range compared to the left tail. Although the solution space might be even bigger than that of the mid-level range of $L$ (e.g., between 1.0 to 1.4), it is much faster than that case because good upper bounds are found quickly (with for example fewer (or even one) centers) which helps skipping a large number of feasible solutions that use bigger number of centers. 
For the same reasons we can see a similar behavior for the optimality gap 
and when $L$ is either small or large, the gap between the algorithm and the optimal solution is much smaller than the results of the mid-level cases of drone range. 

\subsection{Comparison with Traditional Delivery System}
\label{subsec:comparisonTSP}

\begin{figure}[t]
    \centering
    \begin{subfigure}[b]{0.22\linewidth}
        \includegraphics[width=\linewidth]{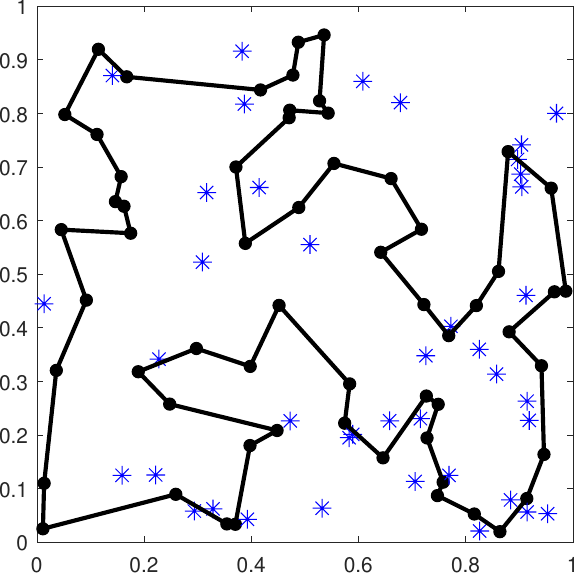}
        \caption{Total Time = 6.1880}
    \end{subfigure}
    \quad
    \begin{subfigure}[b]{0.22\linewidth}
        \includegraphics[width=\linewidth]{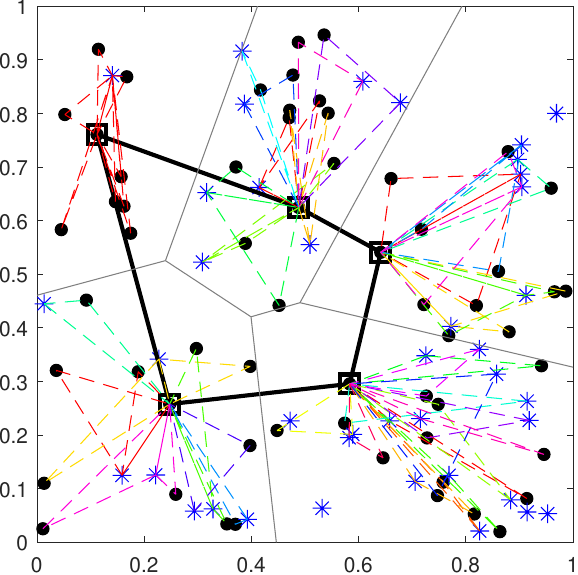}
        \caption{Total Time = 5.0042}
    \end{subfigure}
    \qquad
    \begin{subfigure}[b]{0.22\linewidth}
        \includegraphics[width=\linewidth]{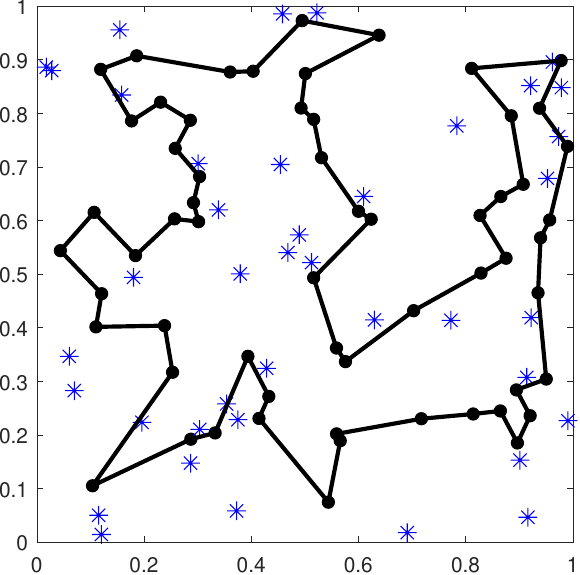}
        \caption{Total Time = 5.8312}
    \end{subfigure}
    \quad
    \begin{subfigure}[b]{0.22\linewidth}
        \includegraphics[width=\linewidth]{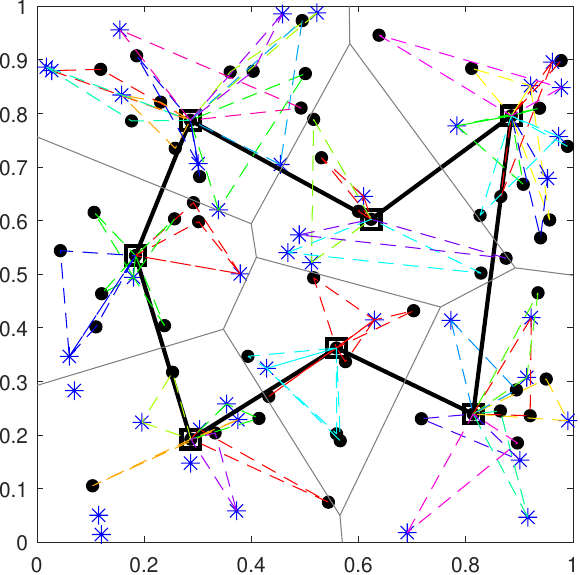} 
        \caption{Total Time = 4.5983}
    \end{subfigure}
    \\
     \begin{subfigure}[b]{0.22\linewidth}
        \includegraphics[width=\linewidth]{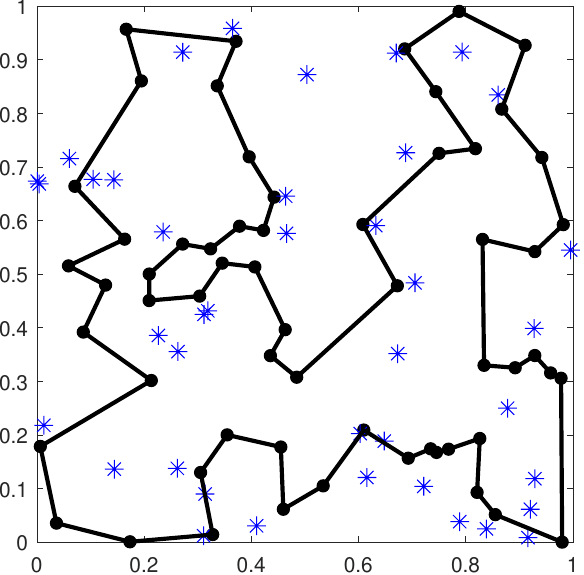}
        \caption{Total Time = 6.5363}
    \end{subfigure}
    \quad
    \begin{subfigure}[b]{0.22\linewidth}
        \includegraphics[width=\linewidth]{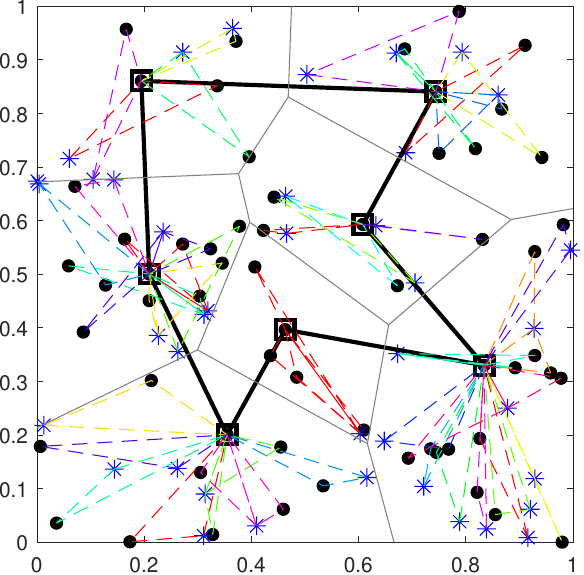} 
        \caption{Total Time = 4.9806}
    \end{subfigure}
    \qquad
    \begin{subfigure}[b]{0.22\linewidth}
        \includegraphics[width=\linewidth]{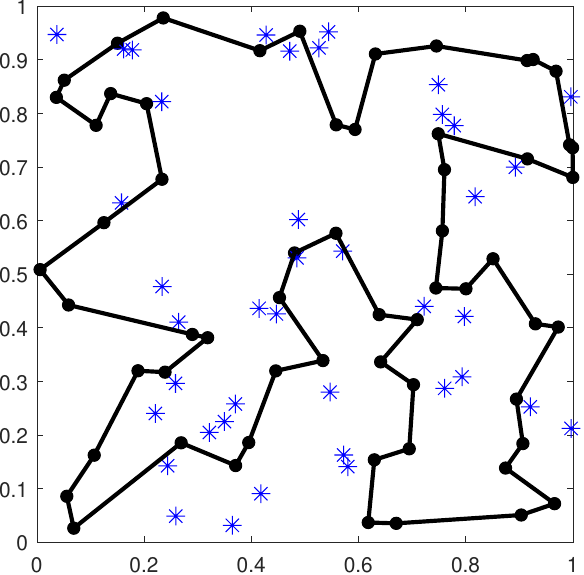}
        \caption{Total Time = 6.2583}
    \end{subfigure}
    \quad
    \begin{subfigure}[b]{0.22\linewidth}
        \includegraphics[width=\linewidth]{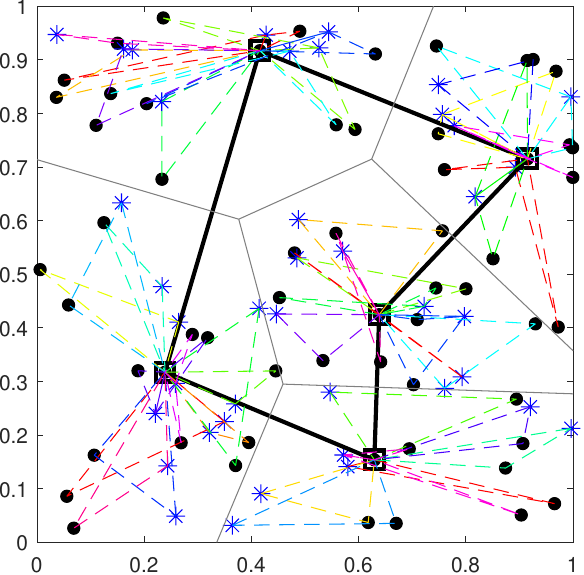} 
        \caption{Total Time = 3.9290}
    \end{subfigure}
    \\
    \begin{subfigure}[b]{0.35\linewidth}
       \includegraphics[width=\linewidth]{legend.pdf}
    \end{subfigure}
    \caption{Sample solutions for a multiple center with recharging problem with 60 customer nodes and 40 drone nodes for $L=0.8$ obtained by our Algorithm \ref{alg:RouteFinderMC3}, representing a decentralized delivery system, for different random examples in a unit box with their objective function values are illustrated side by side to their corresponding TSP solutions (centralized truck-only delivery system) for the same instances obtained by the LKH algorithm.}
    \label{fig:TSV-vs-SharedDrones-60-40}
\end{figure}

We also compared our approach with the traditional centralized truck-only delivery system in which a single truck would deliver all packages, i.e., a pure TSP model. We ran both recharging and revisiting versions of Algorithm \ref{alg:RouteFinderMC3}, which provides a decentralized alternative, as well as the LKH algorithm for the pure TSP model on several problem instances with 60 customer nodes and 40 drone nodes distributed uniformly at random in a unit square. Again, we set the relative speed of drone to be 2. For drone range we chose $L=0.8$ to ensure having both truck and drone routes.
Figure \ref{fig:TSV-vs-SharedDrones-60-40} shows solutions obtained by our algorithm presented side-by-side with the corresponding TSP solution for each instance. It is easy to see the advantages of our proposed decentralized system over the conventional centralized delivery system.
On these examples, the average delivery time when we combine a truck with crowdsourced drones was 4.6280, while the average delivery time for the same problem if the truck serves all the customers was 6.2035. This crude comparison suggests an almost 25.32\% improvement in the efficiency of the last-mile delivery, in our randomly generated examples, if we combine truck delivery and drone delivery in the context of sharing economy platforms. Certainly, this would highly depend on the chosen input parameters that we investigate further in Sections \ref{sec:SensitivityAnalysis} and \ref{sec:CarbonFootprint}.

\section{Sensitivity Analysis}
\label{sec:SensitivityAnalysis}
In the rest of the paper, we perform an extensive sensitivity analysis with respect to several factors to study their impact on the quality of the solution and savings of the shared delivery system compared to the traditional truck-only delivery. 
These factors include speed of drones, number of available drones, a measure combining speed and number of drones, and customer distribution. Moreover, we analyze the impact of improvements presented in \cref{sec:Improvements} for different values of the input parameters.
Finally, a comparison is made between three models to measure the impact of shared delivery model on carbon footprint. These three models are the traditional truck-only delivery (TSP) model, delivery with a truck and a drone where the truck carries a drone and both deliver packages in a coordinated way (we call this model TSPD), and our truck and shared drone delivery model (TSP-SD). Almost all simulations are done on an instance of the recharging problem with 60 customer nodes and several drone nodes (the number varies based on the input parameters) distributed uniformly at random or according to a distribution (depending on the analysis) in a unit square assuming $L=0.8$. In most cases, the revisiting model delivers similar results but wherever we expect a significant difference we run the simulation on the revisiting problem as well.

To analyze the sensitivity of the proposed models and algorithms that combines a truck with crowdsourced drones in the last-mile delivery operation to input parameters, we consider the impact of several factors including relative speed of drones, relative number of crowdsourced drones, and the customer and drone distribution on the solutions obtained by our algorithm. We expect similar behavior in the optimal solutions. Here we define two parameters as: $\rho_1 = \frac{V_d}{V_t}=\frac{\text{speed of drone}}{\text{speed of truck}}$ and $\rho_2 = \frac{m}{n}=\frac{\text{number of drones}}{\text{number of customers}}$. In terms of customer distributions, a multivariate normal distribution is assumed and its variance will be adjusted. We compare the results with respect to the \emph{savings} in time by using the crowdsourced drones that is 
\[
    \frac{\text{Time cost of TSP (truck only) - Time cost of TSP-SD (crowdsourced drones combined with a truck)}}{\text{Time cost of TSP (truck only)}} \,,
\]
i.e., we compare our solutions with those obtained by the LKH algorithm. 

When comparing Algorithm \ref{alg:RouteFinderMC3} to the results of LKH for pure TSP, one should keep in mind that any savings suggested by our algorithm in despite the fact that for our algorithm with $L=0.8$, we usually have a considerable optimality gap, while the optimality gap of LKH for pure TSP is almost zero. This difference could be seen as additional savings if our problem is solved to optimality.

\subsection{Relative Speed of Drones} 
\label{sec:RelativeDroneSpeed}
To analyze the sensitivity with respect to $\rho_1$, we adjust this parameter in the range $0.75\leq \rho_1 \leq 3$ that means the drone speed can vary from 0.75 times the truck speed to three times the truck speed, while we fix the number of drones to be half or twice of the number of customers, i.e., $\rho_2\in\{0.5,2\}$. 
The choice of the range for $\rho_1$ is driven by the fact that it is known in the literature of the centralized alternative (TSPD) that for having any significant gain by combining a truck and a drone in delivery operations we must have $\rho_1 \geq 2$ \cite{Ferrandez2016optimization,agatz2018optimization} and it is most often assumed to have $\rho_1 \geq 1$.

\begin{figure}[h!]
    \centering
    \begin{subfigure}[b]{0.35\linewidth}
        \includegraphics[width=\linewidth]{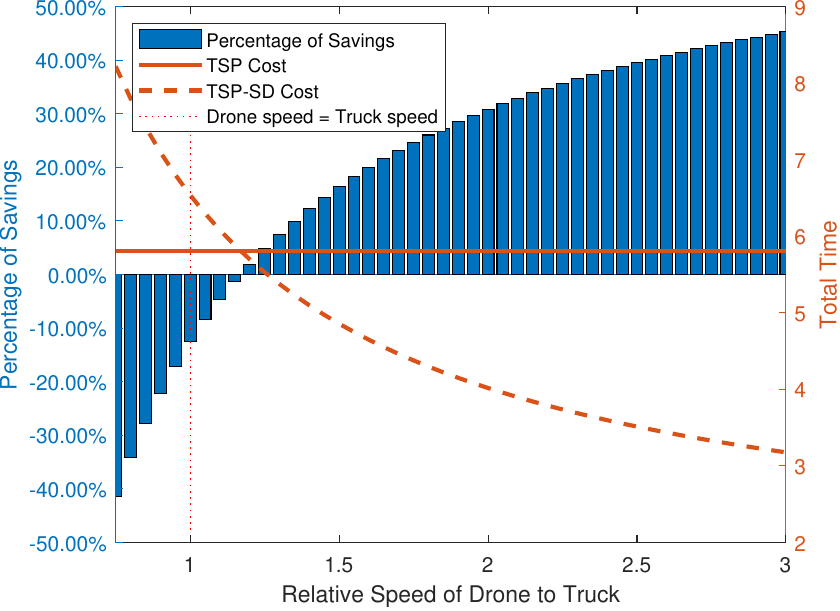}
        \caption{$\rho_2 = 0.5$}
        \label{fig:SA1-rho2Half}
    \end{subfigure}
    \qquad
    \begin{subfigure}[b]{0.35\linewidth}
        \includegraphics[width=\linewidth]{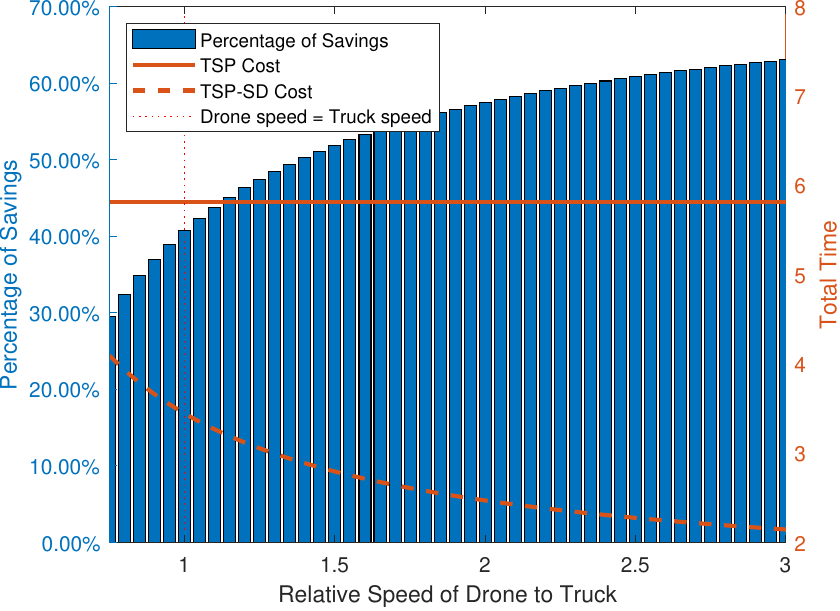}
        \caption{$\rho_2 = 2$}
        \label{fig:SA1-rho2Two}
    \end{subfigure}
    \caption{Sensitivity analysis of the results of Algorithm \ref{alg:RouteFinderMC3} with respect to the relative speed of drones for two cases of drones supply.}
    \label{fig:SA1}
\end{figure}
From \cref{fig:SA1-rho2Half} we can find that, when there is considerable shortage of drones relative to the number of customers, only when the drone speed is less than 0.75 times the speed of the truck ($\rho_1 \leq 1.15$), the TSP model will have less time cost than the model with shared drones. Also, from \cref{fig:SA1-rho2Two} we can observe that 
when relative drone supply increases to 2 the chart shifts to left, so to speak, and we see a positive and significant savings for the shared model for all values of $\rho_1$ in the considered range ($0.75 \leq \rho_1 \leq 3$). 
In both situations, the TSP with sharing drones model has a consistently better performance than the TSP model when drones are reasonably faster than the speed of the truck ($\rho_1 \geq 1.15$). Intuitively, the combined model gains help from the drones by working (potentially) in parallel to each other and to the truck, and on the other hand loses some time by idling the truck at the cluster centers to wait for the drones to fly back and forth and deliver packages. The overall, tradeoff between these two opposite forces shows that the combined model has a significant advantage over the pure TSP model and the proposed model of combined truck and shared drones performs significantly better. This advantage becomes bigger by increasing drones speed relative to the speed of truck ($\rho_1$) with a diminishing return. Both these properties are quite consistent with the results of \cite{carlsson2018coordinated} that shows the gain for a TSPD model over TSP is proportional to $\sqrt{\rho_1}$. Finally, we also see that increasing the supply of drones ($\rho_2$) increases the savings. This phenomena is similar to growing the supply side in other sharing economy (two-sided market) applications.

\subsection{Relative Number of Drones} 
\label{sec:RelativeNumberDrones}
We analyze the sensitivity with respect to $\rho_2$, by varying this parameter in the range $0.3 \leq \rho_2 \leq 2$, while we fix the speed of drones to either be one or two times that of the truck, i.e., $\rho_1\in\{1,2\}$. The choice of the range for $\rho_2$ is motivated by the well-known importance of having a balance between supply and demand in two-sided markets, i.e., $\rho_2=1$.

\begin{figure}[h!]
    \centering
    \begin{subfigure}[b]{0.35\linewidth}
        \includegraphics[width=\linewidth]{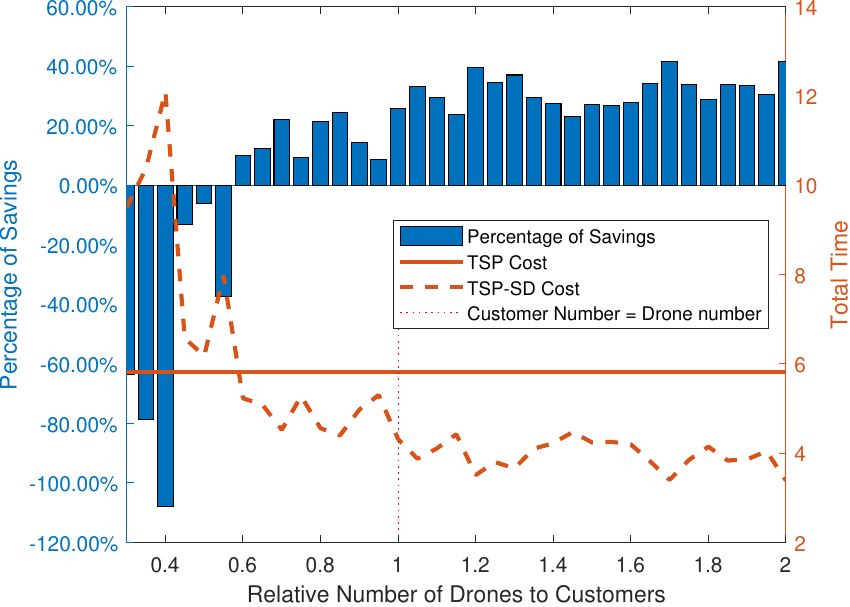}
        \caption{$\rho_1 = 1$}
         \label{fig:SA2rho1One}
    \end{subfigure}
    \qquad
    \begin{subfigure}[b]{0.35\linewidth}
        \includegraphics[width=\linewidth]{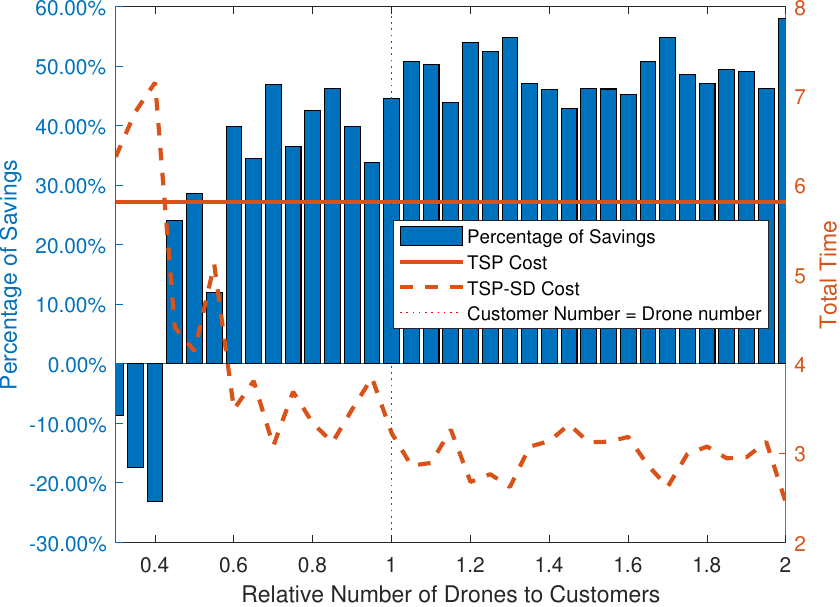}
        \caption{$\rho_1 = 2$}
        \label{fig:SA2rho1Two}
    \end{subfigure}
    \caption{Sensitivity analysis of the results of Algorithm \ref{alg:RouteFinderMC3} with respect to the relative number of drones.}
    \label{fig:SA2}
\end{figure}

From \cref{fig:SA2rho1One}, we can find that if the speed of the drones is equal to the speed of the truck, for the TSP-SD model it is sufficient to have $\rho_2 \geq 0.6$ to consistently outperform the pure TSP model. Increasing the relative speed of drones to 2 (\cref{fig:SA2rho1Two}) both increases the magnitude of savings and the enlarges the range of $\rho_2$ that deliver positive savings ($\rho_2 \geq 0.45$). The same increasing and concave pattern can be observed in $\rho_2$ as well. Comparing \cref{fig:SA2} to \cref{fig:SA1} we see that the sensitivity of the performance of TSP-SD to relative speed of drones is much more than its sensitivity to relative number of drones. It is reasonable to see that the number of drones is not as important as the speed of drones. This is because, if drones are fast enough, we only need a few drones in that area to make the performance of the model with drones better than the normal TSP model, while having huge number of slow moving drones may still force us to prefer truck over the drone on many of the deliveries. 

\subsection{Combining \texorpdfstring{$\rho_1$ and $\rho_2$}{}}
\label{sec:CombiningRhos}
In this section, we try to combine the analysis of $\rho_1$ and $\rho_2$ together to see how much we can improve the delivery performance by including the crowd-owned drones when we consider both parameters at the same time, i.e., comparing TSP-SD with TSP. We let the parameters to vary in the ranges $0.75\leq \rho_1 \leq 3$ and $0.3 \leq \rho_2 \leq 2$. 

\begin{figure}[!h]
    \centering
    \begin{subfigure}[b]{0.35\linewidth}
        \includegraphics[width=\linewidth]{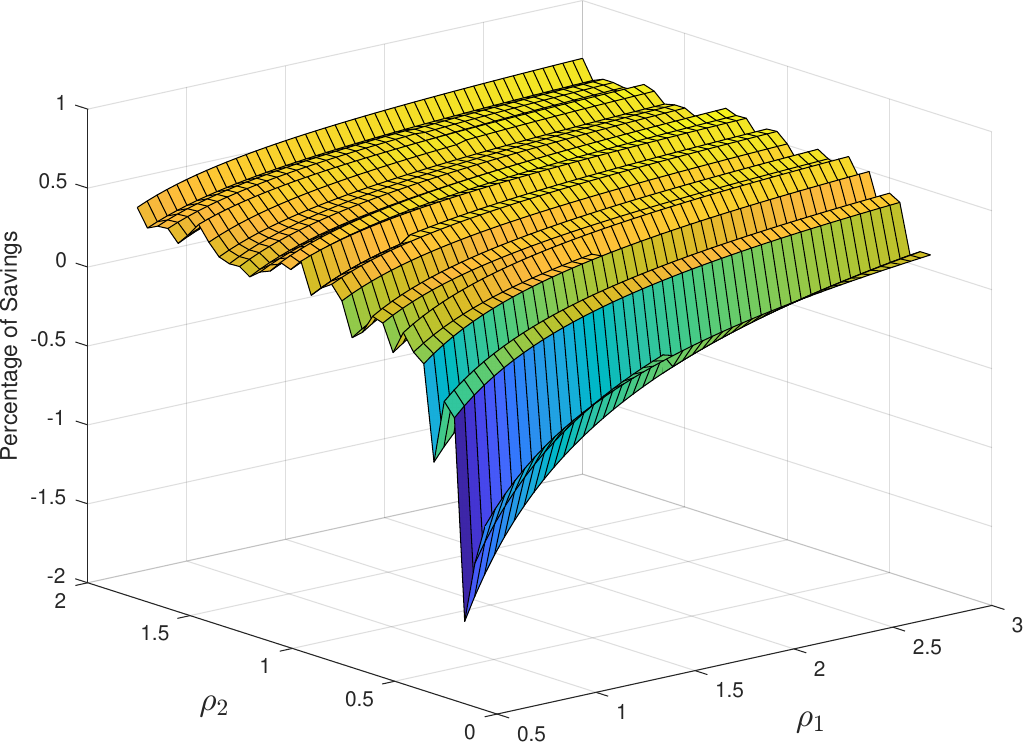}
        \caption{}
         \label{fig:SA3CompleteRange}
    \end{subfigure}
    \qquad
    \begin{subfigure}[b]{0.35\linewidth}
        \includegraphics[width=\linewidth]{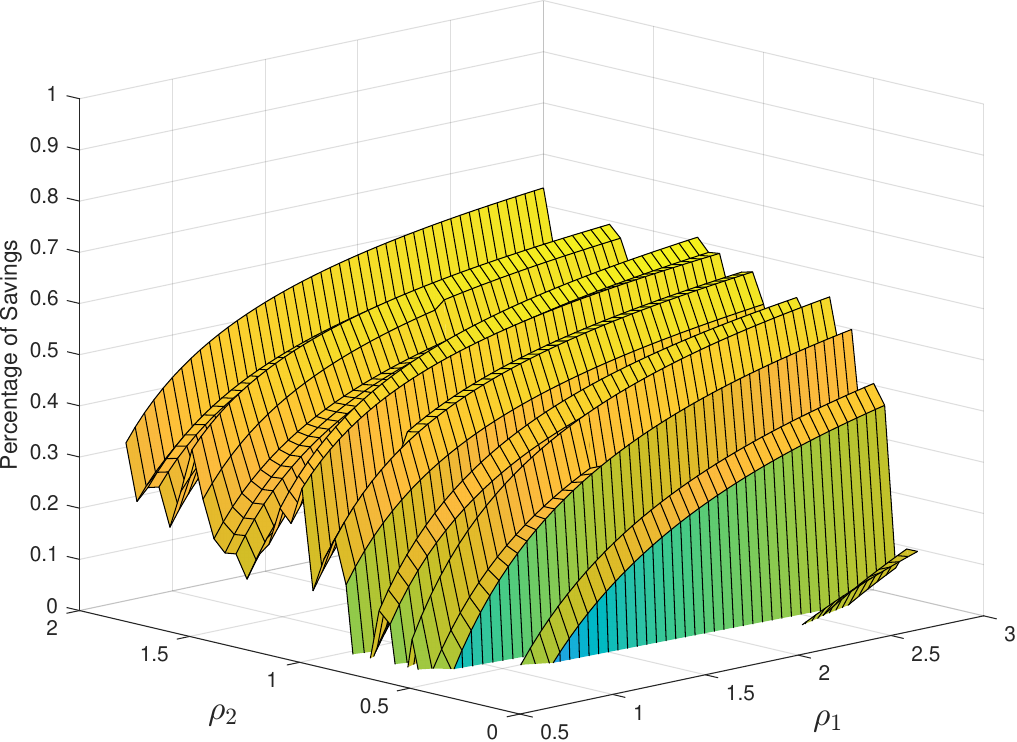}
        \caption{}
         \label{fig:SA3PositiveRange}
    \end{subfigure}
    \caption{The sensitivity analysis of the results of Algorithm \ref{alg:RouteFinderMC3} with respect to the relative speed of drones combined with the relative number of drones. The positive savings in (\subref{fig:SA3CompleteRange}) are shown in (\subref{fig:SA3PositiveRange}).}
    \label{fig:SA3}
\end{figure}

From \cref{fig:SA3CompleteRange}, we can see that if both $\rho_1$ and $\rho_2$ are small, the normal TSP will have much better performance than the model with sharing drones (TSP-SD), which we can also observe from the two previous sections, especially in  \cref{fig:SA2rho1Two}. On the other hand, from \cref{fig:SA3PositiveRange}, which shows the range with positive savings, we can see that the performance of the TSP-SD model improves as $\rho_1$ and $\rho_2$ increase. It reaches its maximum percentage of savings around $63.47\%$ when both $\rho_1$ and $\rho_2$ are around our set upper bounds of 3 and 2, respectively, implying that the improvement could increase further by relaxing this upper end of the range. Meanwhile, the growth rate becomes smaller as $\rho_1$ and $\rho_2$ grow, which is apparent by the concavity of the curve and suggests a diminishing return phenomena. We expect even bigger savings for the revisiting model but the pattern stays the same. For the other centralized alternative (TSPD) it is shown that for $\rho_1 \leq 2$ the gain in delivery time by using a drone is insignificant \cite{Ferrandez2016optimization,agatz2018optimization}. However, our shared drone delivery system shows considerable gain for a wider range of $\rho_1$ that makes our model more applicable across a variety of available drone technologies and regulatory environments.

\subsection{Population Distribution}
\label{sec:PopulationDistribution}

\begin{figure}[h!]
    \centering
    \begin{subfigure}[b]{0.31\linewidth}
        \includegraphics[width=\linewidth]{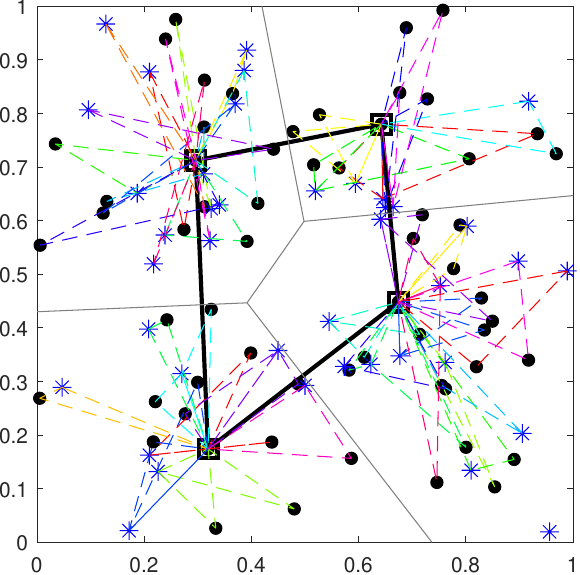}
        \caption{One Population Center with $\sigma^2$ = 0.1}
    \end{subfigure}
    \qquad
    \begin{subfigure}[b]{0.31\linewidth}
        \includegraphics[width=\linewidth]{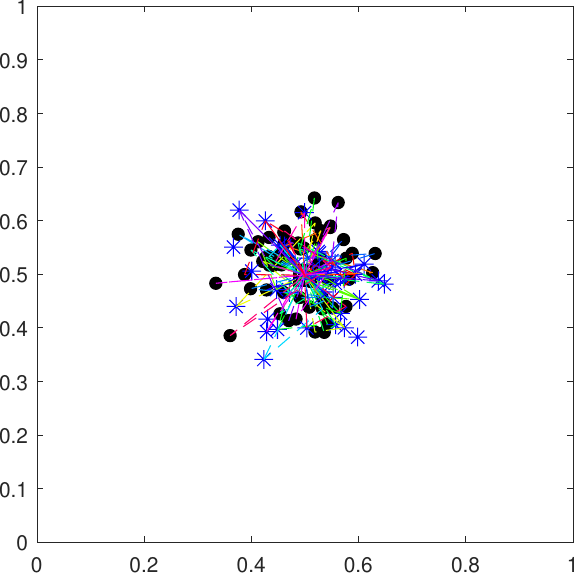}
        \caption{One Population Center with $\sigma^2$ = 0.005}
    \end{subfigure}
    \\
    \begin{subfigure}[b]{0.31\linewidth}
        \includegraphics[width=\linewidth]{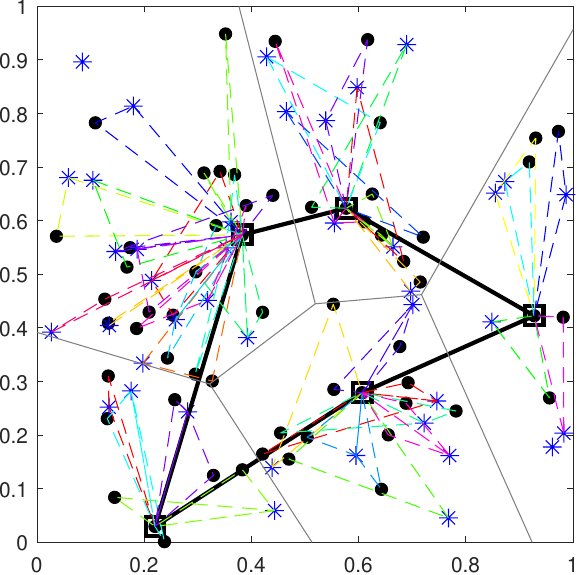}
        \caption{Four Population Centers with $\sigma^2$ = 0.1}
    \end{subfigure}
    \qquad
    \begin{subfigure}[b]{0.31\linewidth}
        \includegraphics[width=\linewidth]{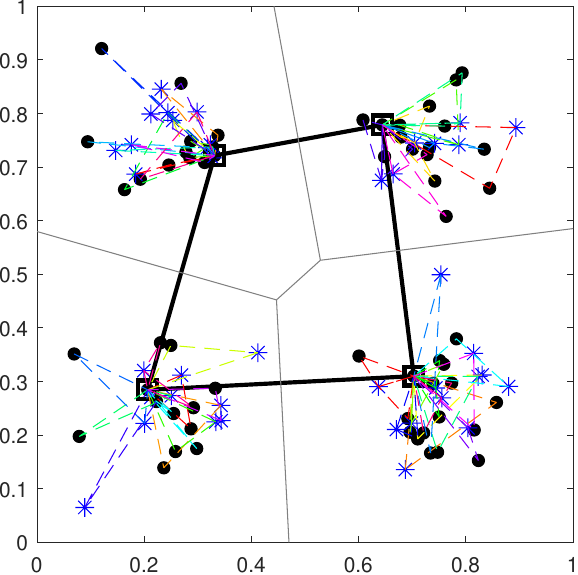}
        \caption{Four Population Centers with $\sigma^2$ = 0.005}
    \end{subfigure}
    \\
    \begin{subfigure}[b]{0.3\linewidth}
    \includegraphics[width=\linewidth]{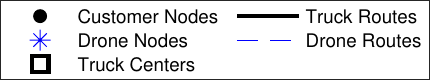}
    \end{subfigure}
    \caption{Two cases of population distribution with two values of $\sigma^2$ and the solutions of Algorithm \ref{alg:RouteFinderMC3} for them. 
    We see that having a clustered population may make a TSP tour a necessity in our solution regardless of how large drone range is relative of the variance. However, we can expect shorter TSP tour when variance is smaller.}
    \label{fig:customer_distribution}
\end{figure}

\begin{figure}[h!]
    \centering
    \begin{subfigure}[b]{0.35\linewidth}
        \includegraphics[width=\linewidth]{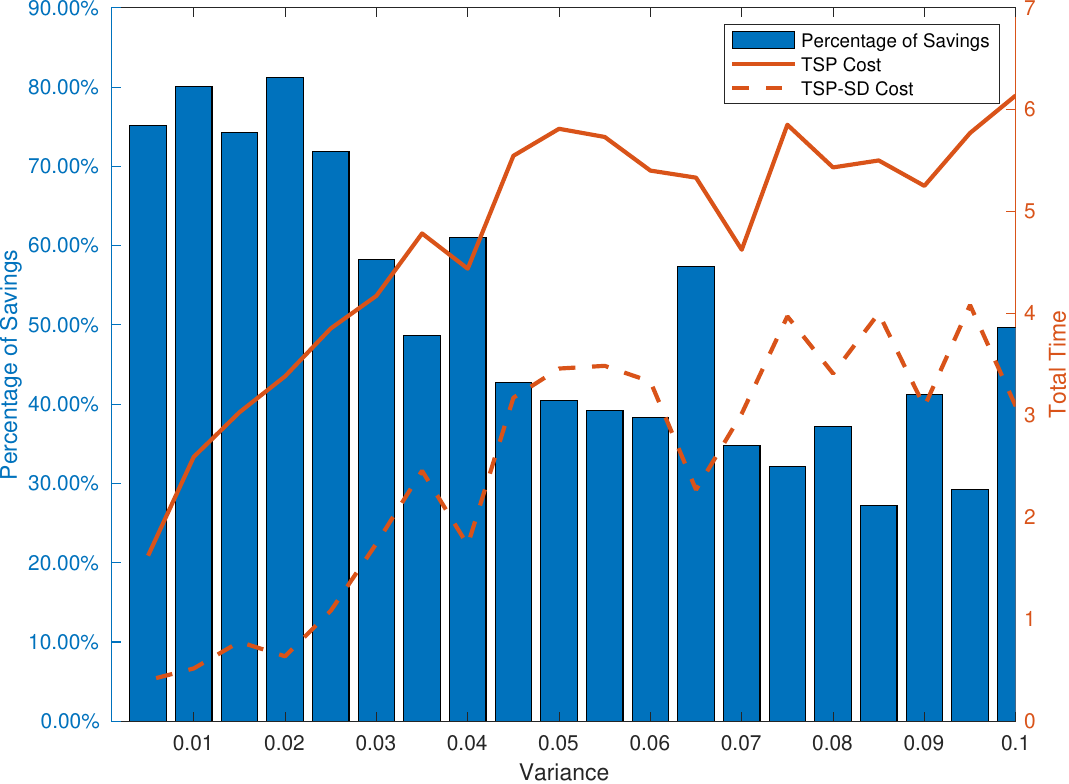}
        \caption{One Population Center}
        \label{fig:SA4OneCenter}
    \end{subfigure}
    \qquad
    \begin{subfigure}[b]{0.35\linewidth}
        \includegraphics[width=\linewidth]{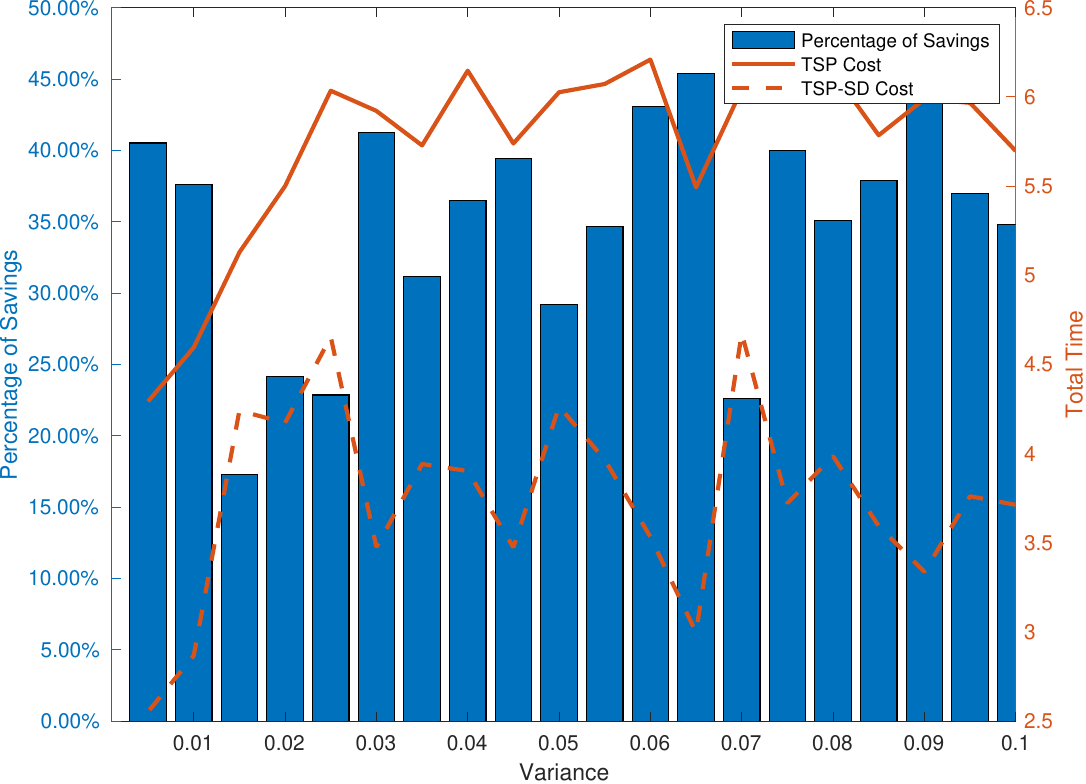}
        \caption{Four Population Centers}
         \label{fig:SA4FourCenter}
    \end{subfigure}
    \caption{Comparison of TSP and TSP-SD models for two cases of population distribution for varying values of variance.}
    \label{fig:SA4}
\end{figure}

We consider two cases with different population distributions. The first case has one population center of 60 customers and 40 drone service providers distributed inside a unit box according to a multivariate Gaussian distribution with mean $\mu=(0.5,0.5)$ and covariance matrices $\Sigma=\usebox{\smlmats}$, where $\sigma^2 \in [0.005,0.1]$ with increments of 0.005.
The second case is a clustered data with four population centers of 60 customers and 40 drone providers (across the clusters) distributed inside a unit box according to an even mixture of four truncated multivariate Gaussian distributions with means $\mu_{1},\mu_{2},\mu_{3},\mu_{4}=(0.25.0.25),(0.25,0.75),(0.75,0.75),(0.25,0.75)$ and covariance matrices $\Sigma_{1}=\Sigma_{2}=\Sigma_{3}=\Sigma_{4}=\usebox{\smlmats}$, where $\sigma^2 \in [0.005,0.1]$ with increments of 0.005. In both cases we assume $L=0.8$ and $\rho_1 = 2$ and compare the results of Algorithm \ref{alg:RouteFinderMC3} with the results of LKH for pure TSP for different values of $\sigma^2$. The \cref{fig:customer_distribution} shows examples of generated instances with the solutions found by the recharging version of our Algorithm \ref{alg:RouteFinderMC3}. Complete results for these instances are presented in \cref{fig:SA4}.

We can see from \cref{fig:SA4} that the proposed combined model of truck and sharing economy drones (TSP-SD) consistently and significantly performs better than the pure TSP (truck-only) model in both cases and across the range of the variance. Its best performance usually happens when customers are concentered in one place, i.e., very small variance, whether we have one center or several centers of population. This is consistent with real-life applications as in the metropolitan areas population distribution is usually nonuniform and residences usually live near each other in a neighborhood or town while having several neighborhoods or towns in the metropolitan area with some distances from each other.
When the variance is low, we can see savings of up to $80\%$ in the delivery time, when compared to TSP, by using drones in the once center case (\cref{fig:SA4OneCenter}). Although this savings becomes smaller when the distribution of the customers gets closer to a uniforms distribution, we can still expect at least $30\%$ savings in the delivery time. 
When there are four centers of population in a large city, we can still see between $17\%$ and $45\%$ savings in the delivery time (\cref{fig:SA4FourCenter}). Unlike the one population center case that shows a negative correlation between the savings and $\sigma^2$, in four population centers case we see that the savings stays more or less flat that is due to the clustered nature of the population. 
These observations on the sensitivity of our proposed model to the customer distribution is consistent with findings of \cite{chen2018multi} that suggests that the benefits of the crowdsourced delivery depends on the spatial characteristics of the network (of pickup/delivery/transfer locations and roads).

\subsection{Impact of Moving Cluster Centers and Merging Centers}
\label{sec:MovingCentersImpact}
\begin{figure}[h!]
    \centering
    \begin{subfigure}[b]{0.3\linewidth}
        \includegraphics[width=\linewidth]{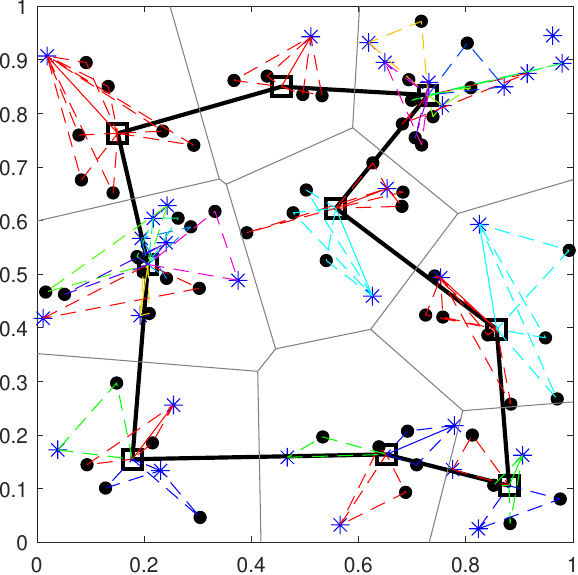}
        \caption{Without Moving Centers}
        \label{fig:SA_OHMC3_OH}
    \end{subfigure}
    \qquad
    \begin{subfigure}[b]{0.3\linewidth}
        \includegraphics[width=\linewidth]{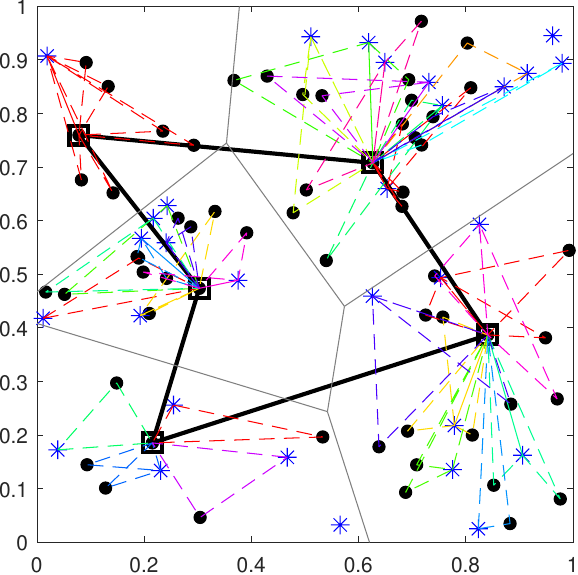}
        \caption{After Moving Centers}
        \label{fig:SA_OHMC3_MC3}
    \end{subfigure}
    \\
    \begin{subfigure}[b]{0.3\linewidth}
    \includegraphics[width=\linewidth]{legend_final.pdf}
    \end{subfigure}
    \caption{A comparison of solution of our original heuristic, presented in \cref{sec:algorithm}, in (\subref{fig:SA_OHMC3_OH}) and the solution obtained by  \cref{alg:RouteFinderMC3} after moving cluster centers to customer nodes and merging centers in (\subref{fig:SA_OHMC3_MC3}). In this instance we have 60 customer nodes, 30 drone nodes ($\rho_2=0.5$), distributed uniformly at random in a unit box with $\rho_1=3$ and $L=0.8$.}
    \label{fig:MCE}
\end{figure}

Here, we compare the performance of our original heuristic, presented in \cref{sec:algorithm} in which the cluster centers (truck nodes) could be anywhere in the continuous space, with the final improved   \cref{alg:RouteFinderMC2}, and 
Algorithm \ref{alg:RouteFinderMC3} discussed in \cref{sec:Improvements} in which the cluster centers (truck nodes) are located at some customer nodes.
As mentioned in Sections \ref{subsec:MovingCentersNearestCustomerNode} and \ref{subsec:MovingCentersNearestCustomerNode2AreaCenter}, moving a cluster center to customer nodes needs to satisfy one condition that is the new 
center needs to be able to still cover all of its assigned customers using its assigned drones. When merging the centers, as discussed in \cref{subsec:MergingCenters}, this feasibility is also satisfied in another way by reallocation of customers. \cref{fig:MCE} compares an example solution for before and after moving the cluster centers customer nodes (i.e., the original heuristic vs. \cref{alg:RouteFinderMC3}).

\begin{figure}[h!]
    \centering
    \begin{subfigure}[b]{0.33\linewidth}
        \includegraphics[width=\linewidth]{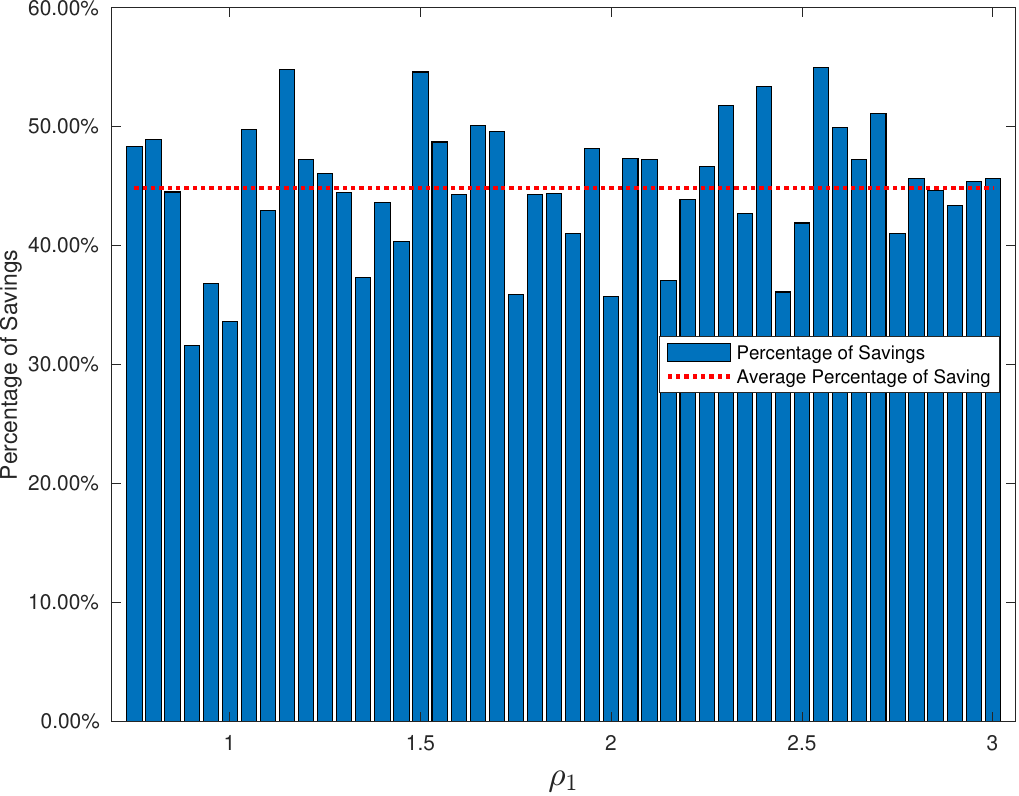}
        \caption{$\rho_2=0.5$}
        \label{fig:SA_MCrho2Half}
    \end{subfigure}
    \begin{subfigure}[b]{0.33\linewidth}
        \includegraphics[width=\linewidth]{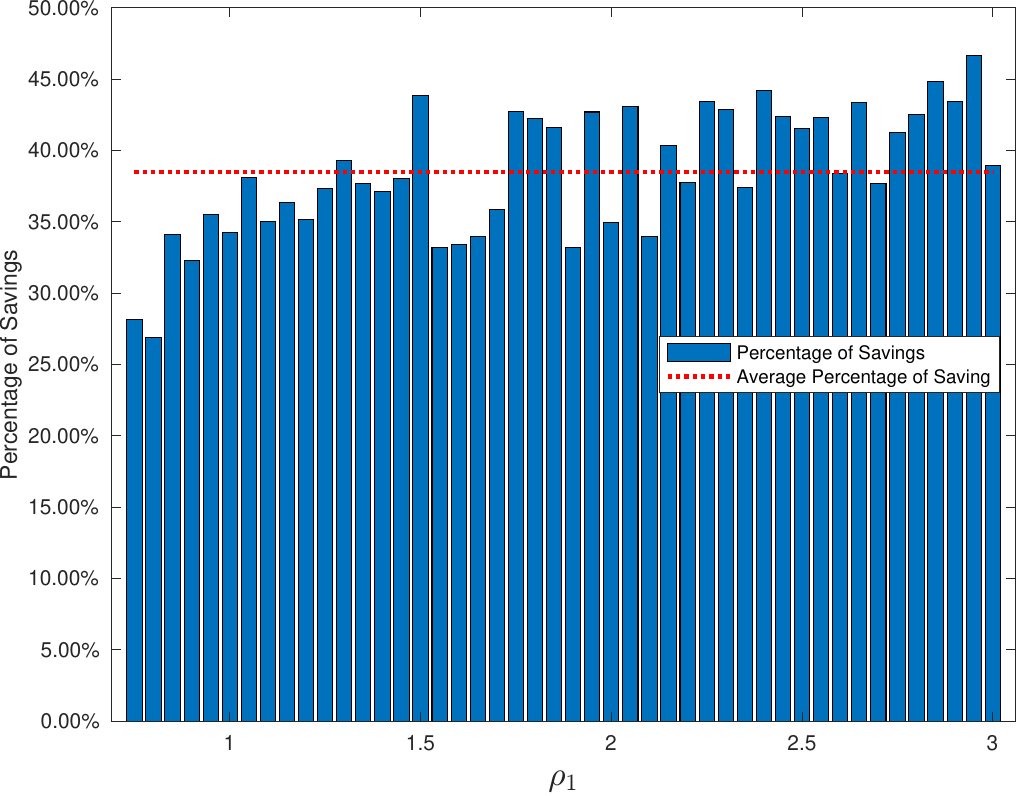}
        \caption{$\rho_2=2$}
        \label{fig:SA_MCrho2Two}
    \end{subfigure}
    \begin{subfigure}[b]{0.33\linewidth}
        \includegraphics[width=\linewidth]{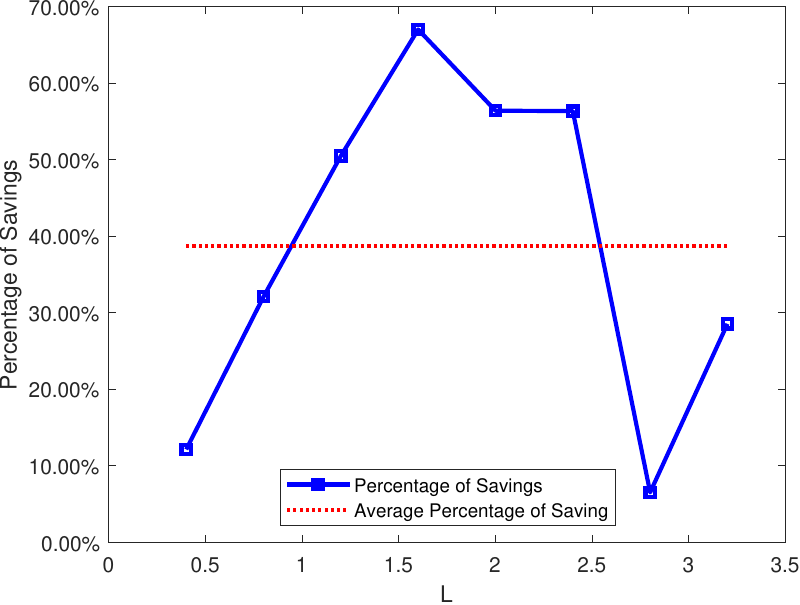}
        \caption{$\rho_1=2,\;\rho_2=1$}
        \label{fig:SA_MC_L}
    \end{subfigure}
    \caption{Sensitivity analysis of the impact of moving cluster centers to customer nodes via \cref{alg:RouteFinderMC1} and \cref{alg:RouteFinderMC2} and further merging them together in \cref{alg:RouteFinderMC3} on the performance when compared to the original heuristic presented in \cref{sec:algorithm}. (\subref{fig:SA_MCrho2Half}) shows the results for $\rho_2=0.5$ and (\subref{fig:SA_MCrho2Two}) presents the results for $\rho_2=2$. The range of relative speed in both cases is $0.75 \leq \rho_1 \leq 3$. The sensitivity of the impact of this improvement with respect to drone range in illustrated in \subref{fig:SA_MC_L} which shows the maximum improvement for middle range $L$.}
    \label{fig:MC}
\end{figure}

To evaluate the overall impact of this change, we use assume $0.75 \leq \rho_1 \leq 3$ and $\rho_2 \in \{0.5,2\}$ for our simulation. The results of comparison between \cref{alg:RouteFinderMC3} and our original heuristic of \cref{sec:algorithm} is presented in \cref{fig:MC}. In Figures \ref{fig:SA_MCrho2Half} and \ref{fig:SA_MCrho2Two} We see about 45\% and 38\% savings on average for two values of $\rho_2$, respectively, which suggest a slight decrease in savings by increasing $\rho_2$. This is because when we have an abundance of drone supply the original heuristic could also perform well making the gap between the two algorithms narrower. However, this average decrease is quite small. The savings are more or less flat across different values of $\rho$ as well. This suggests the robustness of the implemented improvements with respect to $\rho_1$ and $\rho_2$, providing more than roughly 40\% savings in the delivery time when compared to the algorithm before improvement. However, the design of \cref{alg:RouteFinderMC3} hints more sensitivity with respect to drone range $L$ for the purpose of this comparison. This analysis for increments of 0.4 in $L$ is shown in \cref{fig:SA_MC_L}. We see the biggest improvement (up to 68\%) for the middle range $L$ where we needed it the most as already seen in Tables \ref{tab:RechargingResultsModelAlgComp} and \ref{tab:RevisitingResultsModelAlgComp}. For very large $L$ both problems can result to one truck node (serving all customers with drones). If the location of this center is the same in the solutions obtained by these two algorithms, the only remaining factor making a difference in the solution is the drone routes obtained by the tabu search step. This explains the small gap for $L=2.8$.

\subsection{Revisiting versus Recharging}
\label{sec:batteryUtilizationImpact}
\begin{figure}[t!]
    \centering
    \begin{subfigure}[b]{0.35\linewidth}
        \includegraphics[width=\linewidth]{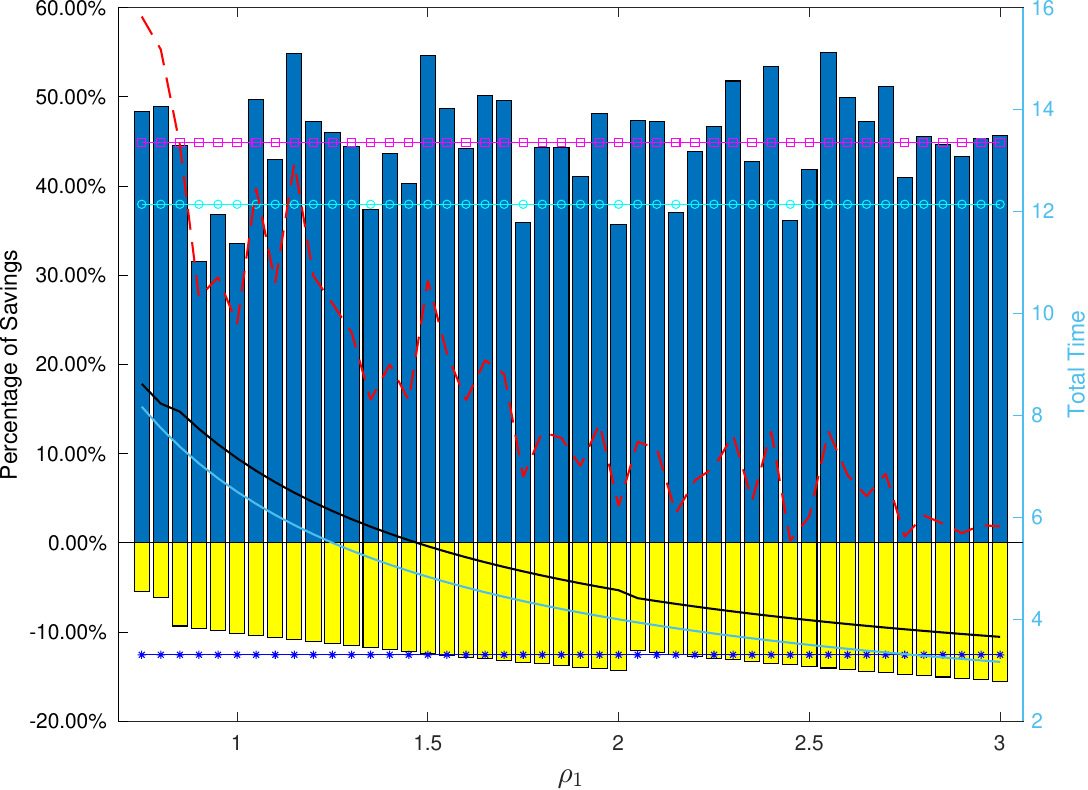}
        \caption{Savings analysis when $\rho_2=0.5$}
         \label{fig:buRho2Half}
    \end{subfigure}
    \qquad
    \begin{subfigure}[b]{0.35\linewidth}
        \includegraphics[width=\linewidth]{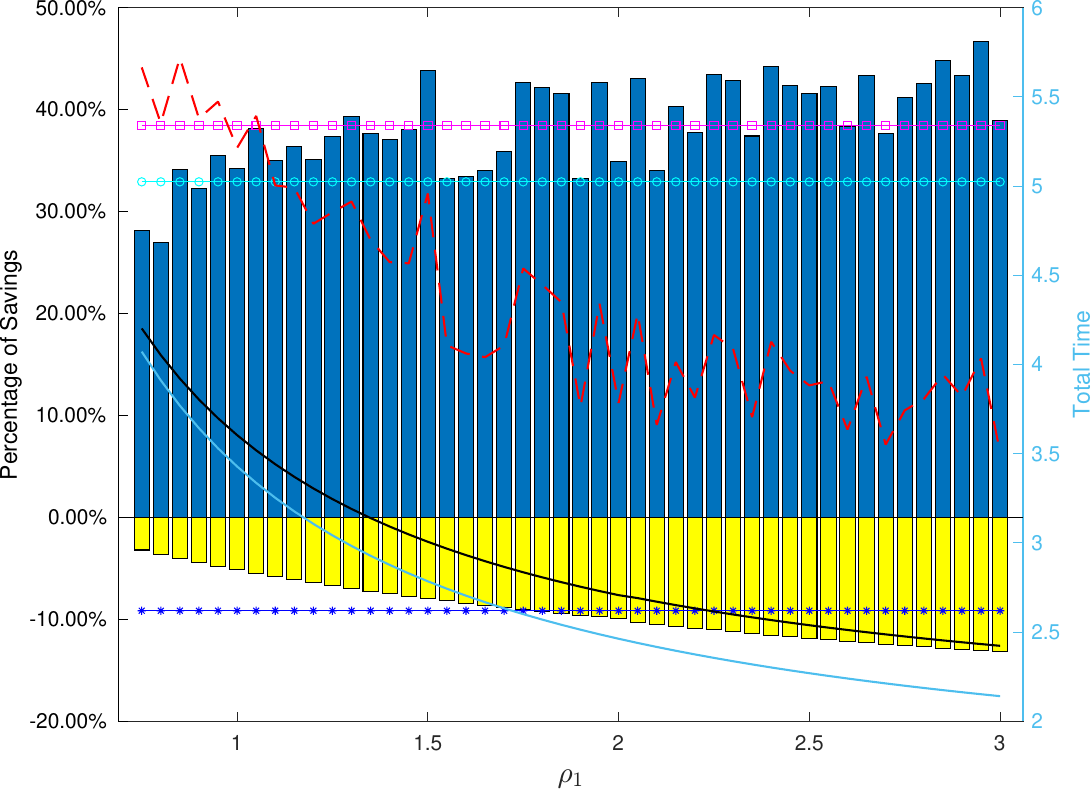}
        \caption{Savings analysis when $\rho_2=2$}
         \label{fig:buRho2Two}
    \end{subfigure}
        \\
    \begin{subfigure}[b]{0.7\linewidth}
    \includegraphics[width=\linewidth]{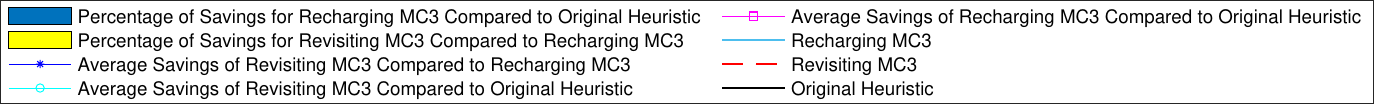}
    \end{subfigure}
    \caption{Illustration of the impact of battery utilization modification that shows savings, for various values of $\rho_1$ with $\rho_2=\{0.5,2\}$, by comparison of three algorithms including recharging and revisiting versions of \cref{alg:RouteFinderMC3} and the original heuristic of  \cref{sec:algorithm}.}
    \label{fig:BU}
\end{figure}

Comparing the revisiting and recharging modes using the results of the algorithm is tricky. When using the proposed algorithm we could see revisiting model as a battery utilization modification to the algorithm that allows the drones to revisit the truck node after a delivery to pick up another package if their battery limits allow that. However, one should keep in mind that achieving any improvement by this modification highly depends on the problem instance and the input parameters. When comparing the results of the optimization models we expect to see further improvement by allowing revisiting for the drones that was identified in \cref{tab:RevisitingResultsModelAlgComp} in several instances. However, comparing the results of \cref{alg:RouteFinderMC3} in Tables \ref{tab:RechargingResultsModelAlgComp} and \ref{tab:RevisitingResultsModelAlgComp} we can see a mixed situation and in some instances the recharging model is providing a better solution. This is due to the additional complexity of the revisiting model that makes it much harder for the algorithm to find better solutions for certain instances of the problem, resulting in a larger optimality gap and potentially weaker solution than the recharging model. In other words, we are dealing with a larger feasibility set potentially containing a better optimal solution but depending on the complexity of the problem instance we may end up with a better or worse solution than the output of the recharging model.
 \cref{fig:BU} shows an example of implementing this modification that results to weaker solutions for the revisiting model. When revisiting is desired we suggest to run both recharging and revisiting versions of the algorithm and pick the best. For smaller size problems we can obviously solve the problem with the optimization model of Problem IV to expect better results. Further adjustment of the parameters of the algorithm to better handle the complexity of the revisiting model for such unfriendly instances may also help to achieve positive savings compared to the recharging version of \cref{alg:RouteFinderMC3}. 

\section{The Impact on Carbon Footprint and Comparison with TSPD}
\label{sec:CarbonFootprint}
It is clear that drone delivery systems such as our proposed combined system and the coordinated but centralized system of a truck and a drone in which the drone is controlled from the truck (the TSPD model\footnote{We follow the name used by \cite{agatz2018optimization}. This is also called horsefly problem, truck-and-drone TSP, etc.}) have less carbon footprint than a truck-only (pure TSP) delivery system since drones are much more energy-efficient than trucks. 

In order to evaluate the energy efficiency of our proposed shared drone delivery model, we compare it with the TSPD model. First, note that if in both models we assume the drones are capable of carrying more than one package, the shared delivery model would have less carbon footprint since it becomes very similar to the comparison of the \emph{generalized traveling salesman problem} (GTSP) and the TSP done in the earlier paper \cite{carlsson2016household} that showed the household-level transportation has its own economies of scale and most likely outperforms a centralized delivery system. 

\begin{figure}[t!]
    \centering
     \begin{subfigure}[b]{0.27\linewidth}
        \includegraphics[width=\linewidth]{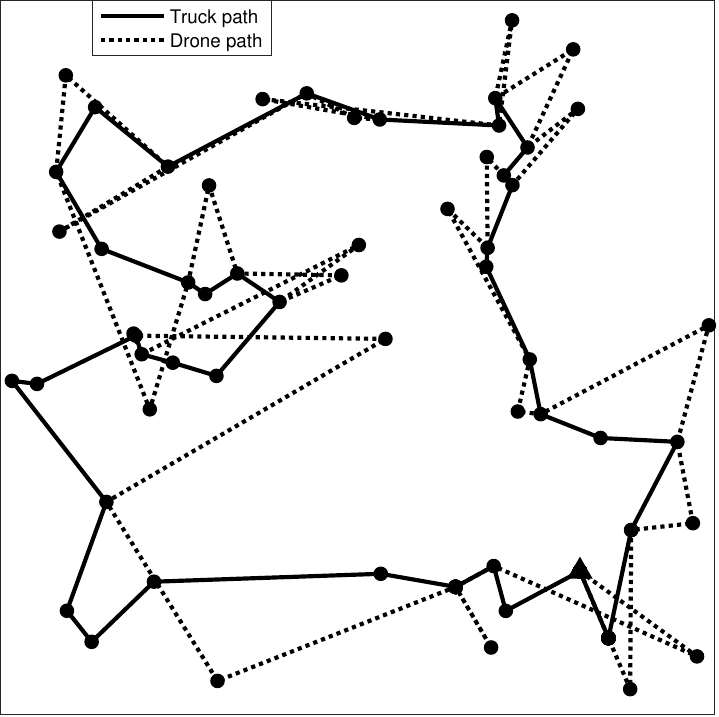}
        \caption{Solution of the TSPD model \\
        					~}
	\label{fig:TSPDvsSharing-TSPDModel}
    \end{subfigure}
        \quad
    \begin{subfigure}[b]{0.27\linewidth}
        \includegraphics[width=\linewidth]{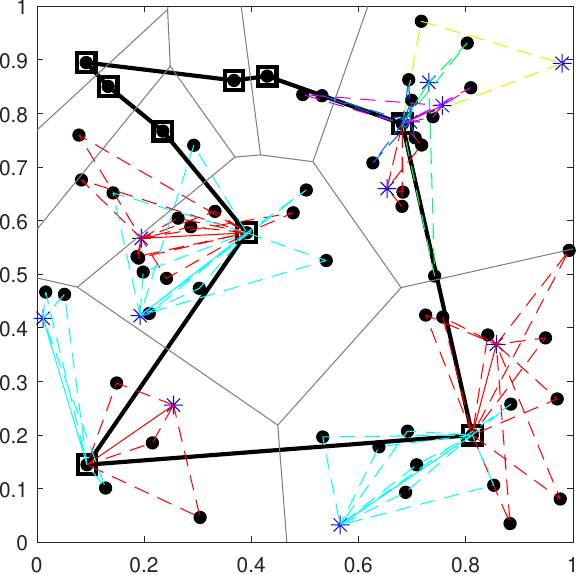}
        \caption{Solution of the TSP-SD model with one drone in each cluster}
        \label{fig:TSPDvsSharing-SharingModel}
    \end{subfigure}
    \quad
    \begin{subfigure}[b]{0.37\linewidth}
        \includegraphics[width=\linewidth]{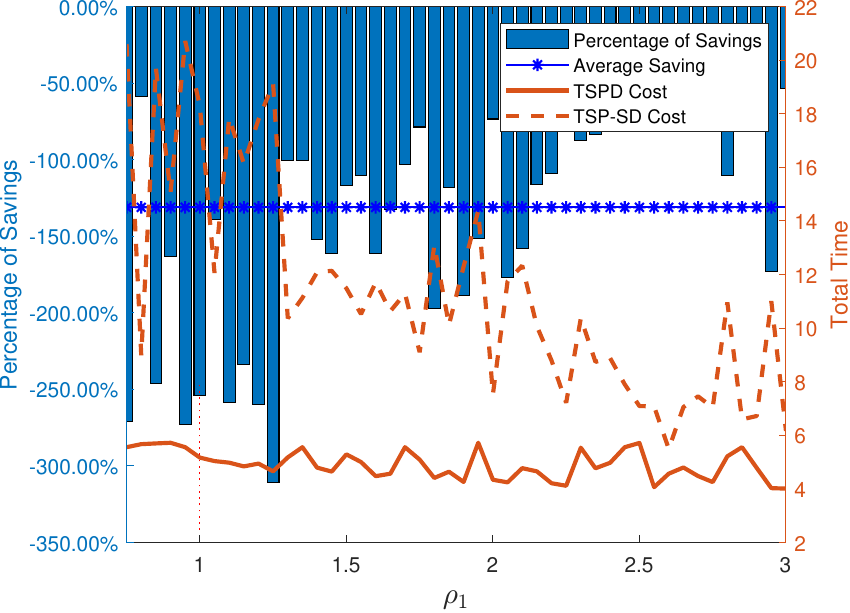}
        \caption{Comparison of TSPD model and TSP-SD model with one drone in each cluster ($\rho_2=0.18$)}
        \label{fig:TSPDvsSharing-Savings}
    \end{subfigure}
    \caption{Comparison of the TSPD model the shared delivery model (TSP-SD) assuming only one drone in each cluster of the shared model.}
    \label{fig:TSPDvsSharing}
\end{figure}

\begin{figure}[t!]
    \centering
    \begin{subfigure}[b]{0.32\linewidth}
        \includegraphics[width=\linewidth]{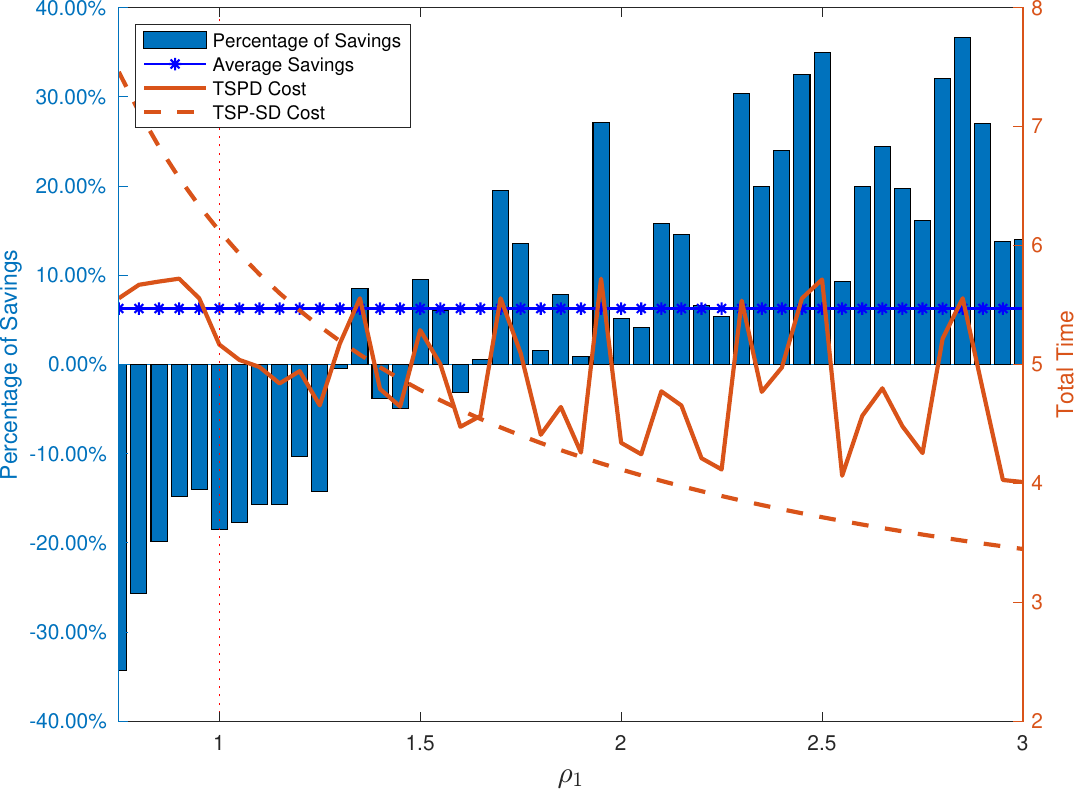}
        \caption{$\rho_2=0.5$}
        \label{fig:TSPDvsSharing-Savings-rhoTwoHalf}
    \end{subfigure}
    \begin{subfigure}[b]{0.32\linewidth}
        \includegraphics[width=\linewidth]{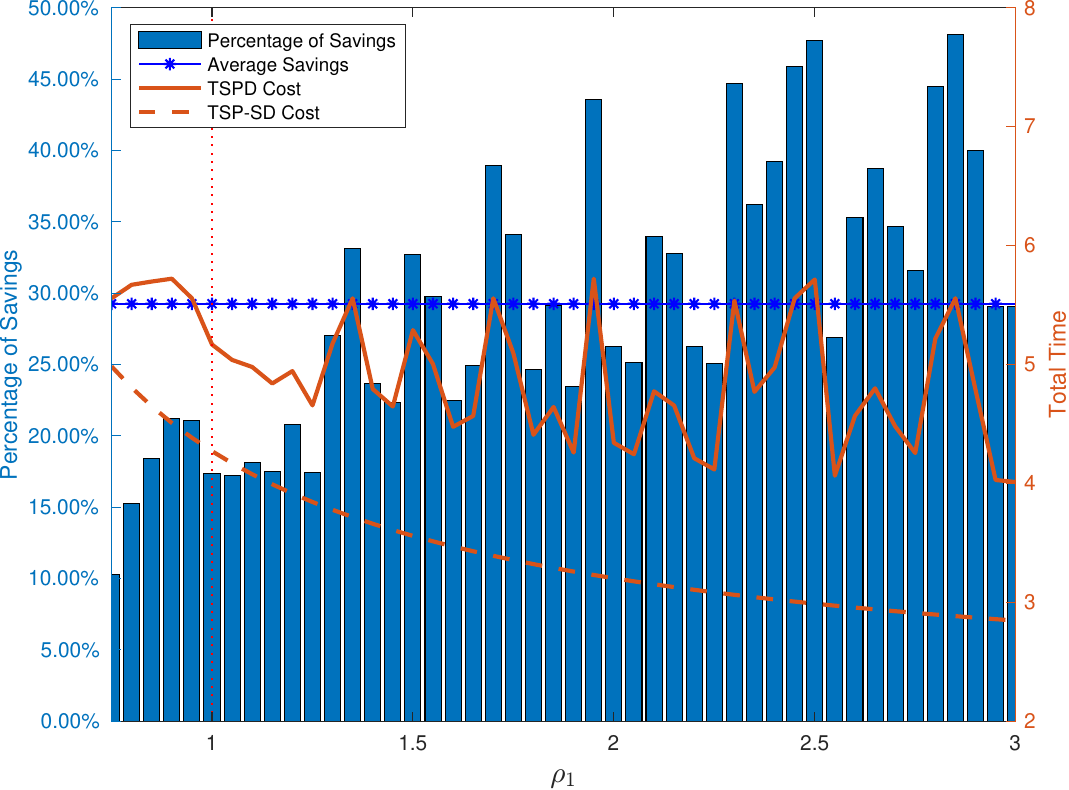}
        \caption{$\rho_2=1$}
        \label{fig:TSPDvsSharing-Savings-rhoTwoOne}
    \end{subfigure}
    \begin{subfigure}[b]{0.32\linewidth}
        \includegraphics[width=\linewidth]{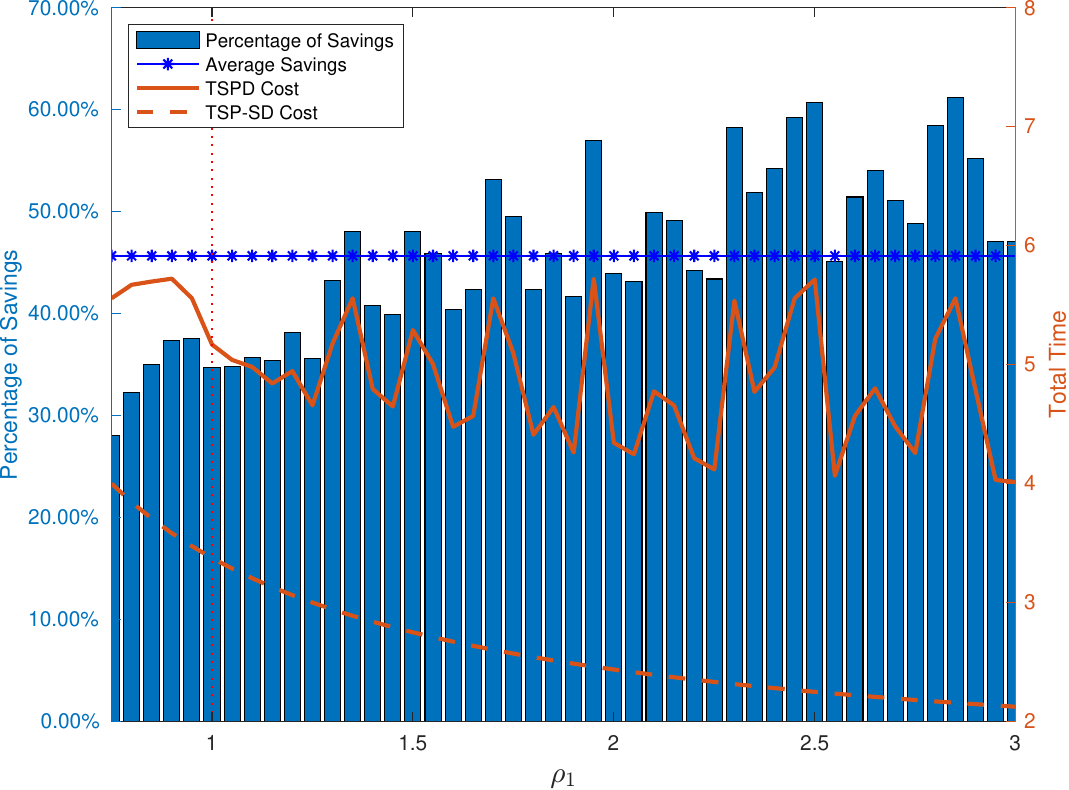}
        \caption{$\rho_2=2$}
        \label{fig:TSPDvsSharing-Savings-rhoTwoTwo}
    \end{subfigure}
    \caption{Comparison of the TSPD model the shared delivery model (TSP-SD)  for different values of relative number and speed of drones.}
    \label{fig:TSPD-Sharing-IncreasingRho2}
\end{figure}

Now assume the drones can only carry one package per trip (which is the case in most current delivery drone technologies) and have to return to their home base (for replacing battery or recharging) after each delivery. The comparison would obviously depend on the balance between the supply of drones and demand in deliveries. The \cref{fig:TSPDvsSharing} shows a comparison between the two systems in two situations when there is only one drone in each cluster ($\rho_2=0.18$ for this instance with 60 customers), i.e., extreme shortage of drone supply, and \cref{fig:TSPD-Sharing-IncreasingRho2} illustrates the results of this comparison when we increase drone supply to more reasonable levels in sharing economy environments, which as we will discuss next makes this comparison more fair as well. We use the ratio $\frac{\text{Time cost of TSPD - Time cost of TSP-SD}}{\text{Time cost of TSPD}}$ to calculate the percentage of savings.  

For the TSPD model we assumed that drone can pickup a package at one stopping point and returns to the truck after delivering it at either the same point or at another stopping point which makes the model more flexible and more efficient. We also assume no restriction of the drone range. Both of these assumptions give an advantage to this model over our TSP-SD. We solved this problem via the algorithms developed in our working papers \cite{ma2023lastmilehorsefly,mohammadi2023lastmilehorsefly} and picked the best result.
Figures \ref{fig:TSPDvsSharing-TSPDModel} and \ref{fig:TSPDvsSharing-SharingModel} show the solution of the TSPD model and the sharing model TSP-SD, respectively, for a problem with 60 customers distributed uniformly at random in a unit box. We have set $\rho_1=3$ and $L=0.8$ for both models. Using the binary search and $k$-means step of \cref{alg:RouteFinderMC1} to find a feasible $k$ we found that $k=11$ cluster centers are needed. The number of clusters is further reduced by \cref{alg:RouteFinderMC3} to 9 via merging centers as visible in \cref{fig:TSPDvsSharing-SharingModel}. Assuming one drone in each of the original 11 clusters gives $\rho_2 = 0.18$. 
 It is clear from \cref{fig:TSPDvsSharing-Savings} that with this restrictive assumption the TSPD model will have much better performance than our TSP-SD model. This is understandable as in the shared model, the truck needs to wait at the cluster centers for one drone to serve all the customers which looks very inefficient when drone supply is limited. The drone-traveling time is huge in this circumstance and the truck idles when it could serve customers. In contrast, in the TSPD model the truck can serve another customer while the drone is delivering the package as well. It is obvious that the TSPD model is faster than the shared model if there is only one drone in each cluster since there is no such huge waiting time for the truck. TSPD can be seen as a model that has one drone for each customer in the neighborhood. To see this more clearly, imagine a very large relative drone speed, i.e., letting $\rho_1 \rightarrow \infty$. In such extreme case the truck will stop at some point in the region and all the customers will be served by the drone, forming a star-shape combined route. This situation is equivalent to a situation where each customer node sends its own drone to the truck to pick up its package and bring it back to customer's home, i.e., $\rho_2=1$. The sharing economy model assumes reasonable drone supply (such as a perfect balance of supply-demand with $\rho_2=1$) with limited relative drone speed $\rho_1$. Obviously, this situation justifies using the truck for serving some of the customers and leaves the question of finding the optimal trade-off. Therefore, to have a fair comparison we have to increase the relative number of drones in the sharing model. Note that the assumption of unlimited drone range ($L \rightarrow \infty$) for the TSPD model also gives a significant advantage to that model. This is one of the reasons behind the fact that TSPD is showing less sensitivity to $\rho_1$ in its considered range as can be observed in \cref{fig:TSPDvsSharing-Savings} and Figures \ref{fig:TSPDvsSharing-Savings-rhoTwoHalf} -- \subref{fig:TSPDvsSharing-Savings-rhoTwoTwo}. Of course, this will change when $\rho_1$ increases significantly beyond this range.

To investigate the impact of the balance between the supply of drones and the demand in deliveries, we did a sensitivity analysis in terms of the number of drones. In \cref{fig:TSPD-Sharing-IncreasingRho2}, it is clear that when the relative number of drones increases to half of the customers, the average savings of the sharing model can be as much as about 6\%, an increase by 137 percentage points compared to the case with one drone in each cluster. We see in \cref{fig:TSPDvsSharing-Savings-rhoTwoHalf} that from $\rho_1=1.35$ we start to see positive savings in this case.
Furthermore, when there is the same number of drones and customers in the neighborhood ($\rho_2 = 1$), the sharing model will have better performance compared to the TSPD model on all values of $\rho_1$, and can achieve about $30\%$ savings on average. If we further increase the drone supply to $\rho_2 = 2$ there is more than $46\%$ savings on average in the delivery time. This is consistent with other sharing economy application and shows that to use shared drones in a delivery system, we need to have enough supply (of drones) relative to the demand in deliveries in order to have a better performance than the centralized alternative of coordinated truck and drone delivery system. Clearly, the gap between the TSPD and our TSP-SD model when we allow drones to revisit the pickup centers becomes wider, at least when we compare the optimal solutions. Assuming a larger $L$ will also increase this gap significantly.

\section{Conclusion}
\label{sec:Conclusion}
In this paper we have developed a shared last-mile delivery model in which a truck carries packages to a neighborhood and then outsources the last piece of the trip to private drone operators whose service can be utilized on a sharing economy platform. We have developed several optimization models for different versions of the problem and proposed efficient algorithms to solve them. Our computational analysis show that the shared delivery model (decentralized model) is much more efficient than the traditional truck-only delivery model (centralized model) in almost all possible scenarios and offers significant savings in the delivery time. This is aligned with the results from \cite{carlsson2016household}. The comparison between the shared delivery model and the coordinated delivery system, in which a truck carries and controls a drone during the delivery operation, depends on other factors such as number of available drones in the platform, their capacity and speed. Our analysis in this case also shows a considerable overall advantage for the decentralization. For future work, one may look into considering different factors such as time windows for delivery to customers, time windows for drone availability, ability of drones to carry multiple packages at the same time, weight capacity for the drones, and combination of the system with crowdsourced drivers. Another avenue for extending this work is to incorporate pickup stations into the problem. Finally, it would be interesting to consider stochastic models: one example is to assume that some of the drones with some probability could fail to deliver and have to bring the package back to the truck and that customer must be then served by the truck.

\section*{Acknowledgments}
The authors gratefully acknowledge support from a Tier-2 grant from Northeastern University (NU) as well as an NSF Planning Grant for establishing an engineering research center (ERC). Authors would also like to thank NU's SHARE Group, led by Ozlem Ergun, for their support and for creating an environment in which this research could be accomplished.

\bibliographystyle{unsrt}
\bibliography{SharedDroneDelivery_Refs.bib}

\pagebreak{}

\appendix

\part*{Online Supplement to ``Last Mile Delivery with Drones and Sharing Economy''}

\section{Alternative Model for Problem I}
\label{sec:AltModelProblemI}
For the alternative formulation of Problem I, we change some of the settings of the model as detailed in Table \ref{tab:Sets-Parameters-Variables-ProblemI-Alternative}.
\begin{table}[!h]
\protect\caption{Sets, parameters, and variables used in the one center model.}
\label{tab:Sets-Parameters-Variables-ProblemI-Alternative}
\begin{center}
\begin{tabular}{c|p{0.72\textwidth}}  
      \hline
      \hline
      \textbf{Sets/Parameters/Variables} & \textbf{Description} \\
      \hline
      $C$ & Customer Nodes \\
      $D$ & Drone Nodes\\
      $T$ & Truck nodes, i.e., only one node here. We call it node $0$. \\
      \hline
      $d_{ij}$ & Distance between node $i$ and $j$ for $i,j\in C\cup D \cup T$ \\
      $v_D$ & Speed of drones\\
      $L$ & Longest distance a drone can travel without charging battery \\
      \hline
       $x_{ij}$ & $\in \N$. The number of trips from node $i$ to node $j$. \\
       $q_{i}$ & Total travel time of the drone with base at node $i$\\
       $Q$ & Maximum time spent by all drones (makespan) \\
       \hline
   \end{tabular}
\end{center}
\end{table}

Then,  an alternative mixed integer linear programming (MILP) model for this problem can be written as follows: 
  \begin{eqnarray}
   \minimize \qquad  Q & &  \quad \qquad \st \nonumber  \\
    & & \nonumber \\
    \sum_{i\in D} x_{ji} &  = & 1 \,, \qquad \forall j \in C    \label{eq:Problem1AlternativeDroneSatisfyAll} \\
    x_{0i} + \sum_{j\in C} x_{ij} & = & 0 \,, \qquad \forall i \in D      \label{eq:Problem1AlternativeNoSelfNoDirect} \\
    x_{0j} & = & 1 \,, \qquad \forall j \in C       \label{eq:Problem1AlternativeCenterSatisfyAll} \\
    x_{i0} - \sum_{j\in C}x_{ji} & = & 0 \,, \qquad \forall i \in D     \label{eq:Problem1AlternativeDroneNodeFlowBalance} \\
        \sum_{i \in D}x_{i0} - \sum_{j\in C}x_{0j} & = & 0\,, \qquad   \label{eq:Problem1AlternativeCenterFlowBalance} \\
    (d_{i0}+d_{0j}+ d_{ji})x_{ji}  & \leq &  L\,,  \qquad \forall i\in D, \; \forall j\in C 		\label{eq:Problem1AlternativeDroneRange} \\
    \frac{1}{v_D}\sum_{j\in C}(d_{i0}+d_{0j}+ d_{ji})x_{ji} & \leq & q_i \,, \qquad \forall i\in D \label{eq:Problem1AlternativeDroneTime} \\
    q_i & \leq & Q \,, \qquad \forall i\in D \label{eq:Problem1AlternativeTotalTime} \\
    x_{0,j}\, \; x_{ji} & \in & \{0,1\} \,, \qquad \forall i,j\in C\cup D \nonumber \\
    x_{i0} & \in & \N \,, \qquad \forall i \in D \nonumber \\
    q_i\,, \; Q & \geq & 0 \,, \qquad \forall i\in D \nonumber
  \end{eqnarray} 
Constraint \eqref{eq:Problem1AlternativeDroneSatisfyAll} ensures all customers are visited once by a drone. 
Constraints \eqref{eq:Problem1AlternativeNoSelfNoDirect} and \eqref{eq:Problem1AlternativeCenterSatisfyAll} ensure all customers are serviced via the center (truck node). 
Constraint \eqref{eq:Problem1AlternativeDroneNodeFlowBalance} and \eqref{eq:Problem1AlternativeCenterFlowBalance} enforce the balance of outbound and inbound flows at each drone node and at the pickup center.
Constraint \eqref{eq:Problem1AlternativeDroneRange} is the drone range limit.
Constraints \eqref{eq:Problem1AlternativeDroneTime} sums the total time (distance) of all routes starting from each drone node. Constraint \eqref{eq:Problem1AlternativeTotalTime} finds the maximum time (distance) among all drones.

\section{Alternative Model for Problem II}
\label{sec:AltModelProblemII}
For the alternative formulation of Problem I, we change some of the settings of the model as detailed in Table \ref{tab:Sets-Parameters-Variables-ProblemII-Alternative}. We know that each drone can make multiple trips and in each trip it can serve multiple customers. This is obviously much more complicated than the model for Problem I. To do this in an easier way, we hypothetically each trip initiated from a drone node is made with a new drone. Therefore, we index different trips as different drones. Since each drone node would need at most $|C|$ trips to serve all customers we define set $K$ as the set of all hypothetical drones in which the first $|C|$ members are associated with drone node 1, the second $|C|$ members are associated with the second drone node, and so on and so forth. Therefore, we have $|K|=|C|\times|D|$.
\begin{table}[!h]
\protect\caption{Sets, parameters, and variables used in the one center model.}
\label{tab:Sets-Parameters-Variables-ProblemII-Alternative}
\begin{center}
\begin{tabular}{c|p{0.72\textwidth}}  
      \hline
      \hline
      \textbf{Sets/Parameters/Variables} & \textbf{Description} \\
      \hline
      $C$ & Customer Nodes \\
      $D$ & Drone Nodes\\
      $K$ & Set of artificial drones representing trips for drones in $D$ \\
      $T$ & Truck nodes, i.e., only one node here. We call it node $0$. \\
      \hline
      $d_{ij}$ & distance between node $i$ and $j$ for $i,j\in C\cup D \cup T$ \\
      $v_D$ & Speed of drones\\
      $L$ & Longest distance a drone can travel without charging battery \\
      \hline
       $x_{ijk}$ & Binary variable. It is $1$ when drone $k$ travels from node $i$ to node $j$ and $0$ otherwise. \\
       $Z_{ijk}$ & Binary variable. It is $1$ if both $x_{i0k},i\in D$ and $x_{0jk},j\in C$ are 1 and 0 otherwise.\\
       $Y_{ijk}$ & Binary variable. It is $1$ if both $x_{i0k},i\in D$ and $x_{j0k},j\in C$ are 1 and 0 otherwise.\\
       $q_{i}$ & Total travel time of the drone with base at node $i$\\
       $Q$ & Maximum time spent by all drones (makespan) \\
       \hline
   \end{tabular}
\end{center}
\end{table}

Then, an alternative mixed integer linear programming (MILP) model for this problem can be written as follows: 
  \begin{eqnarray}
     \minimize \qquad  Q & &  \quad \qquad \st \nonumber  \\
    & & \nonumber \\
    \sum_{i\in D\cup \{0\}}\, \sum_{k\in K} x_{jik} & = & 1\,,   \qquad \forall j\in C  \label{eq:Problem2AlternativeDroneSatisfyAll} \\
     \sum_{k\in K} x_{0ik} + \sum_{j\in C}  \sum_{k\in K} x_{ijk} & = & 0 \,, \qquad \forall i \in D      \label{eq:Problem2AlternativeNoSelfNoDirect} \\
    \sum_{k\in K} x_{0jk} & = & 1 \,, \qquad \forall j \in C       \label{eq:Problem2AlternativeCenterSatisfyAll-1} \\
    \sum_{i\in D\cup C}x_{i0k} - \sum_{j \in C} x_{0jk} & = & 0\,, \qquad \forall k\in K  \label{eq:Problem2AlternativeCenterSatisfyAll-2} \\
     \sum_{k\in K} x_{i0k} - \sum_{j \in C}  \sum_{k\in K} x_{jik} & = & 0\,, \qquad \forall i\in D \label{eq:Problem2AlternativeDroneNodeFlowBalance} \\
    \sum_{i \in D\cup C}  \sum_{k\in K} x_{i0k} - \sum_{j\in D\cup C}  \sum_{k\in K} x_{0jk} & = & 0\,, \qquad   \label{eq:Problem2AlternativeCenterFlowBalance} \\
   \sum_{i \in D} d_{i0}x_{i0k} + \sum_{j\in C}(d_{0j}x_{0jk} + d_{j0}x_{j0k}) + \sum_{j\in C} \sum_{i \in D} d_{ji} x_{jik} & \leq & L\,, \qquad \forall k \in K \label{eq:Problem2AlternativeDroneRange} \\
    \sum_{k\in K}d_{i0}x_{i0k} + \sum_{j\in C}\sum_{k\in K} (d_{0j}x_{i0k}x_{0jk} + d_{j0}x_{i0k}x_{j0k}) + \sum_{j\in C}\sum_{k\in K} d_{ji}x_{jik} & \leq & v_D \, q_i\,, \qquad \forall i\in D \label{eq:Problem2AlternativeDroneTime}\\
    q_i & \leq & Q\,,    \qquad \forall i\in D   \label{eq:Problem2AlternativeTotalTime} \\
        x_{ijk} & \in & \{0,1\} \,, \; \forall i,j\in C\cup D \cup \{0\}\,, \, \forall k\in K \;\; \nonumber \\
    q_i\,, \; Q & \geq & 0 \,, \qquad \forall i\in D \nonumber 
  \end{eqnarray} 
Constraint \eqref{eq:Problem2AlternativeDroneSatisfyAll} ensures all customers are visited once. 
Constraints \eqref{eq:Problem2AlternativeNoSelfNoDirect}, \eqref{eq:Problem2AlternativeCenterSatisfyAll-1}, and \eqref{eq:Problem2AlternativeCenterSatisfyAll-2} ensure that drones visit the center (truck node) first before visiting a customer node. 
Constraint \eqref{eq:Problem2AlternativeDroneNodeFlowBalance} and \eqref{eq:Problem2AlternativeCenterFlowBalance} enforce the balance of outbound and inbound flows at each drone node and at the pickup center, respectively.
Constraint \eqref{eq:Problem2AlternativeDroneRange} ensures that the drone range limitation is not violated before the drone goes back to its base.
Constraints \eqref{eq:Problem2AlternativeDroneTime} sums the total time (distance) of all routes starting from each drone node. Constraint \eqref{eq:Problem2AlternativeTotalTime} finds the maximum time (distance) among all drone nodes.

To linearize the $x_{i0k}x_{0jk}$ term, we add one more dummy variable $Z_{ijk}$. Meanwhile, to linearize the $x_{i0k}x_{j0k}$ term, we add one more dummy variable $Y_{ijk}$. We replace \eqref{eq:Problem2AlternativeDroneTime} with \eqref{eq:Problem2AlternativeDroneTimeLinearized}
\begin{eqnarray}
    \sum_{k\in K}d_{i0}x_{i0k} + \sum_{j\in C}\sum_{k\in K} (d_{0j}Z_{ijk} + d_{j0}Y_{ijk}) + \sum_{j\in C}\sum_{k\in K} d_{ji}x_{jik} & \leq & v_D \, q_i\,, \qquad \forall i\in D \qquad \qquad  \quad \label{eq:Problem2AlternativeDroneTimeLinearized}\\
 \end{eqnarray}
and add six more constraints as:
\begin{eqnarray}
  \hspace{2.7in}  Z_{ijk} +1 & \geq & x_{i0k} + x_{0jk} \qquad \forall i\in D, \; j\in C,\; k \in K,  \label{eq:Problem2AlternativeLinearization1} \\
    Z_{ijk} & \leq & x_{i0k} \qquad \qquad \forall i\in D,\;k \in K	\label{eq:Problem2AlternativeLinearization2} \\
    Z_{ijk} & \leq & x_{0jk} \qquad \qquad \forall  j\in C,\; k \in K	\label{eq:Problem2AlternativeLinearization3} \\
    Y_{ijk} +1 & \geq & x_{i0k} + x_{j0k} \qquad \forall  i\in D, \;j\in C,\;k \in K	\label{eq:Problem2AlternativeLinearization4} \\
    Y_{ijk} & \leq & x_{i0k} \qquad \qquad \forall  i\in D,\;k \in K	\label{eq:Problem2AlternativeLinearization5} \\
    Y_{ijk} & \leq & x_{j0k} \qquad \qquad \forall j\in C,\; k \in K, 	\label{eq:Problem2AlternativeLinearization6} 
\end{eqnarray}

\section{Lloyd's Algorithm}
\label{sec:LloydsAlg}
\begin{figure}[h!]
	   \centering
    \begin{subfigure}[b]{0.3\linewidth}
        \includegraphics[width=\linewidth]{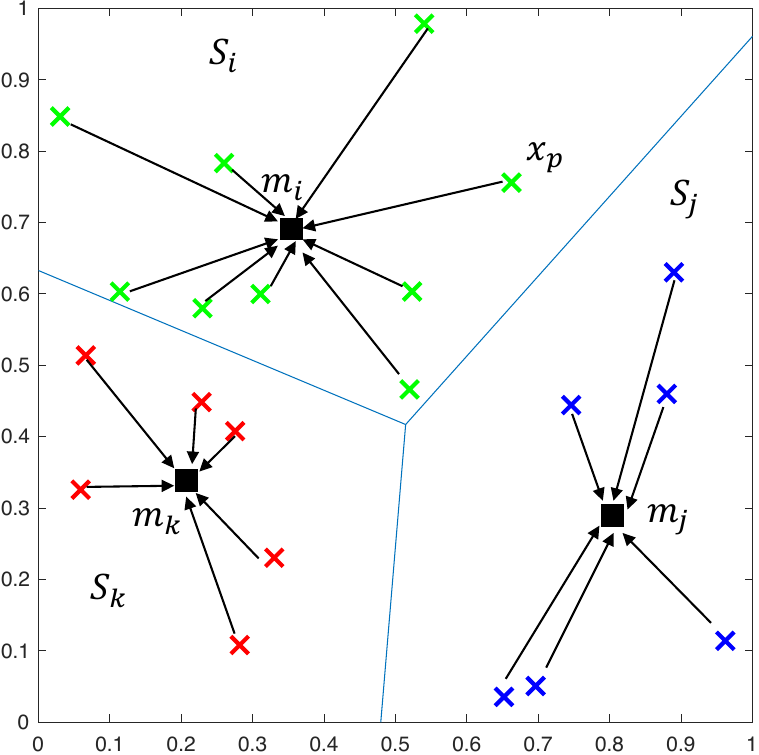}
        \label{fig:kMeansOneIteration}
        \caption{}
    \end{subfigure}
    \qquad
    \begin{subfigure}[b]{0.3\linewidth}
        \includegraphics[width=\linewidth]{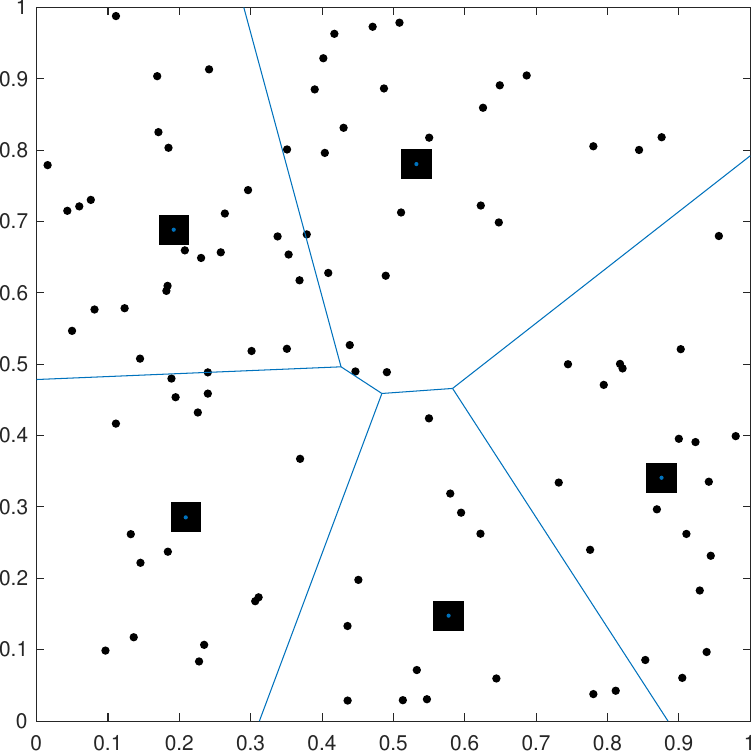}
        \label{fig:kMeansPartitioningExample}
        \caption{}
    \end{subfigure}
    \caption{Assignment of points in each iteration of Lloyd's ($k$-means clustering) algorithm is shown in \subref{fig:kMeansOneIteration}. One random example of $k$-means partitioning with $k=5$ is illustrated in \subref{fig:kMeansPartitioningExample} where the Voronoi cells show the partitions.}
	   \label{fig:kMeansAlg}
\end{figure}

Iteration $t$ of the Lloyd's algorithm works as follows. Given a set of points $P=\{x_1,x_2,\dots,x_n\}\in R$, and an initial set of $k$ means $m_1^{(1)},\dots,m_k^{(1)}$ (which could be set as the same as some of $x_p$'s or be selected randomly in $R$), we want to cluster these points into $k$ sets ${S_1,S_2,\dots,S_k}$ with one mean point in each set such that we assign each point $x_p,\; \forall p$ to the cluster whose mean has the least squared Euclidean distance to $x_p$. This is summarized in \eqref{eq:kMeansClusteringAssignmentStep} and shows in \cref{fig:kMeansAlg}.
\begin{equation}
  S^{(t)}_i = \{x_p \in P: ||x_p-m^{(t)}_i||^2 \leq ||x_p-m^{(t)}_j||^2 \,, \;\;\;  \forall j\,, \;\;\; 1\leq j\leq k \}
  \label{eq:kMeansClusteringAssignmentStep}
\end{equation}
Note that each $x_p$ is assigned to exactly one set $S_i^{(t)}$. The $k$-means clustering algorithm is then followed by an updating step, in which we calculate the new means as the centroids of the points in the each new clusters via \eqref{eq:kMeansClusteringConvergenceStep}.
\begin{equation}
    m_i^{(t+1)} = \frac{1}{|S_i^{(t)}|}\sum_{x_p\in S_i^{(t)}} x_p
    \label{eq:kMeansClusteringConvergenceStep}
\end{equation}
The algorithm stops when $\|m_i^{(t+1)} - m_i^{(t)}\|$ or $\max_p \|x_p^{(t+1)} - x_p^{(t)}\|$ is less than a preset threshold. It is known that if the algorithm converges it converges approximately to the \emph{centroidal Voronoi diagram}.
	     
    \begin{figure}[h!]
    \centering
    \begin{subfigure}[b]{0.24\linewidth}
        \includegraphics[width=\linewidth]{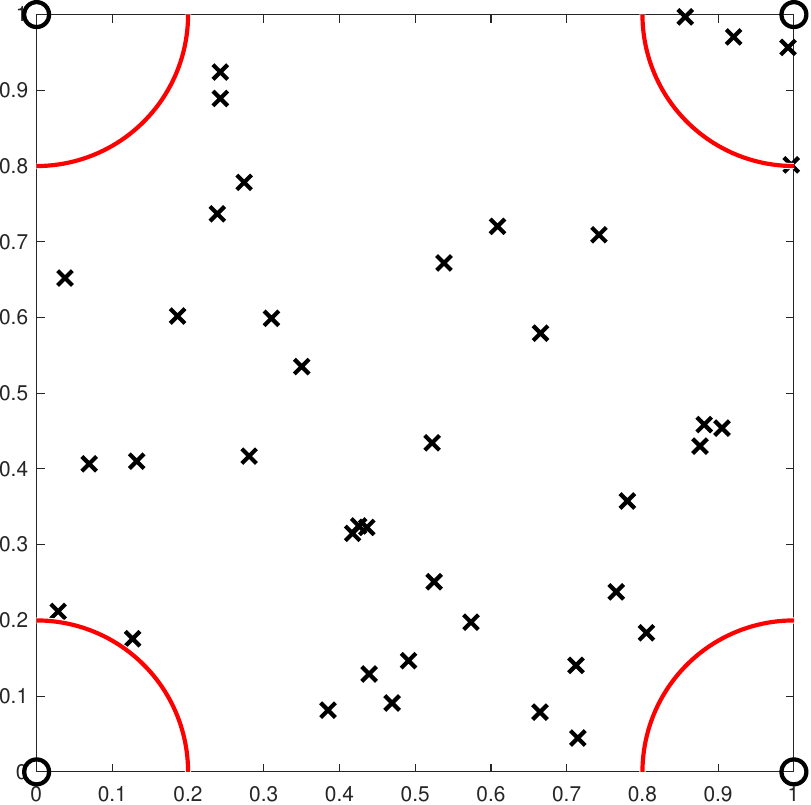}
        \caption{Initial Clusters}
    \end{subfigure}
    \begin{subfigure}[b]{0.24\linewidth}
        \includegraphics[width=\linewidth]{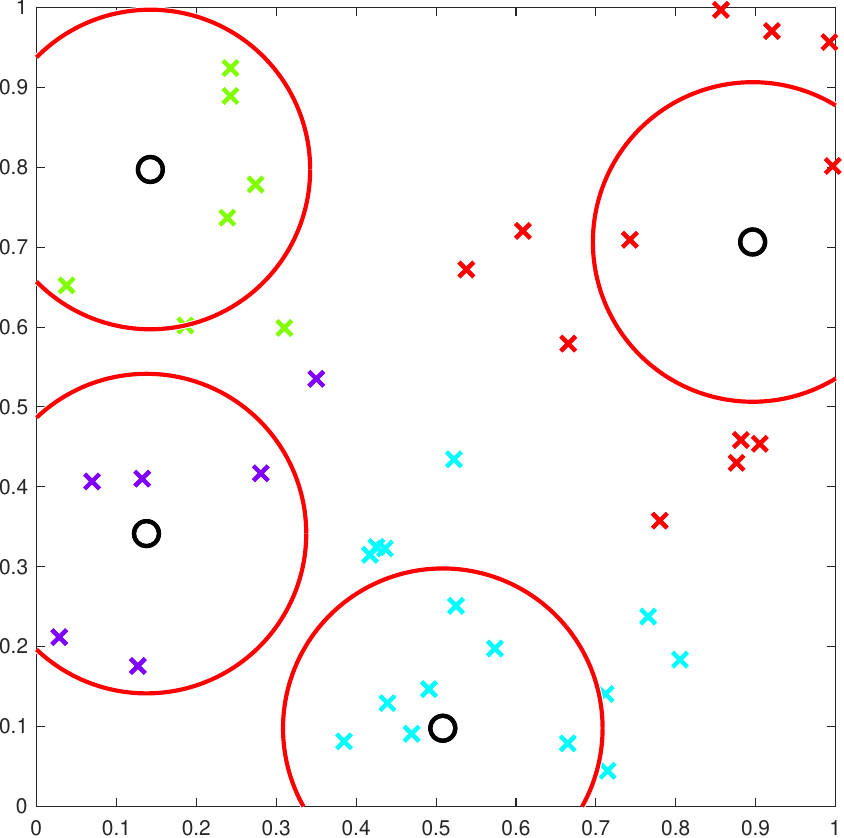}
        \caption{Step One of Cluster}
    \end{subfigure}
    \begin{subfigure}[b]{0.24\linewidth}
        \includegraphics[width=\linewidth]{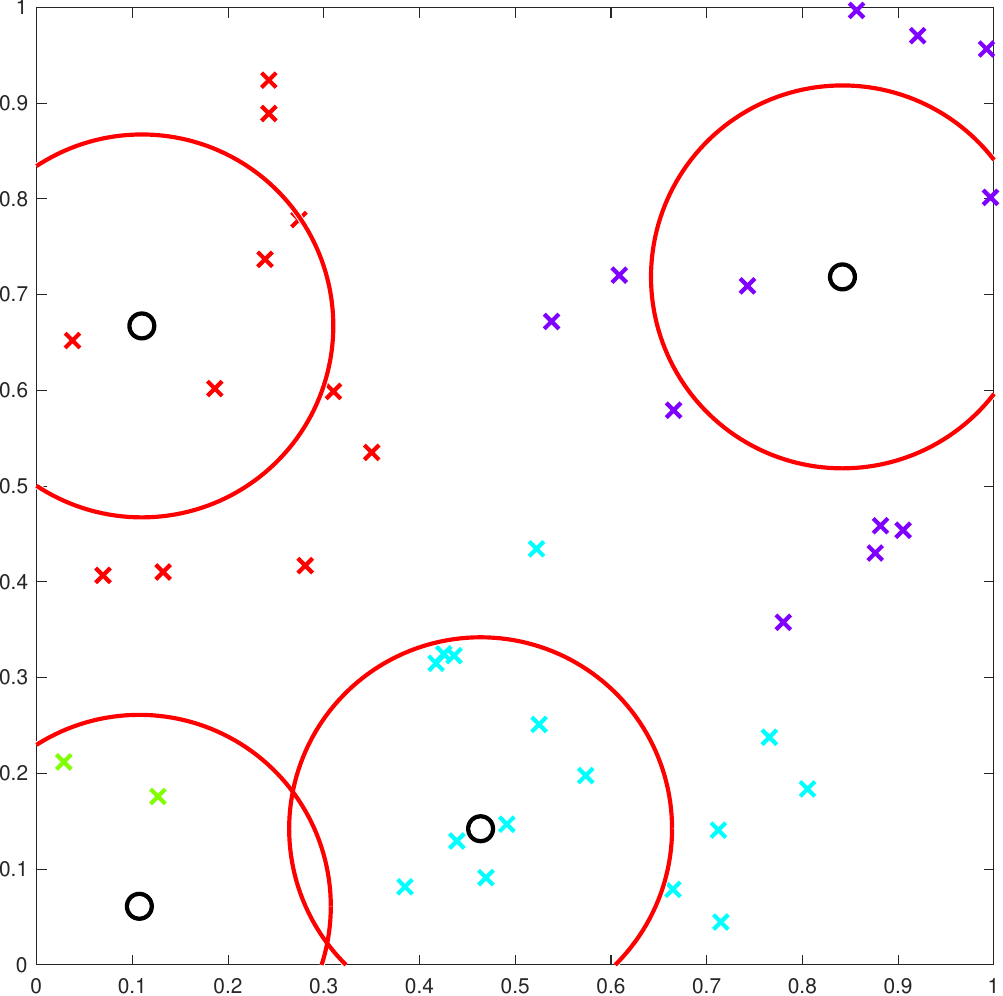}
        \caption{Step Two of Cluster}
    \end{subfigure}
    \begin{subfigure}[b]{0.24\linewidth}
        \includegraphics[width=\linewidth]{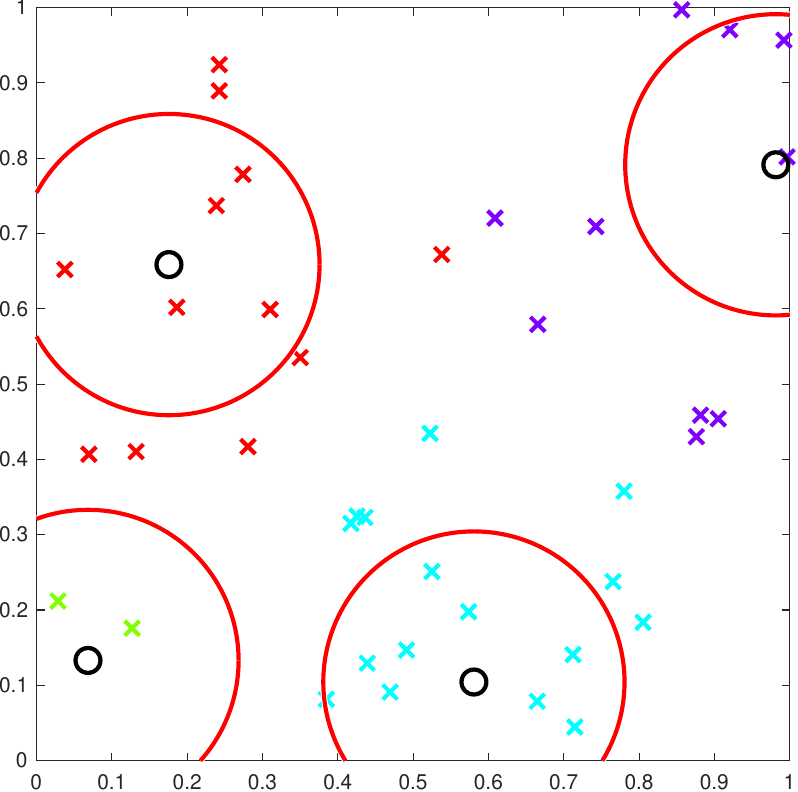}
        \caption{Final Step of Cluster}
    \end{subfigure}
    \caption{Iterations of Lloyd's ($k$-means clustering) algorithm with $k=4$.}
    \label{fig:kmeans}
\end{figure}

\section{Tabu Search Algorithm}
\label{sec:TSAlg}
Tabu Search Algorithm is a meta-heuristic method designed to overcome local optimality in heuristic search procedures. It uses adaptive memory functions to create a flexible searching move between new solutions and past solutions. The three main strategies of Tabu Search are the forbidding strategy, freeing strategy, and short-term strategy. 
	  
In forbidding strategy if an action has been made within a short-term period, it will not be permitted to be reinstated and marked as tabu in the tabu list. Therefore, the algorithm will not consider the possible results repeatedly and will not fall into a locally optimal solution. The potential action will exit the tabu list after a certain period of time, which is called the freeing strategy. Effective tabu search algorithms manage the right balance between the forbidding strategy and freeing strategy to search the areas that have not been visited yet before to avoid local optimal and get as close possible to the optimal solution.
	  
The Tabu Search algorithm works according to the following steps:
 \begin{description}
   \item [Step 1:] Generate an initial random solution and evaluate its cost.
   \item [Step 2:] Form a list of actions that defines the strategies to move from one solution to another (Action List).
   \item [Step 3:] Set the counter of these actions to zero (Tabu Counter).
   \item [Step 4:] Implement the actions in the tabu list that the counter is zero and evaluate its solution and the associated cost.
   \item [Step 5:] Choose the best solution among these evaluated solutions and add its action into Tabu List (set its counter to a certain number). Decrease the counter of other actions by one if they are greater than zero.
   \item [Step 6:] If the stopping criteria are met, stop the algorithm, otherwise, go back to Step 4.
 \end{description}

\end{document}